%     Ce ficher est du plain TeX.   
%     tvd 1990, edition 2005, fichier contenant le tout 
%     Fran\c{c}ois Dubois,  mars,  mai , 15 aout, 4 septembre  2005. 
%                          Pour hal, 15, 16 juin 2010  
%%%%%%%%%%%%%%%%%%%%%%%%%%%%%%%%%%%%%%%%%%%%%%%%%%%%%%%%%%%%%%%%%%%%%%%%%%%%%%%%%%
   \input epsfx.tex
%%%%%%%%%%%%%%%%%%%%%%%%%%%%%%%%%%%%%%%%%%%%%%%%%%%%%%%%%%%%%%%%%%%%%%%%%%%%%%%%%%
  
\overfullrule=0pt
% \nopagenumbers
 
%taille d'agrandissement  :
\magnification=1400
\hsize=12.5cm
\vsize=16.5cm
% marge gauche
% \hoffset=-0.6cm   % bon reglage pour le postscript        mai 2010
% \voffset=0.5cm    % bon reglage pour le postscript 
\hoffset=-0.3cm   % bon reglage pour le pdf 
\voffset=1.0cm    % bon reglage pour le pdf                 mai 2010

\font \nom=cmcsc10 at 10.5 pt

% \font \pbf=cmb10   scaled 200

\font \gcaps=cmbx10 scaled 1333
 
\font \pecaps=cmcsc10 at 9 pt

%  Pagination d'apres Raymond S\'eroul pages 231 et 68 
\newtoks \hautpagegauche  \hautpagegauche={\hfil}
\newtoks \hautpagedroite  \hautpagedroite={\hfil}
\newtoks \titregauche     \titregauche={\hfil}
\newtoks \titredroite     \titredroite={\hfil}
\newif \iftoppage         \toppagefalse   
\newif \ifbotpage         \botpagefalse    
\titregauche={\pecaps    Fran\c{c}ois Dubois }
\titredroite={\pecaps   Nonlinear Interpolation and Total Variation Diminishing Schemes }
\hautpagegauche = { \hfill \the \titregauche  \hfill  }
\hautpagedroite = { \hfill \the \titredroite  \hfill  }
\headline={ \vbox  { \line {  
\iftoppage    \ifodd  \pageno \the \hautpagedroite  \else \the
\hautpagegauche \fi \fi }     \bigskip  \bigskip  }}
\footline={ \vbox  {   \bigskip  \bigskip \line {  \ifbotpage  
\hfil {\oldstyle \folio} \hfil  \fi }}}

%  petite boullette
\def \smb {{\scriptstyle \bullet }}

% Les nombres mathematiques

\def\R{{\rm I}\! {\rm R}}
\def\Z{{\rm Z}\!\! {\rm Z}}

%  indice bas
\def\ib#1{_{_{\scriptstyle{#1}}}}

\def\mod#1{\setbox1=\hbox{\kern 3pt{#1}\kern 3pt}%
\dimen1=\ht1 \advance\dimen1 by 0.1pt \dimen2=\dp1 \advance\dimen2 by 0.1pt
\setbox1=\hbox{\vrule height\dimen1 depth\dimen2\box1\vrule}%
\advance\dimen1 by .1pt \ht1=\dimen1
\advance \dimen2 by .01pt \dp1=\dimen2 \box1 \relax}

\def\nor#1{\setbox1=\hbox{\kern 3pt{#1}\kern 3pt}%
\dimen1=\ht1 \advance\dimen1 by 0.1pt \dimen2=\dp1 \advance\dimen2 by 0.1pt
\setbox1=\hbox{\kern 1pt  \vrule \kern 2pt \vrule height\dimen1 depth\dimen2\box1
\vrule
\kern 2pt \vrule \kern 1pt  }%
\advance\dimen1 by .1pt \ht1=\dimen1
\advance \dimen2 by .01pt \dp1=\dimen2 \box1 \relax}

%carre de fin de demonstration
\def\sqr#1#2{{\vcenter{\vbox{\hrule height.#2pt \hbox{\vrule width .#2pt height#1pt 
\kern#1pt \vrule width.#2pt} \hrule height.#2pt}}}}
\def\square{\mathchoice\sqr64\sqr64\sqr{4.2}3\sqr33} 
 
%  frise 

%  guillemets francais
%   on fait ce qu on peut puisque les bons symboles francais 
%   n ont pas l air d etre dans les fontes.

% saut de ligne
\def\br {\break}

%  page
\def \page #1{\unskip\leaders\hbox to 1.3 mm {\hss.\hss}\hfill {  $\,\,$ {\oldstyle #1}}}

%%%%%%%%%%%%%%%%%%%%%%%%%%%%%%%%%%%%%%%%%%%%%%%%%%%%%%%%%%%%%%%%%%%%%%%%%%%%%%%%%%%%%%  

$~$

% \bigskip
\centerline{\gcaps Nonlinear Interpolation and }
\bigskip
\centerline{\gcaps  Total Variation Diminishing Schemes  $\,$\footnote 
{$ ^{\square}$}{\rm Rapport de recherche Aerospatiale Espace et Defense, 
n$^{\rm o}$ ST/S 46 195,\br
 6 juin 1990, reproduit  au chapitre 4.4 de la th\`ese d'habilitation de l'auteur
 (16~d\'e\-cem\-bre 1992). Version abr\'eg\'ee publi\'ee dans les
 {\it Proceedings of the 3rd International Conference on Hyperbolic Problems}, 
 Bj\"orn Engquist and  Bertil Gustafsson Editors, 
Chartwell-Bratt, volume 1, pages 351-359, 1991. 
Edition \TeX  ~ du 15 ao\^ut  2005, mise en pages du 19 juin 2010. }}
\bigskip 
\bigskip

\centerline { \nom  Fran\c{c}ois Dubois  $\,$\footnote 
{$ ^{\displaystyle \ast}$}{\rm Conservatoire National des Arts et
M\'etiers (Paris) et Universit\'e Paris Sud (Orsay), 
francois.dubois@math.univ-psud.fr.}}

\bigskip 
\bigskip 
\noindent {\bf Abstract}  

\noindent
	The Van Leer approach for the approximation of nonlinear scalar conservation 
laws is studied in one space dimension. The problem can be reduced to a nonlinear
interpolation and we propose a convexity property for the interpolated values. We prove
that under general hypotheses the method of lines in well posed in  $\, \ell^{\infty} \,
\cap \, {\rm BV} \,$  and we give precise sufficient conditions to establish that the total
variation is diminishing. We observe that the second order accuracy can be maintained even
at non sonic extrema. We establish also that both the TVD property and second order
accuracy can be maintained after discretization in time with the second order accurate Heun
scheme. Numerical illustration for the advection equation is presented.

\bigskip   
\bigskip 
\centerline {\bf Interpolation non lin\'eaire}
% \smallskip 
\centerline {\bf   et  sch\'emas \`a variation totale d\'ecroissante}
\bigskip 
\noindent {\bf R\'esum\'e}  

\noindent 
Nous \'etudions le sch\'ema de Van Leer pour l'approximation de lois de conservation
scalaires \`a une dimension d'espace. Nous proposons une propri\'et\'e de convexit\'e de
l'interpolation non lin\'eaire associ\'ee \`a cette m\'ethode. Nous prouvons que sous des
hypoth\`eses g\'en\'erales, le sch\'ema continu en temps conduit \`a un probl\`eme bien
pos\'e et \`a variation totale d\'ecroissante. Ce r\'esultat s'\'etend au cas
discr\'etis\'e en temps et des tests num\'eriques mettant en \'evidence divers ordres de
convergence sont pr\'esent\'es.



%  \vfill \eject 
 
\bigskip 
\noindent {\bf Contents }

\noindent  1) \quad  Introduction  \page {$\,\,$2}

\noindent  2) \quad Some Limitations for the Interface Values \page {$\,\,$5}

\noindent  3) \quad Decrease of the Total Variation  \page {$\,\,$15}

\noindent   4) \quad Convergence, Order of Accuracy \page {$\,\,$23}

\noindent  5) \quad Time Discretization \page {$\,\,$29}

\noindent  6) 	\quad Numerical Tests with the Advection Equation \page {$\,\,$39}

\noindent  7)  \quad Conclusion \page {$\,\,$50}

\noindent  8) \quad References \page {$\,\,$50}

% \titredroite={\pecaps  Introduction  }
\toppagetrue  
\botpagetrue    

\bigskip  
\noindent {\bf 1. \quad Introduction}

\smallskip \noindent
We consider the following scalar conservation law :

\smallskip \noindent   (1.1)  $\qquad   \displaystyle 
{{\partial u}\over{\partial t}} \,+\, {{\partial}\over{\partial x}} f(u) \,\,=\,\, 0
\,\,,\qquad x \in \R \,,\quad t > 0 \,,\quad u(x,\,t) \in \R \, \,$

\smallskip \noindent 
where the flux function $\,f \,$  is supposed to be a convex regular (${\cal C}^2$ class)
real function. The initial condition takes the form

\smallskip \noindent   (1.2)  $\qquad   \displaystyle 
u(x,\,0) \,\,=\,\, u_0(x) \,,\qquad x \in \R \,$

\smallskip \noindent 
and we assume that the initial datum $\, u_0 \,$ satisfies :

\smallskip \noindent   (1.3)  $\qquad   \displaystyle 
u_0 \in {\rm L}^{\infty}(\R) \, \cap \, {\rm BV}(\R) \,.\, $

\bigskip \noindent $ \bullet \qquad $ 
Even for smooth initial datum $\, u_0 ,\,$ a classical solution  of the Cauchy problem
(1.1).(1.2) does not exist in general ({\it e.g.} Smoller [1983]) and we consider in this
paper weak solutions of the problem (1.1).(1.2), {\it i.e.} functions 
$\, u \in {\rm L}^{\infty}(\R \times [0,\, +\infty[) \,$   such that

\smallskip \noindent   (1.4)  $\qquad   \displaystyle 
\int_{\R \times \R_+} \, \Big[ \, u \, {{\partial \varphi}\over{\partial t}} \,+\, f(u) \, 
 {{\partial \varphi}\over{\partial x}} \, \Big] \,  {\rm d}x \, {\rm d}t \,\,=\,\, 0 \,,\qquad
\forall \varphi \in {\cal C}_0^1(\R \times   ]0,\, +\infty[) \,$ 

\smallskip \noindent  
and 

\smallskip \noindent   (1.5)  $\qquad   \displaystyle 
\lim_{\displaystyle  t \longrightarrow 0} \, 
\mid \mid    u(\smb,\,t) \,-\, u_0(\smb)  \mid \mid\ib{{\rm L}^1(\R)}  
 \,\,=\,\, 0 \,.\,$

\bigskip \noindent $ \bullet \qquad $ 
The uniqueness of a weak solution satisfying (1.4).(1.5) is guaranteed
 if we consider (Lax [1971]) only entropy solutions of the conservation law (1.1), {\it
i.e.} weak solutions also satisfying 

\smallskip \noindent   (1.6)  $\quad   \displaystyle 
\int_{\R \times \R_+}  \!\!\! \Big[ \, \eta(u)  \, {{\partial \varphi}\over{\partial t}} \,+\,
\xi(u)  \, {{\partial \varphi}\over{\partial x}} \,\Big] {\rm d}x \, {\rm d}t \, \leq\,  0
\,,\quad \forall \varphi \in {\cal C}_0^1(\R \times   ]0,\, +\infty[) \,,\,  \varphi \geq
0 \,  $ 

\smallskip \noindent 
where  $\, \eta \,:\, \R \longmapsto \R \,$   is a strictly convex entropy and  $\,
\xi\,:\, \R \longmapsto \R \,$  is the associated entropy flux :

\smallskip \noindent   (1.7)  $\qquad   \displaystyle 
\xi'(u) \,\,=\,\, \eta'(u) \, \smb \, f'(u) \,.\,$

\smallskip \noindent 
With the particular choice 

\smallskip \noindent   (1.8)  $\qquad   \displaystyle 
\eta_c(u) \,\,=\,\, \mid u-c \mid \,,\qquad c \in \R \,, \,$

\smallskip \noindent 
Kruskov [1970] proved that a weak solution of (1.1).(1.2) satisfying the 
inequalities (1.6) for all the choices of entropies (1.8) is necessarily unique. More
recently, Di Perna [1983] proved that uniqueness is implied if the inequality (1.6) is
valid for a single strictly convex entropy $\, \eta .\,$ 

\bigskip \noindent $ \bullet \qquad $ 
The approximation of (1.4)-(1.6) by finite volume numerical schemes is a very classical
problem. In this paper, we first look at approximations that  are discrete in space and
continuous in time (method of lines). Let  $\, h > 0 \,$   be a real parameter and 
$\, x_j = j \, h \,$   be the associated mesh points. Following Godunov [1959], we search an
approximation  $\, u_j(t) \,$   which is constant in the interval

\smallskip \noindent     $  \displaystyle 
{\rm I}_j \, \equiv \, \Big] \, \big( j-{{1}\over{2}}\big) \,h \,,\,
\big(j+{{1}\over{2}}\big) \,h \, \Big[ \,$ : 

\smallskip \noindent   (1.9)  $\qquad   \displaystyle 
u_j(t) \,\,\approx\,\, {{1}\over{h}} \, \int_{ ( j-{{1}\over{2}} )
\,h}^{ ( j+{{1}\over{2}} ) \,h} u(x,\,t) \, {\rm d}x \,,\qquad j \in
\Z \,,\quad t> 0 \,.\, \,$

\smallskip \noindent   
As was recognized by Lax-Wendroff [1960], 
we require $\,u_j(t) \,$ to satisfy a conservative numerical scheme :

\smallskip \noindent   (1.10)  $\qquad   \displaystyle 
{{{\rm d}u_j}\over{{\rm d}t}} \,+\,  {{1}\over{h}} \, \Big( f_{j+{{1}\over{2}}} \,-\, 
 f_{j-{{1}\over{2}}} \Big) \,\,=\,\, 0 \,,\qquad  j \in \Z \,,\quad t> 0 \,.\, \,$

\bigskip \noindent $ \bullet \qquad $ 
The system of ordinary differential equations (1.10)  is now completely defined if we
prescribe the initial conditions and the numerical flux $\, f_{j+{{1}\over{2}}}.\,$ 
 Three point numerical schemes are parameterized by the numerical flux function $\, \Phi
\, $~:
\smallskip \noindent   (1.11)  $\qquad   \displaystyle 
f_{j+{{1}\over{2}}} \,\, =\,\, \Phi(u_j \,,\, u_{j+1} ) \,$ 

\smallskip \noindent 
which is supposed to be Lipschitz continuous and consistent with the flux $\, f\,$~:

\smallskip \noindent   (1.12)  $\qquad   \displaystyle 
\Phi(u \,,\, u  ) \,\,\equiv\,\, f(u) \,,\qquad u \in \R \,.\,$ 

\smallskip \noindent  
The numerical flux realizes some approximation of the Godunov flux (we refer to
Harten-Lax-Van Leer [1983] for the details). Following Crandall-Majda [1980], the flux $f$ of
relation (1.11) is monotone if $\, (u,\,v) \longmapsto \Phi(u,\,v) \,$ 
 is a nondecreasing (respectively
nonincreasing) function of $u$ (resp. $v$).  Kusnezov-Volosin [1976] and
Crandall-Majda [1980]  have extended to  monotone schemes the proof of convergence towards
the entropy solution developed by   Leroux [1976] for the Godunov scheme (when the equation
(1.10) is also discretized in time). Sanders [1983] proved that on irregular meshes, the
solutions of the ODE (1.10) also converges to the entropy solution of (1.4)-(1.6) as $\,h
\,$ tends to zero. Osher [1984] with the ``E-schemes'' extended the notion of monotone schemes
to take into account non necessarily convex fluxes. The major default of three point
monotone schemes is that they can have at best first order accuracy.

\bigskip \noindent $ \bullet \qquad $ 
 High resolution schemes (that are in general five point numerical schemes) have been proposed
to overcome this lack of precision. On one hand, Harten [1983] introduced the notion of total
variation diminishing (TVD) schemes ; if the discrete total variation 

\smallskip \noindent   (1.13)  $\qquad   \displaystyle 
{\rm TV}(t) \,\, \equiv \,\, \sum_{  j \in \Z} \, \mid u_{j+1}(t) - u_j(t) \,
\mid \,$ 

\smallskip \noindent 
is a nonincreasing function of time, then the scheme (1.10) is said to
be TVD. Monotone schemes are particular cases of TVD schemes and Harten gave sufficient
conditions (generalized by Sanders [1983] for the method of lines) to ensure that a
numerical scheme is TVD. Unfortunately, three point TVD schemes are at most first order
accurate and Harten [1983] constructed second order numerical schemes by modifying the
first order flux (1.11) into a five point scheme in a way analogous to the Flux Corrected
Transport  (FCT) of Boris-Book [1973].

\bigskip \noindent $ \bullet \qquad $ 
On the other hand, the ``Multidimensional Upstream centered Scheme for Conservation Laws''
(MUSCL) approach proposed by Van Leer [1979] (see also Colella [1985]) assumes that a
piecewise linear interpolation is reconstructed from the mean values$\, u_j(t) \,$ : 

\smallskip \noindent   (1.14)  $\qquad   \displaystyle 
u_h(x,\,t) \,\,=\,\, u_j(t) \,+\, s_j(t) \, (x- x_j) \,,\qquad x \in {\rm I}_j \,,\quad j
\in \Z \,.\,$ 

\smallskip \noindent  
The slopes $\, s_j  \,$  are chosen as nonlinear functions of 
the three mean values at the points $\,  j-1 , \, j ,\,  j+1 \,$ ~:

\smallskip \noindent   (1.15)  $\qquad   \displaystyle 
s_j \,\,=\,\, S ( u_{j-1} \,,\, u_j \,,\, u_{j+1} ) \,$ 

\smallskip \noindent 
in such a way that the reconstruction (1.14) satisfies some monotonicity 
restrictions (Van Leer [1977], see also Sweby [1984] for a unified vision of the flux
limiters). The function $\, x \longmapsto u_h(x,\,t) \,$    is {\it a priori} discontinuous
at the points $\, x_{j+{{1}\over{2}}}  \,$~;  the two values $\, u_{j+{{1}\over{2}}}^{\pm}
\,$  are defined on each side of $\, x_{j+{{1}\over{2}}}  \,$~:

\smallskip \noindent   (1.16)(a)  $\qquad   \displaystyle 
 u_{j+{{1}\over{2}}}^{-} \,\,=\,\, u_j(t) \,+\, {{h}\over{2}} \, S (  u_{j-1} \,,\, u_j
\,,\, u_{j+1} ) \,$

\smallskip \noindent   (1.16)(b)  $\qquad   \displaystyle 
 u_{j+{{1}\over{2}}}^{+} \,\,=\,\, u_{j+1}(t) \,-\, {{h}\over{2}} \, S (   u_j \,,\,
u_{j+1}\,,\, u_{j+2}  ) \,. \,$

\smallskip \noindent   
The flux function $\, f_{j+{{1}\over{2}}} \,$   is computed as in the Godunov scheme, but
with the latter two values as arguments :

\smallskip \noindent   (1.17)  $\qquad   \displaystyle 
f_{j+{{1}\over{2}}} \,\,=\,\, \Phi \big(  u_{j+{{1}\over{2}}}^{-} \,,\, 
u_{j+{{1}\over{2}}}^{+} \big) \,,\qquad t \geq 0 \,,\quad j \in \Z \,.\,$

\smallskip \noindent  
The theoretical convergence properties of MUSCL schemes have been
studied by Osher [1985]. He gave sufficient conditions to establish that the scheme
(1.10).(1.17) is TVD and proved under some restrictions the convergence to the entropy
solution (using a discrete entropy inequality proposed in Osher [1984]). Moreover second
order accuracy is realized even at regular extrema which are not sonic points,  {\it i.e.} values
of u such that 

\smallskip \noindent   (1.18)  $\qquad   \displaystyle 
f'(u) \,\, \neq \,\, 0 \,.\,$ 

\smallskip \noindent  
This fact was recognized by Osher-Chakravarthy [1984]. The MUSCL approach
 has been extended to higher orders by Collela-Woodward [1984] (PPM method) and by Harten
et al [1987] with the so-called ENO schemes.

\bigskip \noindent $ \bullet \qquad $ 
The main purpose of this paper is to show how to choose nonlinear interpolation formulae
like (1.16) in order to get a TVD scheme. In the second part of this paper we propose some
natural limitations on the interpolated values between two mesh points and we establish in
the third part sufficient conditions that prove that the resulting scheme is TVD. In part
IV, we demonstrate that the second order accuracy is obtained in {\bf all} smooth parts of
the flow. In part V we dicretize the method of lines in time and we precise conditions to
maintain both the TVD property and second order accuracy in smooth regions. In the last
part we present some numerical experiments for the advection equation with periodic
boundary conditions.

% \titredroite={\pecaps  Some Limitations for the Interface Values  }
\toppagetrue  
\botpagetrue    

\bigskip   \bigskip  
\noindent {\bf 2. \quad Some Limitations for the Interface Values}

\smallskip \noindent
We recall that a sequence of values $\, u_j \,$   is supposed to be
given at the points  $\, x_j \,$  of a regular mesh on the real axis. We notice a small
difference with the finite volume method where the variables  $\, u_j \,$   represent the
mean values of some function in each interval $\, {\rm I}_j .\,$   We interpolate
those values at the points $\, x _{j+{{1}\over{2}}} \,$   by two different values 
$\,  u_{j+{{1}\over{2}}}^{-} \,$    and $\,  u_{j+{{1}\over{2}}}^{+} \,$    on each side of
the interface  $\, x _{j+{{1}\over{2}}} .\,$ By analogy with the MUSCL approach (1.16), we
set

\smallskip \noindent   (2.1)  $\qquad   \displaystyle 
 u_{j+{{1}\over{2}}}^{-} \,\,=\,\,L ( u_{j-1} \,,\, u_j \,,\, u_{j+1} ) \,$ 

\smallskip \noindent   (2.2)  $\qquad   \displaystyle 
 u_{j+{{1}\over{2}}}^{+} \,\,=\,\, R ( u_j \,,\, u_{j+1} \,,\, u_{j+2}  ) \,. \,$ 

\smallskip \noindent  
The interpolation defined by (2.1).(2.2) is a priori {\bf nonlinear} 
and in the following of
this part we give some restrictions concerning the functions $\, L(\smb,\,\smb,\,\smb) \,$
and $\, R(\smb,\,\smb,\,\smb) .\,$

\bigskip 
% \vfill \eject         %%%   fd 15 juin 2010 
\noindent  {\bf Property 2.1} \quad   {\bf Homogeneity}  

\noindent
We suppose that the multiplication of all the arguments of (2.1).(2.2) by 
the same real number $\, \lambda \,$  multiplies the interpolate values
$\,  u_{j+{{1}\over{2}}}^{-} ,\,$ $\,  u_{j+{{1}\over{2}}}^{+} \,$   by the same argument,
{\it i.e.} that the functions L and R are {\bf homogeneous of degree 1}~:

\smallskip \noindent   (2.3)  $\qquad   \displaystyle 
J( \lambda u  \,,\, \lambda v  \,,\, \lambda w ) \,\,=\,\, \lambda \, J(u,\,v,\,w) 
\,,\qquad  J=L \,\, {\rm or} \,\, R \,,\quad \lambda \in \R \,.\,$ 

\bigskip \noindent  {\bf Property 2.2} \quad   {\bf Translation invariance}  

\noindent
We suppose that by addition of some constant  $\, \lambda \,$ to the three arguments of
(2.1).(2.2), the resulting interpolate values are translated by the same value :

\smallskip \noindent   (2.4)  $\qquad   \displaystyle 
J(u+ \lambda \,,\, v + \lambda  \,,\, w + \lambda) \,\,=\,\, J(u,\,v,\,w) \,+\, \lambda
\,,\quad  J=L \,\, {\rm or} \,\, R \,,\quad \lambda \in \R \,.\,$ 

\bigskip \noindent  {\bf Property 2.3} \quad   {\bf Left-Right symmetry}  

\noindent
We suppose that by exchange of the two extreme arguments of (2.1) 
({\it i.e.} we reverse left
and right) the left interpolation at  $\, x_{j+{{1}\over{2}}} \,$   is transformed into the
right interpolation at $\, x_{j-{{1}\over{2}}} \,$ (Figure 2.1) :

\smallskip \noindent   (2.5)  $\qquad   \displaystyle 
R ( w ,\, v ,\, u ) \,\, =\,\, L( u ,\, v ,\, w ) \,. \,$

\bigskip 
\centerline { \epsfysize=5cm    \epsfbox  {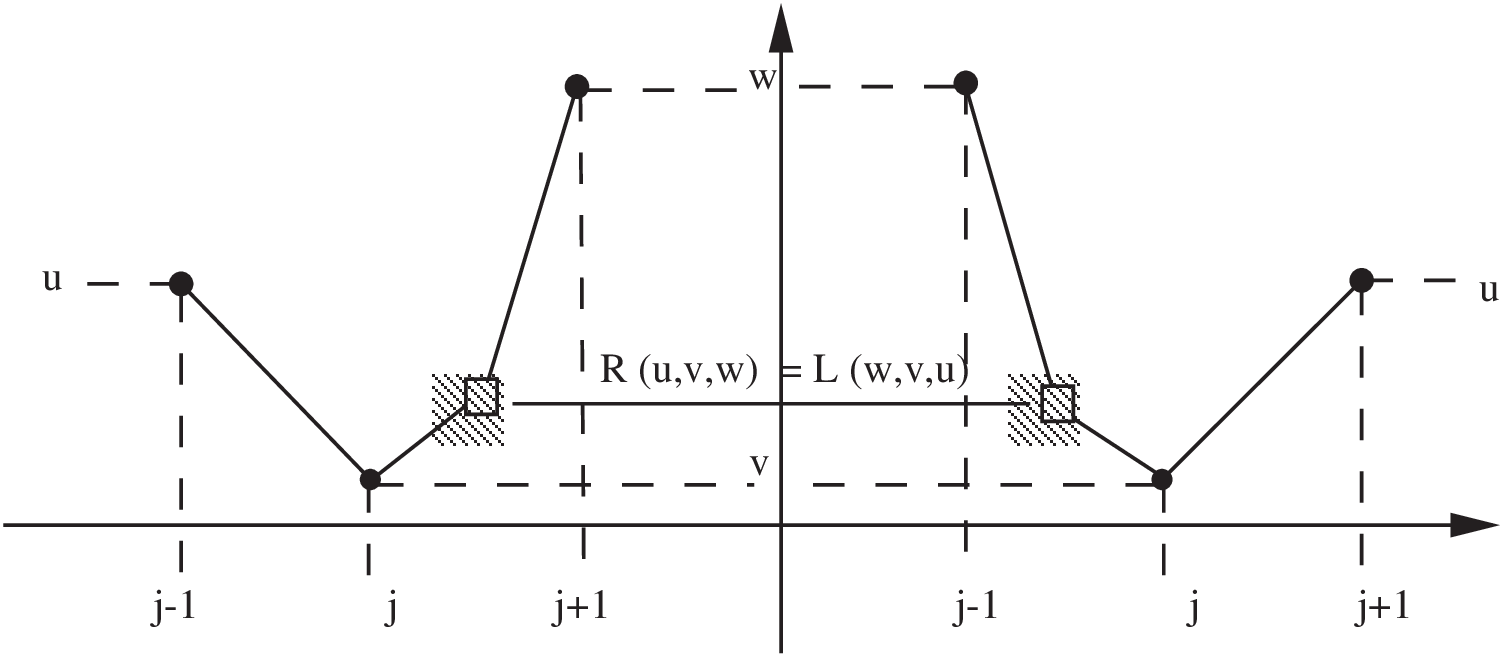} }
\smallskip  \smallskip

\centerline { {\bf Figure 2.1.} \quad  Left-right symmetry.} 

\bigskip \noindent  {\bf Proposition 2.1}  

\noindent
If the properties 2.1 to 2.3 are assumed, there exists some limiter function  $\, \psi
\,:\, \R \longrightarrow \R \,$   such that

\smallskip \noindent   (2.6)  $\qquad   \displaystyle 
L ( u ,\, v ,\, w ) \,\,=\,\, v \,+\, {{1}\over{2}} \, \psi \Big( {{w-v}\over{v-u}} \Big)
\, (v-u) \,$
\smallskip \noindent   (2.7)  $\qquad   \displaystyle 
R ( u ,\, v ,\, w ) \,\,=\,\, v \,-\, {{1}\over{2}} \, \psi \Big( {{v-u}\over{w-v}} \Big)
\, (w-v) \,. \,$

\smallskip \noindent  {\bf Proof of Proposition 2.1}  

\noindent 
We have simply 

\smallskip \noindent  $   \displaystyle 
 L(u ,\, v ,\, w )   \,\,=\,\,  v \,+\, L( u-v ,\,0 ,\, w-v)  \,\,=\,\,
v \,+\, (v-u) \, \, L  \Big( -1 ,\, 0 ,\, {{w-v}\over{v-u}} \, \Big) \,$ 

\smallskip \noindent 
and  $  \qquad  \displaystyle 
  R (u ,\, v ,\, w )   \,\,=\,\, L(w ,\, v ,\, u )  \,\,=\,\,  v\,+\, {{1}\over{2}} \, 
\psi \Big( {{u-v}\over{v-w}} \Big) \,\,   (v-w) \, . \, $ $\hfill \square $

\bigskip \noindent  {\bf Property 2.4} \quad   {\bf Monotonicity}  

\noindent
We suppose that the interpolated values $\, \smash { u_{j+{{1}\over{2}}}^{-} }\,$ and $\, 
u_{j+{{1}\over{2}}}^{+} \,$   belong to the interval $\, 
[\min (u_j  ,\, u_{j+1}) \,,\,  \max (u_{j}  ,\, u_{j+1})], \,$  {\it i.e.}

\smallskip \noindent   (2.8)  $\qquad   \displaystyle 
\min (v,\, w) \,\, \leq \,\, L(u,\, v ,\, w ) \,\, \leq \,\, \max (v,\, w ) \,$ 
\smallskip \noindent   (2.9)  $\qquad   \displaystyle 
\min (u,\,v ) \,\, \leq \,\, R(u,\, v ,\, w ) \,\, \leq \,\, \max (u,\,v ) \,. \,$ 

\smallskip \noindent  
We remark that the inequalities (2.9) are a simple consequence of (2.5) and (2.8).

\bigskip \noindent  {\bf Proposition 2.2}  

\noindent
	If the properties 2.1 to 2.4 are valid, the limiter $\, \psi \,$ 
 satisfies the condition 

\smallskip \noindent   (2.10)  $\qquad   \displaystyle 
0 \,\, \leq \,\, {{\psi(\lambda)}\over{\lambda}} \,\, \leq \,\, 2 \,,\qquad \lambda \in \R
\,$ 

\smallskip \noindent  
and if the properties 2.1 to 2.3 are supposed, the inequalities (2.10) 
imply the property 2.4.

\smallskip \noindent  {\bf Proof of Proposition 2.2}  

\noindent 
It is a straightforward consequence of (2.6), written under the form :

\smallskip \noindent   (2.11)  $\qquad   \displaystyle 
L(u,\, v ,\, w) \,\, = \,\, v   \,+\, {{1}\over{2}} \, {{\psi \Big(
\displaystyle {{w-v}\over{v-u}} \Big)} \over { \displaystyle  {{w-v}\over{v-u}} }} \, 
\, (w-v) \,. \,$  $\hfill \square $ 

\bigskip \noindent $ \bullet \qquad $ 
We introduce now an original notion concerning the nonlinear interpolation at the 
interfaces, first proposed in Dubois [1988].

\bigskip 
\vfill \eject  %%%%%%%%%%    fd 15 juin 2010 
\noindent  {\bf Property 2.5} \quad   {\bf Convexity} 
 
\noindent 
We suppose that for each integer $\, j,\,$  the sequence of the interpolate values 

\smallskip \noindent   (2.12)  $\qquad   \displaystyle 
u_{j-1} \,,\, u_{j-{{1}\over{2}}}^{+} \,,\, u_j \,,\, u_{j+{{1}\over{2}}}^{-} \,,\, u_{j+1}
\,,\,$ 

\smallskip \noindent  
is the restriction to the nodes 

\smallskip \noindent   (2.13)  $\qquad   \displaystyle 
x_{j-1} \,,\, x_{j-{{1}\over{2}}}  \,,\, x_j \,,\, x_{j+{{1}\over{2}}}  \,,\, x_{j+1}
\,,\,$ 

\smallskip \noindent  
of some convex or concave function, depending on the concavity of the sequence

\smallskip \noindent   (2.14)  $\qquad   \displaystyle 
u_{j-1} \,,\,   u_j \,,\,  \,,\, u_{j+1} \,.\,$ 

\smallskip \noindent  
More precisely, if the inequality 

\smallskip \noindent   (2.15)  $\qquad   \displaystyle 
u_j \,-\, u_{j-1} \,\, \leq \,\, u_{j+1} \,-\, u_j \,$ 

\smallskip \noindent  
holds (convexity of the sequence (2.14)), we suppose that the discrete 
increments of (2.12) realizes a nondecreasing sequence :

\smallskip \noindent   (2.16)  $\qquad   \displaystyle  
u_{j-{{1}\over{2}}}^{+} \,-\, u_{j-1} \,\, \leq \,\, u_j \,-\, u_{j-{{1}\over{2}}}^{+}
\,\, \leq \,\, u_{j+{{1}\over{2}}}^{-}\,-\, u_{j} \,\, \leq \,\,  u_{j+1}  \,-\, 
 u_{j+{{1}\over{2}}}^{-} \,$ 

\smallskip \noindent 
and if we change ``$\leq$'' into ``$\geq$'' in the relations (2.15) then we must 
do it also in
(2.16) (concavity of the sequence (2.16)), see Figure 2.2.

\bigskip 
\centerline { \epsfysize=5cm    \epsfbox  {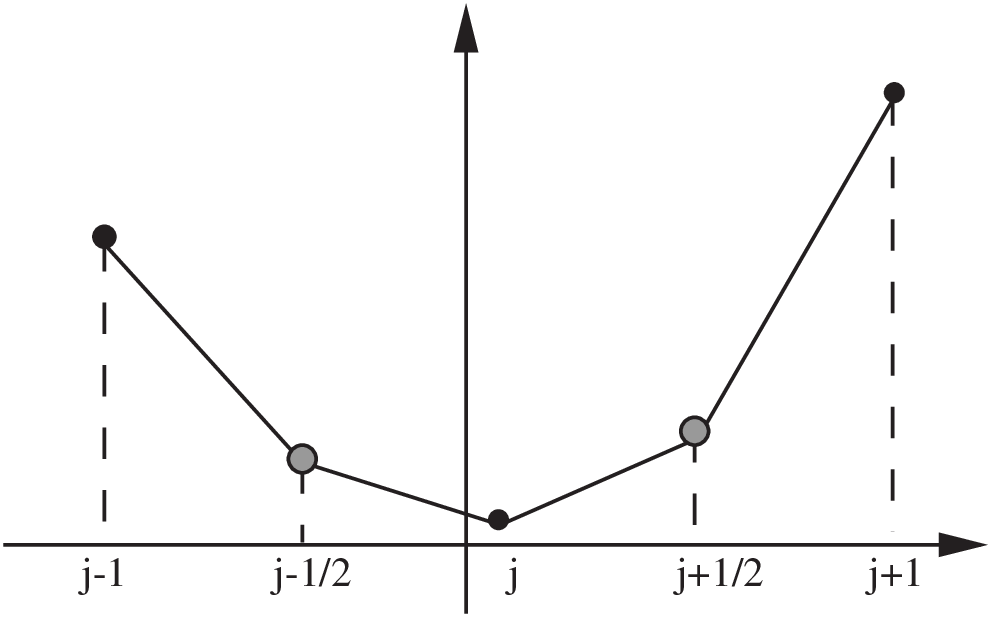} }
\smallskip  \smallskip

% \centerline {
\noindent  {\bf Figure 2.2.} \quad Convexity of the interpolated 
value  at the intermediate points.

\bigskip \noindent  {\bf Example 2.1.}

\noindent 
The classical Lagrange interpolation by a polynomial of degree 2 inside the 
cell $\, ] x_{j-{{1}\over{2}}} ,\,  x_{j+{{1}\over{2}}} [ \,$ corresponds to the limiter
function  $\, \psi(\lambda)  = {{1 + 3 \, \lambda}\over{4}} .\,$    It satisfies
the properties of homogeneity, translation invariance, convexity but not monotonicity 
(see Van
Leer [1973]).

\bigskip 
\vfill \eject  %%%%%%%%%%    fd 15 juin 2010 
\noindent  {\bf Property 2.3}   

\noindent
	If the properties 2.1 to 2.3 hold then the property 2.5 is equivalent to the following 
inequalities concerning the limiter function~:

\smallskip \noindent   (2.17)  $\qquad   \displaystyle  
1 \,\, \leq \,\, \lambda \, \psi \Big( {{1}\over{\lambda}} \Big) \,\, \leq \,\, 
\psi(\lambda ) \,\, \leq \,\,  \lambda \,,\qquad \lambda \geq 1 \,$ 
\smallskip \noindent   (2.18)  $\qquad   \displaystyle  
\lambda\,\, \leq \,\, \psi(\lambda)  \,\, \leq \,\,  \lambda \,\, \psi \Big(
{{1}\over{\lambda}} \Big) \,,\qquad \qquad \,\, \lambda \leq 0 \,. \,$ 

\smallskip \noindent  {\bf Proof of Proposition 2.3}  

\noindent 
We first suppose that the property 2.5 is valid. We set 

\smallskip \noindent   (2.19)  $\qquad   \displaystyle  
w \,\, = \,\, v \,+\, \lambda \, (v - u ) \,.\, $

\smallskip \noindent  
Then  the inequality (2.15) is satisfied if and only if 

\smallskip \noindent   (2.20)  $\qquad   \displaystyle 
(\lambda - 1) \,\, \, (v - u ) \,\, \geq \,\, 0 \,.\,$ 

\smallskip \noindent 
Then the inequalities (2.16) take the form

\setbox21=\hbox {$\displaystyle  
\Big( 1 - {{1}\over{2}} \, \lambda \, \psi(\lambda) \Big) \, (v-u) \,\, \leq \,\, 
 - {{1}\over{2}} \, \lambda \, \psi(\lambda) \, (v-u) \,\, \leq \,\, \dots  $}
\setbox22=\hbox {$\displaystyle  \qquad \qquad \qquad \qquad 
\dots \,\, \leq \,\,    {{1}\over{2}}  \, \psi(\lambda) \, (v-u) \,\, \leq \,\, 
\Big( \lambda \,-\, \psi(\lambda) \Big) \,  (v-u) \, .\,$}
\setbox30= \vbox {\halign{#&# \cr \box21 \cr \box22    \cr   }}
\setbox31= \hbox{ $\vcenter {\box30} $}
\setbox44=\hbox{\noindent  (2.21) $\displaystyle  \qquad  
\left\{ \box31 \right. $}   
\smallskip \noindent   $ \box44 $ 

\smallskip \noindent 
If we suppose that $\, (v-u) \geq 0, \,$  we obtain easily (2.17). If we reverse the
 inequality (2.20), we have also to reverse (2.21), thus (2.17). If $\, 0 \leq \lambda
\leq 1 \,$  the change of variable $\, \mu = {{1}\over{\lambda}} \,$  gives 

\smallskip \noindent   $  \displaystyle 
1 \,\, \geq \,\, {{1}\over{\mu}} \, \psi(\mu) \,\, \geq \,\,  \psi \Big( {{1}\over{\mu}}
\Big) \,\, \geq \,\, {{1}\over{\mu}} \,,\qquad \lambda \,=\, {{1}\over{\mu}} \,,\quad 0 \,<
\lambda \, < \, 1 \,$ 

\smallskip \noindent 
which is exactly (2.17). When $\, \lambda \leq 0 ,\,$  the inequality (2.21) remains
valid but under the hypothesis $\,  v-u < 0. \,$  Then we have 

\smallskip \noindent   (2.22)  $\qquad   \displaystyle 
1 \,\, \geq \,\, \lambda \,  \psi \Big( {{1}\over{\lambda}} 
\Big) \,\, \geq \,\,   \psi ( \lambda )  \,\, \geq \,\,  \lambda  \,,\qquad \lambda  \leq 0
\,$ 

\smallskip \noindent 
which establishes (2.18).

\smallskip \noindent $ \bullet \qquad $ 
We suppose now that the inequalities (2.17).(2.18) hold. We claim that (2.22) is valid because
if we change $\lambda $ into $ {{1}\over{\lambda}}\,$   in the inequality (2.18) we get 

\smallskip \noindent  $  \displaystyle 
{{1}\over{\lambda}}  \,\, \leq \,\,  \psi \Big( {{1}\over{\lambda}} \Big)  \,\, \leq \,\,  
{{1}\over{\lambda}} \, \psi(\lambda) \,,\qquad \lambda \leq 0 \,$

\smallskip \noindent 
which establishes the two first inequalities of (2.22). The end of the proof can be
obtained without difficulty. $\hfill \square $ 

\bigskip \noindent 
We summarize the last two properties (relations (2.10), (2.17), (2.18)) into 
the following one :

\bigskip \noindent  {\bf Proposition  2.4}   

\noindent
	If we suppose that the properties 2.1 to 2.5 are valid, the limiter $\, \psi \,$  defined in
(2.6) satisfies 

\smallskip \noindent   (2.23)  $\qquad   \displaystyle 
1 \,\, \leq \,\, \lambda \,\, \psi \Big( {{1}\over{\lambda}} \Big)  
 \,\, \leq \,\, \psi( \lambda )  \,\, \leq \,\,\lambda  \,\qquad \lambda  \,\, \geq \,\,  1
\,$ 

\smallskip \noindent   (2.24)  $\qquad   \displaystyle 
\lambda \,\, \leq \,\,  \psi( \lambda )   \,\, \leq \,\,  2 \,  \lambda   \qquad \qquad 
\qquad  \quad  0 \,\, \leq \,\, \lambda \,\, \leq \,\, 1 \,$ 

\smallskip \noindent   (2.25)  $\qquad   \displaystyle 
\lambda \,\, \leq \,\,  \psi( \lambda )   \,\, \leq \,\, 0 \,  \qquad \qquad \qquad \quad
\,\, \lambda \,\, \leq \,\, 0 \,. \,$ 

\bigskip \noindent 
Figure 2.3 represents the constraints associated with (2.23).(2.25) ( except the second 
inequality of (2.23) !).

\bigskip 
\centerline { \epsfysize=6cm    \epsfbox  {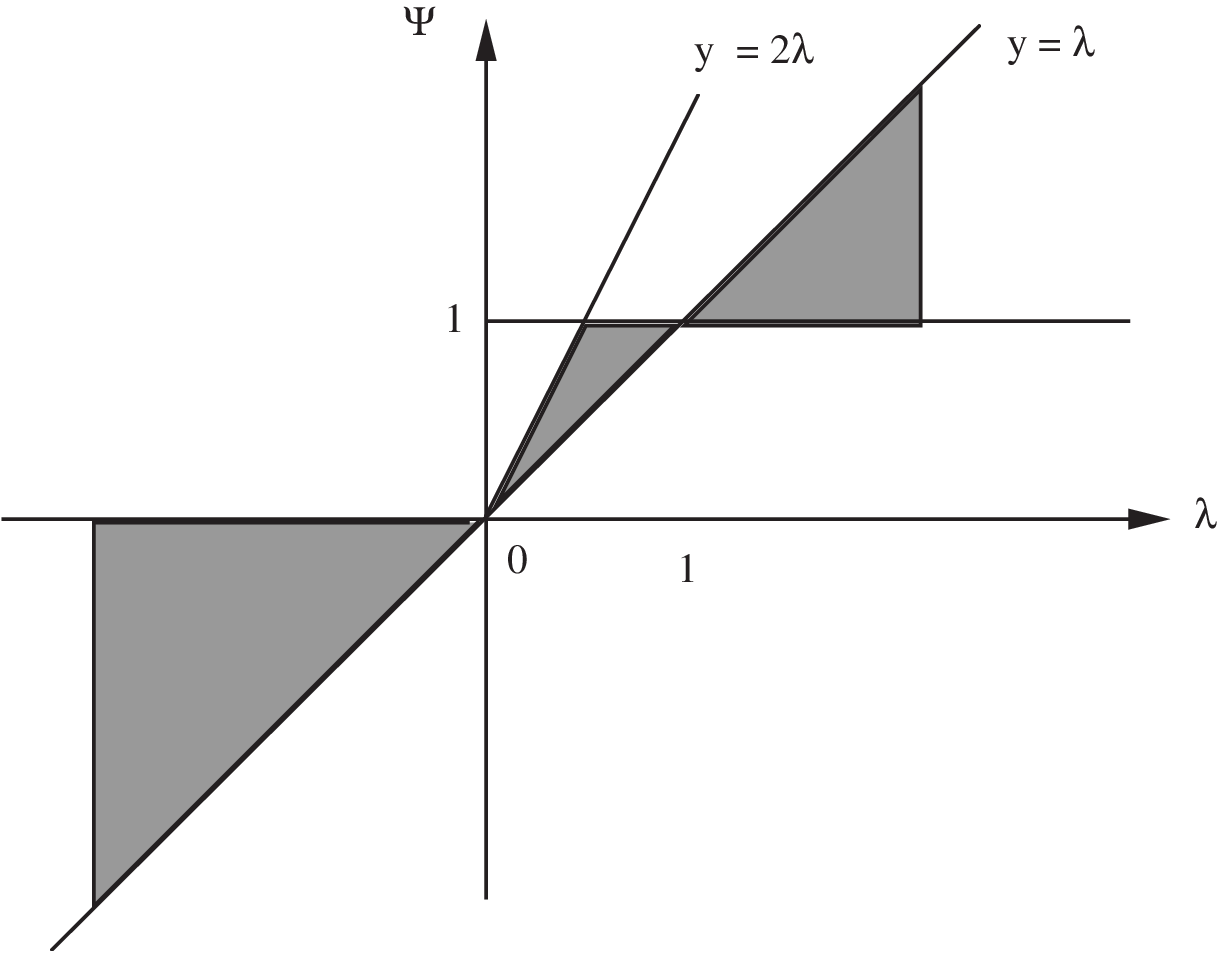} }
\smallskip  \smallskip

% \centerline {
\noindent  {\bf Figure 2.3.} \quad  Constraints of monotonicity and 
convexity for the limiter  function $\, \psi \,$  defined
at (2.6). The graph of $\, \lambda \longmapsto \psi(\lambda) \,$  lies inside the
shaded region.

\bigskip \noindent  {\bf Remark 2.1.}

\noindent
The interpolated values $\, u_{j-{{1}\over{2}}}^- \,$ and $\, u_{j-{{1}\over{2}}}^+ \,$ 
satisfy (according to (2.1), (2.2), (2.6), (2.7)) :
 
\smallskip \noindent   (2.26)(a)  $\qquad   \displaystyle 
u_{j-{{1}\over{2}}}^- \,\,=\,\, u_j \,+\, {{1}\over{2}} \, \psi(\lambda_j) \, ( u_j - 
u_{j-1} ) \,$ 
\smallskip \noindent   (2.26)(b)  $\qquad   \displaystyle 
u_{j-{{1}\over{2}}}^+ \,\,=\,\, u_j \,-\, {{1}\over{2}} \, \lambda_j \, 
\psi \Big( {{1}\over{\lambda_j}} \Big)  \, ( u_j -  u_{j-1} ) \,$ 

\smallskip \noindent   
with  $\qquad   \displaystyle \lambda_j \,\, =\,\, {{u_{j+1} - u_j}\over{u_{j} - u_{j-1}}}
\,$ 

\smallskip \noindent   
and following {\it e.g.} Sweby [1984] we can introduce the inverse $\, r_j \,$    of
the ratio $\, \lambda_j \,$ :

\smallskip \noindent   (2.27)  $\qquad   \displaystyle 
 r_j \,\,=\,\, {{ u_j -  u_{j-1} } \over{u_{j+1} - u_j}} \,$ 

\smallskip \noindent 
and rewrite the limiter $\, \psi \,$   in terms of the variable $\, r ,\,$   by defining 

\smallskip \noindent   (2.28)  $\qquad   \displaystyle 
\varphi(r) \,\, \equiv \,\, r \, \, \psi \Big( {{1}\over{r}} \Big) \,. \,$ 

\smallskip \noindent  
We emphasize that there is a priori no reason to suppose$\, \varphi \equiv \psi .\,$  The two
first equalities of (2.26) take the form 

\smallskip \noindent   (2.29)(a)  $\qquad   \displaystyle 
 u_{j+{{1}\over{2}}}^{-} \,\,=\,\, u_j  \,+\, {{1}\over{2}} \, \varphi(r_j) \, (u_{j+1} -
u_j) \,$
\smallskip \noindent   (2.29)(b)  $\qquad   \displaystyle 
 u_{j-{{1}\over{2}}}^{+} \,\,=\,\, u_j  \,-\, {{1}\over{2}} \, \varphi \Big( {{1}\over{r_j}}
\Big) \,  \, (u_{j} - u_{j-1}) \, . \,$

\smallskip \noindent
Moreover, the inequalities derived in Proposition 2.4 can be translated in terms of the
limiter $\, \varphi \,$  defined in (2.28).

\bigskip \noindent  {\bf Proposition  2.5}   

\noindent
Let $\, \varphi \,$  be defined by (2.28). The limiter $\, \psi \,$  satisfies the
inequalities (2.23)-(2.25) if and only if the function $\, \varphi \,$  verifies :

\smallskip \noindent   (2.30)  $\qquad   \displaystyle 
\max \, (0 ,\, r ) \,\, \leq \,\, \varphi(r) \,\, \leq \,\, 1 \,,\qquad -\infty \,< \, r \,
\leq \, 1 \,$  
\smallskip \noindent   (2.31)  $\qquad   \displaystyle 
1 \,\, \leq \,\, \varphi(r)  \,\, \leq \,\,  \min \, (r ,\, 2 )\,,\qquad \quad   
r \, \geq \, 2 \,$  
\smallskip \noindent   (2.32)  $\qquad   \displaystyle 
\varphi(r)  \,\, \leq \,\,   r \,\, \varphi \big( {{1}\over{r}} \big) \,,\qquad \qquad \qquad
\,\,\,  r \, \geq \, 1 \,.  \,$  

\smallskip \noindent  {\bf Proof of Proposition 2.5}  

\noindent  $ \bullet \qquad $ 
We first suppose $\, r \leq 0 .\,$   Thus the inequality (2.25) joined with (2.28) shows
that we have  $ \, 0 \leq \varphi(r) \leq 1 \,$  and (2.30) is established in this case.  If
$\, r \geq 0 ,\,$  the convexity property reduces to (2.17) and the monotonicity property
corresponds to (2.10). By changing $\, \lambda \,$   into $\, {{1}\over{r}} ,\,$    we get 

\smallskip \noindent   (2.33)  $\qquad   \displaystyle 
r \,\, \leq \,\, r \,\,  \varphi \big( {{1}\over{r}} \big) \,\, \leq \,\, 
 \varphi(r) \,\, \leq \,\, 1 \,,\qquad 0 \, \leq \, r \, \leq \, 1 \,$  
 
\smallskip \noindent 
which achieves the proof of (2.30). 	

\smallskip  \noindent  $ \bullet \qquad $ 
If $\, r \geq 0 ,\,$   we get  from (2.33) :

\smallskip \noindent   (2.34)  $\qquad   \displaystyle 
1 \,\, \leq \,\, \varphi(r)  \,\, \leq \,\, r \,,\qquad 0 \, \leq \, r \, \geq \, 1 \,$  

\smallskip \noindent  
and also (2.32). On the other hand condition (2.10) can be rewritten by changing
$\, \lambda \,$   into $\, {{1}\over{r}} \, ; \,$ we deduce :

\smallskip \noindent  $  \displaystyle 
\varphi(r) \,\, \equiv \,\, r \,\, \psi \big( {{1}\over{r}} \big) \,\, \leq \,\, 2 \,,\qquad  
  r \, \geq \, 1 \,. \,$  

\smallskip \noindent
The latter inequality, joined with (2.34) proves that (2.31) holds.
 The other side of the equivalence is established by similar arguments.   $\hfill \square $

\bigskip \noindent  {\bf Example 2.2.}

\noindent 
Let $\, \psi \,$  be defined by the relations

\smallskip \noindent   (2.35)(a)  $\qquad   \displaystyle 
\psi(\lambda) \,\, = \,\, 0 \,,\qquad \qquad \quad \,\, \lambda \, \leq \, 0 \,$ 
\smallskip \noindent   (2.35)(b)  $\qquad   \displaystyle 
\psi(\lambda) \,\, = \,\, 2 \, \lambda  \,,\qquad \qquad \,\, \,
0 \,\leq\, \lambda \, \leq \, {{1}\over{3}} \,$ 
\smallskip \noindent   (2.35)(c)  $\qquad   \displaystyle 
\psi(\lambda) \,\, = \,\, {{1}\over{2}} \,+\, {{1}\over{2}} \,  \lambda  \,,\qquad 
{{1}\over{3}} \,\leq\, \lambda \, \leq \,3 \,$ 
\smallskip \noindent   (2.35)(d)  $\qquad   \displaystyle 
\psi(\lambda) \,\, = \,\, 2\,     \,,\qquad  \qquad  \quad \, \lambda \, \geq \,3 \, . \,$ 

\smallskip \noindent 
This corresponds to the original MUSCL extrapolation (Van Leer [1977]) and the 
$\, \kappa \,$  limiters af Anderson, Thomas, Van Leer [1986] (see also Chakravarthy et al
[1985]) correspond to the choice

\smallskip \noindent   (2.36)(a)  $\qquad   \displaystyle 
\psi_{\kappa}(\lambda) \,\, = \,\, 0 \,,\qquad \qquad \quad \,\, \lambda \, \leq \, 0 \,$ 
\smallskip \noindent   (2.36)(b)  $\qquad   \displaystyle 
\psi_{\kappa}(\lambda) \,\, = \,\, {{1}\over{2}} \, \bigg \{ \, (1-\kappa) \, \min \, \Big( 1
\,,\, {{3-\kappa}\over{1-\kappa}} \, \lambda \Big) \,+\, 
(1+\kappa) \, \min \, \Big( \lambda \,,\, {{3-\kappa}\over{1-\kappa}} \Big) \, \bigg\}
\,$

\smallskip \noindent  
if $\, \lambda \geq 0 \,.\,$ All of those limiters  $\, (-1 \leq \kappa \leq 1 )
\,$  satisfy $\, \varphi \equiv \psi \,$   (see (2.28)) and the inequalities (2.23).(2.25) ;
we remark that in the latter relations, we do {\bf not} impose the values of $\, \psi \,$ 
corresponding to negative arguments to be null. This last condition is necessary if we
interpolate u by a {\bf linear} function inside the interval $\, {\rm I}_j \,$  and suppose
moreover the condition of monotonicity. We remark also that the so-called Van Albada limiter
[1982], defined by

\smallskip \noindent   (2.37)  $\qquad   \displaystyle 
\psi(\lambda) \,\,=\,\, {{\lambda \, (1 + \lambda)}\over{(1 + \lambda)^2}} \,\, \equiv \,\,
\lambda \,\, \psi \Big( {{1}\over{\lambda}} \Big) \,,\qquad \lambda \in \R \,$

\smallskip \noindent 
satisfies neither (2.23).(2.25)  [for $\lambda \leq  -1$]  nor (2.30).(2.32)  [for  $-1 \leq
\lambda \leq 0$].

\bigskip \noindent  {\bf Proposition  2.6}   

\noindent
	The limiter $\, L(u,v,w) \,$ satisfying the properties 2.1 to 2.3 is Lipschitz continuous,
{\it i.e.}

\setbox21=\hbox {$\displaystyle  
\exists \, C > 0 \,, \, \forall \, u,\, v,\, w ,\, u' ,\, v' ,\, w' \in \R , \,  $}
\setbox22=\hbox {$\displaystyle  
\mid L(u,\,v,\,w) - L( u',\, v' ,\, w' ) \mid \,\, \leq \,\, C \, \big( \mid u-u' \mid \,+\, 
\mid v -v' \mid \,+\, \mid w - w' \mid \big) \,$}
\setbox30= \vbox {\halign{#&# \cr \box21 \cr \box22    \cr   }}
\setbox31= \hbox{ $\vcenter {\box30} $}
\setbox44=\hbox{\noindent  (2.38) $\displaystyle  \,\,  
\left\{ \box31 \right. $}   
\smallskip \noindent   $ \box44 $ 

\smallskip \noindent 
if and only if the limiter functions $\, \psi \,$  (defined at (2.6)) and $\, \varphi \,$
(defined at (2.28)) satisfy a similar property~:

\smallskip \noindent   (2.39)(a)  $\qquad   \displaystyle 
\exists \, K > 0 \,, \, \forall \,\lambda ,\, \lambda'  \in \R , \quad 
\mid \psi (\lambda') - \psi (\lambda) \mid \, \,\leq \,\, K \, \mid \lambda' - \lambda \mid $
\smallskip \noindent   (2.39)(b)  $\qquad   \displaystyle 
\exists \, K > 0 \,, \, \forall \,r ,\, r'  \in \R , \quad 
\mid \varphi (r') - \varphi (r) \mid \, \,\leq \,\, K \, \mid r' - r \mid  \, . \, $

\bigskip \noindent  {\bf Proof of Proposition  2.6}   

\noindent $ \bullet $ \qquad 
We suppose first that (2.38) holds. We take $u = -1 $, $v = 0 $, $w = \lambda $ and analogous
relations for the ``prime'' variables. Then we have (according to (2.6)) :

\smallskip \noindent  $ \displaystyle 
{{1}\over{2}} \, \big( \psi (\lambda') - \psi (\lambda)  \big) \,\,=\,\, L ( -1 ,\, 0 ,\, 
\lambda' - 1 ) \,-\, L ( -1 ,\, 0 ,\, \lambda - 1 ) \,.\,$

\smallskip \noindent 
 Similarly, we have :

\smallskip \noindent  $ \displaystyle 
{{1}\over{2}} \, \big( \varphi (r') - \varphi (r)  \big) \,\,=\,\, L ( -r' ,\, 0 ,\,1 ) 
\,-\, L ( -r ,\, 0 ,\,  1 ) \, \,$

\smallskip \noindent
which establishes (2.39).	

\smallskip  \noindent $ \bullet $ \qquad 
We suppose that (2.39) is realized. For a set of given real numbers
$ u,\, v,\, w$ (respectively $u',\, v',\, w'$) , we define $\, \lambda ,\, r \,$  (resp. $\,
\lambda' ,\, r'  $) by the formulae  $\, w = v + \lambdaÊ(v - u)  $ ; $    vÊ=ÊuÊ+ÊrÊ(w-v)   $
(resp. $ w'Ê=Êv'Ê+Ê\lambda'Ê(v' - u')   $ ; $  v'Ê=Êu'Ê+Êr'Ê(w'-v') ) $ and we decompose the
left hand side of (2.38) as follows :

\setbox21=\hbox {$\displaystyle  
L( u',\, v' ,\, w' )  - L(u,\,v,\,w) \,\,=\,\,  v' - v  \,+\,   $}
\setbox22=\hbox {$\displaystyle  \qquad \qquad 
\,+\, {{1}\over{2}} \, \big(  \varphi (r') - \varphi (r)  \big) \, (w' - w) \,+\, {{1}\over{2}}
\, \varphi (r)  \, \big( (w'-w) - (v'-v) \big) \,$}
\setbox30= \vbox {\halign{#&# \cr \box21 \cr \box22    \cr   }}
\setbox31= \hbox{ $\vcenter {\box30} $}
\setbox44=\hbox{\noindent  (2.40) $\displaystyle  \quad 
\left\{ \box31 \right. $}   
\smallskip \noindent   $ \box44 $ 

\smallskip \noindent 
If $ r'Ê=Ê\infty $  (resp. $ rÊ=Ê\infty $) then $w'Ê=Êv'$ (resp. $wÊ=Êv$) and (2.38) follows 
directly from (2.40) and (2.30).(2.31). In the other cases, the estimate 

\smallskip \noindent   (2.41)  $\qquad   \displaystyle 
\mid r' - r \mid \,\, \mid w'-v' \mid \,\, \leq \,\,  \mid u' - u  \mid \,+\, (1 + \mid r
\mid) \, \mid v' - v \mid \,+\, \mid r \mid \, \mid w' - w \mid \,$

\smallskip \noindent 
joined to (2.39) establishes the property (2.38) under the restriction 

\smallskip \noindent   (2.42)  $\qquad   \displaystyle 
\mid r \mid \,\, \leq \,\, 1 \,.\,$ 

\smallskip \noindent 
The same conclusion holds if we reverse the roles of $r$ and $r'.$ 
 Inequality (2.38) is established if we suppose

\smallskip \noindent   (2.43)  $\qquad   \displaystyle 
\mid \lambda \mid \,\, \leq \,\, 1 \,\, $ and \quad $ \mid \lambda' \mid \,\, \leq \,\, 1 \,.\,
$

\smallskip \noindent  
We have in this case

\setbox21=\hbox {$\displaystyle  
L( u',\, v' ,\, w' )  - L(u,\,v,\,w) \,\,=\,\,  v' - v  \,+\,   $}
\setbox22=\hbox {$\displaystyle  \qquad  \qquad 
\,+\, {{1}\over{2}} \,   \psi ( \lambda) \, (v'-v \,+\, u' - u ) \,+\, 
{{1}\over{2}} \, \big(   \psi ( \lambda') - \psi ( \lambda) \big) \, (v'-u') \,$}
\setbox30= \vbox {\halign{#&# \cr \box21 \cr \box22    \cr   }}
\setbox31= \hbox{ $\vcenter {\box30} $}
\setbox44=\hbox{\noindent  (2.44) $\displaystyle  \quad 
\left\{ \box31 \right. $}   
\smallskip \noindent   $ \box44 $ 

\smallskip \noindent 
and according to (2.10) and (2.39) it is sufficient to remark that we have 

\smallskip \noindent   $  \displaystyle 
\mid \lambda' - \lambda \mid \,\, \mid v' - v \mid \,\, \leq \,\, (1 + \mid  \lambda \mid )
\, \mid v' - v \mid \,+\, \mid \lambda \mid \,  \mid u' - u \mid \,+\, \mid w' - w \mid  \,$

\smallskip \noindent 
to obtain finally :

\smallskip \noindent   $  \displaystyle 
\mid  L( u',\, v' ,\, w' )  - L(u,\,v,\,w) \mid \,\, \leq \,\, \Big( 1 + {{K}\over{2}} \Big)
\, ( \mid u' - u \mid \,+\,2 \, \mid v' - v \mid \,+\,  \mid w' - w  \mid ) \,$ 

\smallskip \noindent 
which ends the proof. $\hfill \square $

\bigskip 
\vfill \eject    %%%%%%%%%%%%%   fd le 15 juin 2010  
\noindent  {\bf Remark 2.2.}

\noindent
	The total variation of the ``reconstructed sequence''

\smallskip \noindent   (2.45)  $\qquad   \displaystyle 
\widetilde{u} \,\, \equiv \,\, \cdots ,\, u_{j-1} \,,\, u_{j-{{1}\over{2}}}^- \,,\, 
 u_{j-{{1}\over{2}}}^+ \,,\, u_j \,,\, u_{j+{{1}\over{2}}}^-  \,,\,  u_{j+{{1}\over{2}}}^+ 
\,,\,   u_{j+1} \,,\, \cdots \qquad (j \in \Z) \,$

\smallskip \noindent  
is naturally defined by

\smallskip \noindent   (2.46)  $\qquad   \displaystyle 
{\rm TV}(\widetilde{u}) \,\, \equiv \,\, \sum_{j\in \Z} \big( \mid u_j - 
u_{j-{{1}\over{2}}}^+ \mid \,+\, \mid  u_{j+{{1}\over{2}}}^- - u_j \mid \big) 
\,+\,  \sum_{j\in \Z} \mid u_{j+{{1}\over{2}}}^+ - u_{j+{{1}\over{2}}}^- \mid \,.\, $

\smallskip \noindent 
The properties 2.1 to 2.5 do {\bf not} imply that the total variation of the reconstructed
sequence is smaller than TV(u) ({\it c.f.} (1.23)). For example, in the particular case of the
following simple sequence 

\smallskip \noindent   (2.47)(a)  $\qquad   \displaystyle 
u_j \,\,=\,\, -3 \, h \,,\qquad j \leq -1 \,$
\smallskip \noindent   (2.47)(b)  $\qquad   \displaystyle 
u_0 \,\,=\,\, 0 \,,\qquad u_1 \,\,=\,\,   h \,$
\smallskip \noindent   (2.47)(c)  $\qquad   \displaystyle 
u_j \,\,=\,\, 4 \, h \,,\qquad j \geq 2 \,$

\smallskip \noindent  
the reconstructed sequence associated with the MUSCL scheme (2.35) 
admits a total variation given by

\smallskip \noindent   (2.48)  $\qquad   \displaystyle 
{\rm TV}(\widetilde{u}) \,\, = \,\, {\rm TV}(u) \,+\, 2 \, h \,.\,$ 

\smallskip \noindent  
(see Figure 2.4) 

\bigskip 
\centerline { \epsfysize=5cm    \epsfbox  {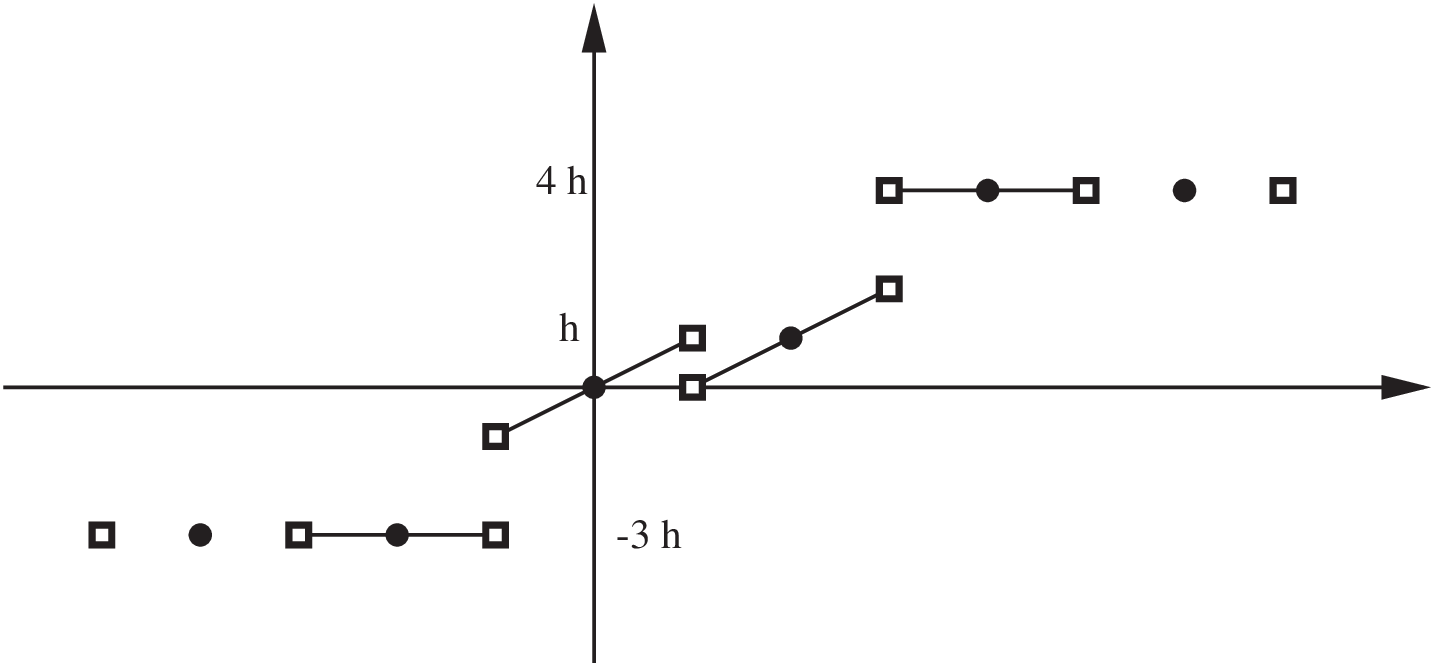} }
\smallskip  \smallskip

% \centerline {
\noindent  {\bf Figure 2.4.} \quad 
The total variation of the reconstructed sequence (2.38) with
the MUSCL scheme (squares) is greater than the TV of the initial sequence (points).

% \titredroite={\pecaps  Decrease of the Total Variation   }
\toppagetrue  
\botpagetrue    

\bigskip   \bigskip  
\vfill \eject    %%%%%%%%%%%%%   fd le 15 juin 2010 
\noindent {\bf 3. \quad Decrease of the Total Variation  }

\smallskip \noindent
In this section, we suppose that the interpolated values at the interfaces\br
  $\,  (j+{{1}\over{2}}) \, h \,$  are given by the relations (2.27).(2.29) and are
submitted to the restrictions described in the Part II  (properties 2.1 to 2.5). We study the
following differential system

\smallskip \noindent   (3.1)  $\qquad   \displaystyle 
{{{\rm d}u_j}\over{{\rm d}t}} \,+\,  {{1}\over{h}} \, \Big( f_{j+{{1}\over{2}}} \,-\, 
 f_{j-{{1}\over{2}}} \Big) \,\,=\,\, 0 \,,\qquad  j \in \Z \,,\quad t> 0 \, \,$

\smallskip \noindent 
with the flux at the interface  $\,  (j+{{1}\over{2}}) \, h \,$ given by a monotone flux
function $ \Phi $ :

\smallskip \noindent   (3.2)  $\qquad   \displaystyle 
f_{j+{{1}\over{2}}} \,\,=\,\, \Phi (u_{j+{{1}\over{2}}}^- \,,\, u_{j+{{1}\over{2}}}^+ )
\,,\qquad  j \in \Z \,. \,$

\smallskip \noindent 
We suppose that the initial datum

\smallskip \noindent   (3.3)  $\qquad   \displaystyle 
u_j(0) \,\,=\,\, u_j^0 \,,\qquad  j \in \Z \,  \,$

\smallskip \noindent  
satisfies

\smallskip \noindent   (3.4)  $\qquad   \displaystyle 
u^0 \, \in \ell^{\infty}(\Z) \, \cap \, {\rm BV} (\Z) \,,\,$ 

\smallskip \noindent  
{\it id est } 

\smallskip \noindent   (3.5)(a)  $\qquad   \displaystyle 
\sup_{j \in \Z} \, \mid  u_j^0 \mid \,\, < \,\, \infty \,$ 
\smallskip \noindent   (3.5)(b)  $\qquad   \displaystyle 
\sum_{j \in \Z} \, \mid u_{j+1} - u_j^0 \mid  \,\, < \,\, \infty \,. \,$ 

\smallskip \noindent  
Then the method of lines is well defined, at least locally in time.

\bigskip 
% \vfill \eject 
\noindent  {\bf   Proposition  3.1}   

\noindent 
The semi-discrete numerical scheme defined by the differential system (3.1).(3.2) and
(2.27).(2.29) associated with the initial condition (3.3).(3.4) has a unique solution
$\, u_j(t) ,\,$  $\,  j \in \Z, \,$ $  tÊ\inÊ[0 ,T], \,$ for some $\, TÊ>Ê0, \,$ if the flux
function $\, \Phi \,$  is Lipschitz continuous on the compact sets of $\, \R \, :$

\setbox21=\hbox {$\displaystyle  
\forall \, M > 0 ,\,  \exists \, C > 0 \,, \, \forall \, u,\, v,   \, u' ,\, v' \in \R , \, 
\, \big(  \mid u \mid \leq M ,\,  \mid v \mid \leq M ,\, \mid u' \mid \leq M  ,\, \, $}
\setbox22=\hbox {$\displaystyle  \quad 
 \mid v' \mid \leq M \big)\, \Longrightarrow \, 
\mid \Phi(u',\,v'  ) - \Phi(u,\,v) \mid \,\, \leq \,\, C \, \big( \mid u-u' \mid
\,+\,  \mid v -v' \mid  \big) \,. \,$}
\setbox30= \vbox {\halign{#&# \cr \box21 \cr \box22    \cr   }}
\setbox31= \hbox{ $\vcenter {\box30} $}
\setbox44=\hbox{\noindent  (3.7) $\displaystyle  \quad  
\left\{ \box31 \right. $}   
\smallskip \noindent   $ \box44 $

\bigskip  
\noindent  {\bf Proof of Proposition  3.1}   

\noindent  $ \bullet \qquad $ 
From the classical Cauchy-Lipschitz theorem it is sufficient to prove that the difference of
the flux

\smallskip \noindent   (3.8)  $\qquad   \displaystyle 
F_j(u) \,\,=\,\, -{{1}\over{h}} \, \Big( f_{j+{{1}\over{2}}} \,-\, 
 f_{j-{{1}\over{2}}} \Big) \,,\qquad  j \in \Z \,,\quad u \in \ell^{\infty}(\Z) \,$

\smallskip \noindent  
is a Lipshitz continuous $\,  \ell^{\infty}(\Z) \longrightarrow  \ell^{\infty}(\Z)  \,$ 
function. According to (3.5), we suppose that $\, \mid u_j \mid \leq M .\,$  Then the
same property is realized for the intermediate states $\, u_{j+{{1}\over{2}}}^\pm \,$  
  according to the property 2.4 (monotonicity). We have :

\smallskip \noindent  $  \displaystyle 
\mid F_j(u)   \mid  \,\, \leq \, {{1}\over{h}} \, \Big[ \, \mid \Phi ( u_{j+{{1}\over{2}}}^-
\,,\, u_{j+{{1}\over{2}}}^+ ) -  \Phi ( u_{j+{{1}\over{2}}}^- \,,\, u_j )   \mid  
+ \mid \Phi ( u_{j+{{1}\over{2}}}^- \,,\, u_{j} ) -  \Phi (  u_{j} \,,\,  u_{j}  )  \mid \,$ 

\smallskip \noindent  $  \displaystyle \hfill 
+ \mid \Phi (  u_{j} \,,\,  u_{j}  ) - \Phi ( u_{j} \,,\, u_{j-{{1}\over{2}}}^+  )   \mid 
+ \mid \Phi ( u_{j} \,,\, u_{j-{{1}\over{2}}}^+  )  -  \Phi (
 u_{j-{{1}\over{2}}}^- \,,\,  u_{j-{{1}\over{2}}}^+ )  \mid  \,  \Big] \,  $

\smallskip \noindent  $  \displaystyle \qquad 
\leq \,\,  {{C}\over{h}} \, \Big[ \, \mid  u_{j+{{1}\over{2}}}^+ - u_j \mid 
+ \mid  u_{j+{{1}\over{2}}}^- - u_j \mid +  \mid  u_{j-{{1}\over{2}}}^+ - u_j \mid 
+ \mid  u_{j-{{1}\over{2}}}^- - u_j \mid \,     \Big] \,  $

\smallskip \noindent  $  \displaystyle \qquad 
\leq \,\,  {{C}\over{h}} \, \Big[ \, \mid 1 - {{1}\over{2}} \varphi \big( {{1}\over{r_{j+1}}}
\big) \mid \,  \mid u_{j+1} - u_j \mid \,  \,\,+\,\,    {{1}\over{2}} \varphi \big(
{{1}\over{r_j}} \big) \,  \mid u_{j+1} - u_j \mid \,\,+\,\,    $

\smallskip \noindent  $  \displaystyle \qquad \qquad \qquad 
+ \,  {{1}\over{2}} \varphi (r_j) \, \mid u_j - u_{j-1} \mid \,\,+\,\,  
\mid 1 - {{1}\over{2}} \varphi \big( {{1}\over{r_{j-1}}}
\big) \mid \,\,   \mid u_{j} - u_{j-1} \mid  \,     \Big] \,  $

\smallskip \noindent  $  \displaystyle \qquad 
\leq \,\, 16 \, {{C}\over{h}} \,    M \,$ 

\smallskip \noindent
that establishes that $ F $  is a function $\,  \ell^{\infty}(\Z) \longrightarrow 
\ell^{\infty}(\Z)  .\,$

\smallskip  \noindent  $ \bullet \qquad $ 
We prove that  $ F $   is Lipschitz continuous on the bounded sets of $\, 
\ell^{\infty}(\Z). \,$    Let  $\, (u_j)_{j\in\Z} ,\, $ $   (v_j)_{j\in\Z} \,$ be two bounded
sequences in $\,  \ell^{\infty}(\Z) \, :  \,$

\smallskip \noindent   (3.9)  $\qquad   \displaystyle 
\mid u_j\mid  \,\,\leq\,\,  M \,,\quad \mid v_j\mid  \,\,\leq\,\, M \,,\qquad 
\forall \, j \in \Z \,.\,$ 

\smallskip \noindent 
According to the monotonicity property (2.8), the interpolated values
$\, u_{j+{{1}\over{2}}}^\pm \,,\,$ $\, v_{j+{{1}\over{2}}}^\pm \,,\,$  satisfy the same
property. From the inequality (3.7) and the proposition 2.6 we deduce 

\smallskip \noindent  $  \displaystyle 
 \mid \Phi ( u_{j+{{1}\over{2}}}^- \,,\, u_{j+{{1}\over{2}}}^+ ) -   
\Phi ( v_{j+{{1}\over{2}}}^- \,,\, v_{j+{{1}\over{2}}}^+ )  \mid \,\, \leq \,\, \big( \mid 
u_{j+{{1}\over{2}}}^- - v_{j+{{1}\over{2}}}^- \mid \,+\, 
\mid  u_{j+{{1}\over{2}}}^+ - v_{j+{{1}\over{2}}}^+ \mid \big) \, $

\smallskip \noindent  $  \displaystyle \qquad 
\leq \,\, K \, C \, \big( \mid u_{j-1} - v_{j-1} \mid  \,+\, 2 \, \mid u_{j} - v_{j} \mid 
\,+\, 2 \, \mid u_{j+1} - v_{j+1} \mid  \,+\,  \mid u_{j+1} - v_{j+1} \mid  \big) \,$ 

\smallskip \noindent  $  \displaystyle \qquad 
\leq \,\, 6 \, K \, C \, \, \sup_{\ell \in \Z} \, \mid u_{\ell} - v_{\ell} \mid \,.\,$ 

\smallskip \noindent 
Then we have the estimate

\smallskip \noindent   (3.10)  $\qquad   \displaystyle 
\mid F_j(u) - F_j(v) \mid \,\, \leq \,\, {{12 \, K \, C } \over{h}} \, 
 \sup_{\ell \in \Z} \, \mid u_{\ell} - v_{\ell} \mid \, \,$ 

\smallskip \noindent 
and the proposition is established.
$\hfill \square$

\bigskip \noindent  {\bf Remark 3.1.} 

\noindent
From the Proposition 3.1, we know that the function $\,  [0,T] \ni t \longmapsto 
(u_j(t))_{j\in\Z} \,$  which belongs to $\, \ell^{\infty}(\Z)  \,$   is derivable~; in
particular, the sequence $\, (u_j(t))_{j\in\Z} \,$    is bounded uniformly in space and
time~:

\smallskip \noindent   (3.11)  $\qquad   \displaystyle 
\exists \, M > 0 \,,\, \forall \, t \in [0 ,\, T] \,,\, \forall \, j \in \Z \,,\, \mid u_j(t)
\mid \, \leq \, M \, .\,$ 

\bigskip \noindent $ \bullet \qquad $ 
We focus now on the total variation of $\, (u_j(t)) \,$   ({\it c.f.} (1.13)). We suppose
moreover the following restrictions:

\bigskip \noindent  {\bf  Property 3.1.} 

\noindent
We fix a real $\, \alpha \,$  such that $\, 1 \leq \alpha \leq 2 .\,$ We suppose that the
limiter $\, \varphi \,$  related to the interpolated values $\,
u_{j+{{1}\over{2}}}^\pm \,$   according to (2.29) satisfies the following inequalities~:

\smallskip \noindent   (3.12)  $\qquad   \displaystyle 
\varphi(r) \,\, \leq \,\, (\alpha - 2) \, r \,,\qquad r \leq 0 \,$ 
\smallskip \noindent   (3.13)  $\qquad   \displaystyle 
\varphi(r) \,\, \leq \,\,  \alpha  \,,\qquad \qquad \quad r >0 \, . \,$

\bigskip \noindent  {\bf Remark 3.2.} 

\noindent
If $\, \varphi \,$   satisfies the conditions (2.30).(2.32), (3.12).(3.13) for the particular
value $ \alpha = 2 $  and the complementary restriction : 

\smallskip \noindent   (3.14)  $\qquad   \displaystyle 
\varphi(r) \,\, \leq \,\, 2 \, r \,,\qquad 0 \leq r \leq 1 \,, \,$ 

\smallskip \noindent  
 then the restrictions that we obtain for $\, \varphi\,$   correspond {\bf exactly} to the
``second order TVD region'' derived previously by SwebyÊ[1984]. For the other values of $\, 
\alpha, \,$  ($  1 \leq \alpha < 2 $) the estimates (3.12).(3.13) are new. We can transcribe
the Property 3.1  and the inequality (3.14) in terms of the $\psi$-limiter ({\it c.f.}
(2.26))~:

\smallskip \noindent   (3.15)(a)  $\qquad   \displaystyle 
\psi(\lambda) \,\, \geq \,\,  \alpha - 2   \,,\qquad \lambda \leq 0 \,$ 
\smallskip \noindent   (3.15)(b)  $\qquad   \displaystyle 
\psi(\lambda) \,\, \leq \,\,  \alpha \, \lambda  \,,\qquad \quad 0 \leq  \lambda \leq 1 \,$ 
\smallskip \noindent   (3.15)(c)  $\qquad   \displaystyle 
\psi(\lambda) \,\, \leq \,\, 2    \,,\qquad \qquad   \lambda \geq 1 \, . \,$

% \vfill \eject

% \titredroite={\pecaps  Decrease of the Total Variation   }
\toppagetrue  
\botpagetrue

\bigskip \noindent  {\bf   Theorem  3.1. \quad Decrease of the total variation}   

\noindent 
Let  $\, (u_j(t)) \, (j\in\Z,\, t \in [0,\,T]) \,$   be defined by the method of
lines (the hypotheses are the ones given at Proposition 3.1). We suppose also that the limiter
function $\, \varphi \,$  that defines the interpolated values at mid-points according to
(2.29), satisfies (2.30).(2.32) and Property 3.1 (relations (3.12).(3.13)).  Moreover the flux
function $\, \Phi $ of relation (3.2) is supposed to be monotone 

\smallskip \noindent   (3.16)  $\qquad   \displaystyle 
\Phi_u(u,\,v) \geq 0 \,,\quad \Phi_v(u,\,v) \leq 0 \, $ 

\smallskip \noindent 
and absolutely continuous relatively to each variable~: 

\smallskip \noindent  $  \displaystyle 
\Phi(u,\,v) \,\,=\,\, \Phi(u_0 \,,\,v) \,+\, \int_{u_0}^{u} \Phi_u (\xi,\,v) \,{\rm d}\xi
\,,\quad 
\Phi(u,\,v) \,\,=\,\, \Phi(u \,,\, v_0) \,+\, \int_{v_0}^{v} \Phi_v (u,\,\xi) \,{\rm d}\xi
\,.\,$

\smallskip \noindent 
Then, the total variation $\, {\rm TV}(t) \,$ of $\, (u_j(t)) \,$ defined in (1.13) is a
decreasing function of time~:

\smallskip \noindent   (3.17)  $\qquad   \displaystyle 
{\rm TV}(t+\theta) \,\, \leq \,\, {\rm TV}(t) \,\qquad 0 \, \leq \, t 
\, \leq \, t+\theta \, \leq \,  T \,.\,$

\bigskip 
\centerline { \epsfysize=6cm    \epsfbox  {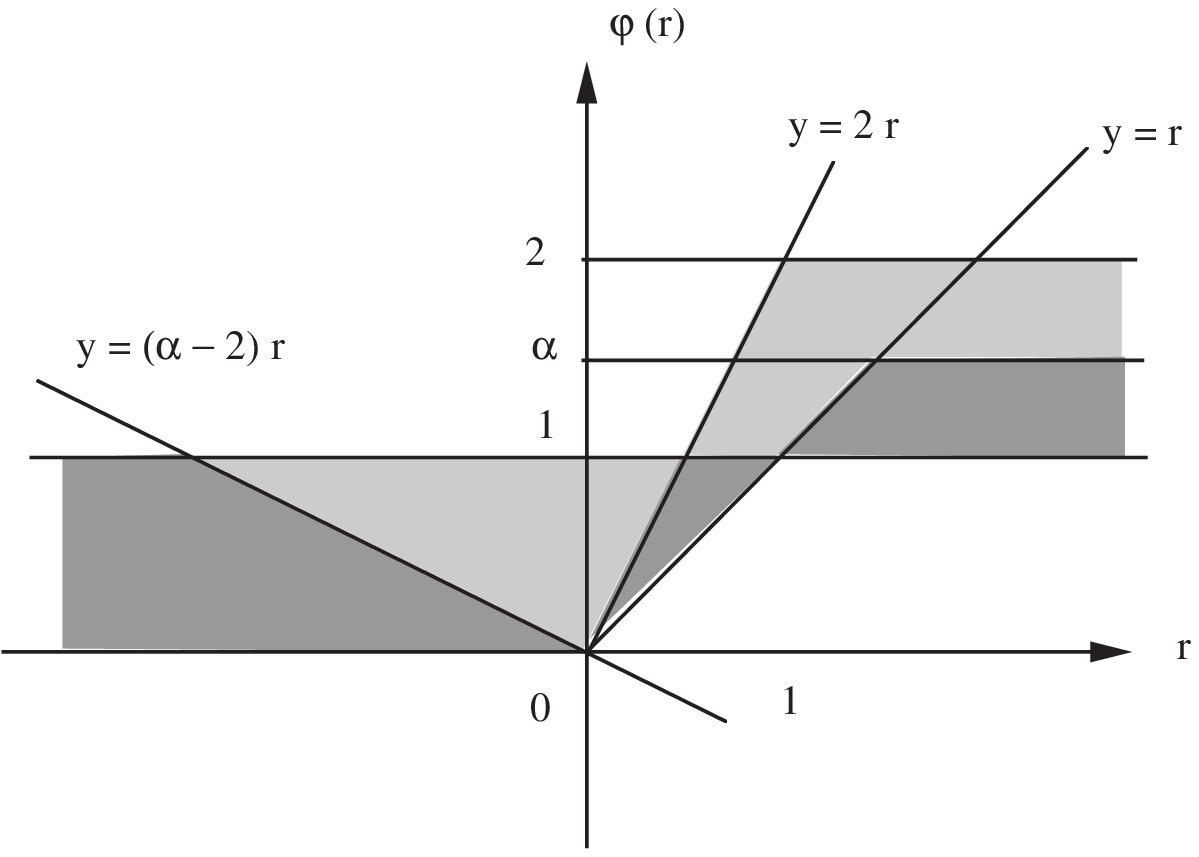} }
\smallskip  \smallskip

% \centerline {
\noindent  {\bf Figure 3.1.} \quad Admissible regions for the limiter$\, \varphi .\,$ 
 The dark shaded region corresponds to the conditions (3.12).(3.14) 
and the clear shaded set completes
it to 	describe (2.30) (2.31).

\bigskip 
% \vfill \eject         %%%%%%%%   fd 15 juin 2010 
\noindent  {\bf   Proof of Theorem  3.1.}

\noindent $   \bullet $ \qquad 
We fix some notations :

\smallskip \noindent   (3.18)  $\qquad   \displaystyle 
\varphi_j \,\, \equiv \,\, \varphi(r_j) \,,\quad \widetilde{\varphi}_j \,\, \equiv \,\,
\varphi \big( {{1}\over{r_j}} \big) \,,\qquad j \in \Z  \,$ 
\smallskip \noindent   (3.19)  $\qquad   \displaystyle 
\Delta_{j+{{1}\over{2}}} \,\, \equiv \,\, u_{j+1} - u_j \,,\qquad  j \in \Z  \,. \,$ 

\smallskip \noindent 
Then (2.29) admits the following form~:

\smallskip \noindent   (3.20)(a)  $\qquad   \displaystyle 
u_{j+{{1}\over{2}}}^- \,\, =\,\, u_j \,+\, {{1}\over{2}} \, \varphi_j \,\,
\Delta_{j+{{1}\over{2}}} \,$
\smallskip \noindent   (3.20)(b)  $\qquad   \displaystyle 
u_{j+{{1}\over{2}}}^+ \,\, =\,\, u_{j+1} \,-\, {{1}\over{2}} \, \widetilde{\varphi}_{j+1}  \,\,
\Delta_{j+{{1}\over{2}}} \,. \,$

\smallskip \noindent 
We also denote by $\, \alpha_{j-{{1}\over{2}}} \,$  (respectively  $\, \beta_{j+{{1}\over{2}}
} $) the mean value of $\, \Phi_u \,$   (respectively   $\, \Phi_v  $) over the domain 
$\, [\, u_{j-{{1}\over{2}}}^- \,,\, u_{j+{{1}\over{2}}}^- \, ] \times  $  $\,
\{u_{j-{{1}\over{2}}}^+ \}
$   (resp. $\, \{u_{j+{{1}\over{2}}}^- \} $  $\times  
[\, u_{j-{{1}\over{2}}}^+ \, ,\, u_{j+{{1}\over{2}}}^+ \, ]) $ :

\smallskip \noindent   (3.21)  $\qquad   \displaystyle  
\Phi ( u_{j+{{1}\over{2}}}^- \,,\, u_{j-{{1}\over{2}}}^+ ) \,-\, 
\Phi ( u_{j-{{1}\over{2}}}^- \,,\, u_{j-{{1}\over{2}}}^+ ) \,\, \equiv \,\, 
\alpha_{j-{{1}\over{2}}} \,( u_{j+{{1}\over{2}}}^- -  u_{j-{{1}\over{2}}}^- ) \,$

\smallskip \noindent   (3.22)  $\qquad   \displaystyle  
\Phi ( u_{j+{{1}\over{2}}}^- \,,\, u_{j+{{1}\over{2}}}^+ ) \,-\, 
\Phi ( u_{j+{{1}\over{2}}}^- \,,\, u_{j-{{1}\over{2}}}^+ ) \,\, \equiv \,\, 
\beta_{j+{{1}\over{2}}} \,( u_{j+{{1}\over{2}}}^+ -  u_{j-{{1}\over{2}}}^+ ) \,. \,$

\smallskip \noindent 
According to the monotonicity (3.16) of the flux $\, \Phi$, we have~:

\smallskip \noindent   (3.23)  $\qquad   \displaystyle  
\alpha_{j+{{1}\over{2}}} \,\geq \, 0 \,,\quad \beta_{j+{{1}\over{2}}} \, \leq \, 0 \,,\qquad 
j \in \Z \,.\,$ 

\smallskip \noindent $   \bullet $ \qquad 
When we derive formally the total variation (1.13) with respect to time,
we have

\smallskip \noindent   $  \displaystyle  
{{\rm d}\over{{\rm d}t}} {\rm TV}(t) \,\, = \,\, \sum_{j \in \Z} \, {\rm sgn}\, (u_{j+1} - u_j)
  \, {{\rm d}\over{{\rm d}t}}(u_{j+1} - u_j) \,.\, $

\smallskip \noindent  
Then we consider the following sequence

\smallskip \noindent   (3.24)  $\qquad   \displaystyle  
p_{j+{{1}\over{2}}} \,\,=\,\, -   {\rm sgn}\, (u_{j+1} - u_j)
  \, {{\rm d}\over{{\rm d}t}}(u_{j+1} - u_j) \,,\qquad j \in \Z \,.\,$

\smallskip \noindent  
From the equalities (3.1), (3.2), (3.20).(3.22), we have  the elementary calculus~:

\smallskip \noindent  $ \displaystyle  
-h {{{\rm d}u_j}\over{{\rm d}t}} \,\, =\,\, - 
\Phi ( u_{j+{{1}\over{2}}}^- \,,\, u_{j+{{1}\over{2}}}^+ ) \,+\,
\Phi ( u_{j-{{1}\over{2}}}^- \,,\, u_{j-{{1}\over{2}}}^+ ) \,\,$
\smallskip \noindent  $ \displaystyle  \quad 
\,\, =\,\, -  \Phi ( u_{j+{{1}\over{2}}}^- \,,\, u_{j+{{1}\over{2}}}^+ ) 
\,+\, \Phi ( u_{j+{{1}\over{2}}}^- \,,\, u_{j-{{1}\over{2}}}^+ )  
\, - \, \Phi ( u_{j+{{1}\over{2}}}^- \,,\, u_{j-{{1}\over{2}}}^+ )  
\,+\, \Phi ( u_{j-{{1}\over{2}}}^- \,,\, u_{j-{{1}\over{2}}}^+ ) \,\,$
\smallskip \noindent  $ \displaystyle  \quad  \,\,=\,\, 
- \beta_{j+{{1}\over{2}}} \,( u_{j+{{1}\over{2}}}^+ -  u_{j-{{1}\over{2}}}^+ ) 
\,-\,  \alpha_{j-{{1}\over{2}}} \,( u_{j+{{1}\over{2}}}^- -  u_{j-{{1}\over{2}}}^- ) \,$
\smallskip \noindent  $ \displaystyle  \quad  \,\,=\,\, 
- \beta_{j+{{1}\over{2}}} \, \Big[ \big( u_{j+1} - {1\over2} \widetilde{\varphi}_{j+1} \,
\Delta_{j+{1\over2}} \big) - \big(  u_{j} - {1\over2} \widetilde{\varphi}_{j} \,
\Delta_{j-{1\over2}} \big) \Big] \,$
\smallskip \noindent  $ \displaystyle  \qquad  \qquad  
\,\,-\,\alpha_{j-{1\over2}} \, \Big[ \big( u_j + {1\over2} \varphi_j \, \Delta_{j+{1\over2}} 
\big) - \big(   u_{j-1} + {1\over2} \varphi_{j-1} \, \Delta_{j-{1\over2}} \big) \Big] \,$
\smallskip \noindent  $ \displaystyle  \quad  \,\,=\,\, 
- \beta_{j+{{1}\over{2}}}  \, \Big( 1 - {{1}\over{2}}  \widetilde{\varphi}_{j+1} +
{{1}\over{2}}   \widetilde{\varphi}_{j} \, r_{j} \Big) \, \Delta_{j+{1\over2}}  \, 
-\,  \alpha_{j-{{1}\over{2}}}  \,  \Big( 1 + {{1}\over{2}} \, {{\varphi_{j}}\over{r_{j}}} -
{{1}\over{2}} \varphi_{j+1}   \Big) \, \Delta_{j-{1\over2}}  \Big] \,. \,$

\smallskip \noindent 
Then 

\smallskip \noindent  $ \displaystyle  
p_{j+{{1}\over{2}}} \,\,=\,\, - {{ \displaystyle \mid \Delta_{j+{1\over2}} \mid }\over{
\displaystyle \Delta_{j+{1\over2}} }} \, \Big[ \, {{\rm d}\over{{\rm d}t}}\big( u_{j+1} - u_j
\big) \, \Big] \,$ 
\smallskip \noindent  $ \displaystyle  \quad  \,\,=\,\, 
 {{ \displaystyle \mid \Delta_{j+{1\over2}} \mid }\over{ \displaystyle \Delta_{j+{1\over2}}
}} \, {{1}\over{h}} \, \, \bigg[ \, 
 \beta_{j+{{3}\over{2}}} \, \Big( 1 - {{1}\over{2}}   \widetilde{\varphi}_{j+2}
+{{1}\over{2}}   \widetilde{\varphi}_{j+1} \, r_{j+1} \Big) \, \Delta_{j+{3\over2}} \,$
\smallskip \noindent  $ \displaystyle  \qquad \quad 
+\,  \alpha_{j+{{1}\over{2}}} \, \Big( 1 + {{1}\over{2}} \, {{\varphi_{j+1}}\over{r_{j+1}}} -
{{1}\over{2}} \varphi_j   \Big) \, \Delta_{j+{1\over2}} \,-\, 
  \beta_{j+{{1}\over{2}}}  \, \Big( 1 - {{1}\over{2}}  \widetilde{\varphi}_{j+1} +
{{1}\over{2}}   \widetilde{\varphi}_{j} \, r_{j} \Big) \, \Delta_{j+{1\over2}}\,$
\smallskip \noindent  $ \displaystyle  \qquad \qquad \qquad 
- \, \alpha_{j-{1\over2}}   \Big( 1 + {1\over2} \, {{\varphi_{j}}\over{r_{j}}} -
{1\over2} \varphi_{j+1}   \Big) \,   \Delta_{j-{1\over2}} \, \bigg] \,.  $

\smallskip \noindent 
and

\setbox21=\hbox {$\displaystyle  
p_{j+{{1}\over{2}}} \,\,=\,\, {{ \mid \Delta_{j+{{1}\over{2}}} \mid } \over{h}} \, \, \bigg[ 
\,\, 
{{\beta_{j+{{3}\over{2}}}}\over{r_{j+1}}}\, \Big( 1 - {{1}\over{2}}   \widetilde{\varphi}_{j+2}
+{{1}\over{2}}   \widetilde{\varphi}_{j+1} \, r_{j+1} \Big) \,+\,    $}
\setbox22=\hbox {$\displaystyle  \qquad \qquad \qquad \qquad \qquad  
\,+\,  \alpha_{j+{{1}\over{2}}} \, \Big( 1 + {{1}\over{2}} \, {{\varphi_{j+1}}\over{r_{j+1}}} -
{{1}\over{2}} \varphi_j   \Big)   $}
\setbox23=\hbox {$\displaystyle  \qquad
- \beta_{j+{{1}\over{2}}}  \, \Big( 1 - {{1}\over{2}}  \widetilde{\varphi}_{j+1} +
{{1}\over{2}}   \widetilde{\varphi}_{j} \, r_{j} \Big) \,-\, 
\alpha_{j-{{1}\over{2}}} \,r_j \,  \Big( 1 + {{1}\over{2}} \, {{\varphi_{j}}\over{r_{j}}} -
{{1}\over{2}} \varphi_{j+1}   \Big) \, \bigg] \,.  $}
\setbox30= \vbox {\halign{#&# \cr \box21 \cr \box22    \cr \box23    \cr  }}
\setbox31= \hbox{ $\vcenter {\box30} $}
\setbox44=\hbox{\noindent  (3.25) $\displaystyle  \quad 
\left\{ \box31 \right. $}   
\smallskip \noindent   $ \box44 $ 

\smallskip \noindent 
We integrate (3.25) relatively to time  $\, (0 \leq t \leq \tau \leq t + \theta \leq   T),\,$
we sum the result for $\,   j =Êk,Êk+1,\dots, \, \ell \,$  and we factorize $\,
\alpha_{j+{{1}\over{2}}}  \,$ and $\, \beta_{j+{{1}\over{2}}} \,$  (Abel transformation). We
obtain :

\setbox21=\hbox {$\displaystyle  
h \, \sum_{j=k}^{\ell} \, \big( \mid \Delta_{j+{{1}\over{2}}} (t)  \mid - \mid 
\Delta_{j+{{1}\over{2}}} (t+\theta) \mid \big) \,\,= \,\,    $}
\setbox22=\hbox {$\displaystyle  \,\, = \,\, 
 \sum_{j=k}^{\ell} \,\,   \, \int_t^{t+\theta} \mid \Delta_{j+{{1}\over{2}}} (\tau) 
\mid  \alpha_{j+{{1}\over{2}}} \, \Big( 1 - {{r_j}\over{\mid r_j \mid}} \Big) \, \Big( 1 +
{{1}\over{2}}  {{\varphi_{j+1}}\over{r_{j+1}}} - {{1}\over{2}} \varphi_j \Big) \, {\rm d}
\tau  \,$}
\setbox23=\hbox {$\displaystyle  \qquad \qquad 
+ \, \int_t^{t+\theta} \mid \Delta_{\ell+{{1}\over{2}}} (\tau) 
\mid  \alpha_{\ell+{{1}\over{2}}} \,
\Big( 1 + {{1}\over{2}}  {{\varphi_{\ell+1}}\over{r_{\ell+1}}} - {{1}\over{2}} \varphi_{\ell}
\Big) \, {\rm d} \tau  \,$} 
\setbox24=\hbox {$\displaystyle   \qquad \qquad 
- \, \int_t^{t+\theta} \mid \Delta_{k-{{1}\over{2}}} (\tau) 
\mid  \alpha_{k-{{1}\over{2}}} \,   {{r_k}\over{\mid r_{k} \mid}}  \,
\Big( 1 + {{1}\over{2}}  {{\varphi_{k}}\over{r_{k}}} - {{1}\over{2}} \varphi_{k-1}
\Big)  \, {\rm d} \tau  \,$}
\setbox25=\hbox {$\displaystyle    
+ \,  \sum_{j=k+1}^{\ell} \, \int_t^{t+\theta} \mid \Delta_{j+{{1}\over{2}}} (\tau) 
\mid  \beta_{j+{{1}\over{2}}} \, \Big( {{r_j}\over{\mid r_j \mid}} - 1  \Big) \, \Big( 1 -
{{1}\over{2}} \widetilde{\varphi}_{j+1} +  {{1}\over{2}} \widetilde{\varphi}_{j} r_j  \Big)
\, {\rm d} \tau  \,$}
\setbox26=\hbox {$\displaystyle   \qquad   
+ \, \int_t^{t+\theta} \mid \Delta_{\ell+{{3}\over{2}}} (\tau) 
\mid  \beta_{\ell+{{3}\over{2}}} \,  {{r_{\ell +1}}\over{\mid r_{\ell +1} \mid}}   \,
\Big( 1 -  {{1}\over{2}} \widetilde{\varphi}_{\ell +2 } +  {{1}\over{2}}
\widetilde{\varphi}_{\ell +1} r_{\ell +1}  \Big) \, {\rm d} \tau  \,$}
\setbox27=\hbox {$\displaystyle   \qquad   \qquad  
- \, \int_t^{t+\theta} \mid \Delta_{k+{{1}\over{2}}} (\tau) 
\mid  \beta_{k+{{3}\over{2}}} \, 
\Big( 1 -  {{1}\over{2}} \widetilde{\varphi}_{k +1 } +  {{1}\over{2}}
\widetilde{\varphi}_{k  } r_{k }  \Big) \, {\rm d} \tau  \, . \,$}
\setbox30= \vbox {\halign{#&# \cr \box21 \cr \box22  \cr \box23 \cr \box24 \cr \box25 \cr
\box26 \cr \box27 \cr     }}
\setbox31= \hbox{ $\vcenter {\box30} $}
\setbox44=\hbox{\noindent  (3.26) $\displaystyle  \quad  \left\{ \box31 \right. $}   
\smallskip \noindent   $ \box44 $ 

\smallskip \noindent 
The first, second, fourth and sixth terms of the right hand side of (3.26) are 
positive, due to (3.23) and the following inequalities~:

\smallskip \noindent   (3.27)  $\qquad   \displaystyle  
1 + {{1}\over{2}}  {{\varphi_{j+1}}\over{r_{j+1}}} - {{1}\over{2}} \varphi_j \,\, \geq \,\, 0
\,,\qquad j \in \Z \,$ 
\smallskip \noindent   (3.28)  $\qquad   \displaystyle  
1 - {{1}\over{2}} \widetilde{\varphi}_{j+1} +  {{1}\over{2}} \widetilde{\varphi}_{j}\,  r_j 
\,\, \geq \,\, 0 \,,\qquad j \in \Z \, . \,$ 

\smallskip \noindent $   \bullet $ \qquad 
We establish (3.27). We distinguish two cases, according to the 
sign of $\, r_{j+1} .\,$   If $\, r_{j+1} \geq 0 , \,$  we have, due to (3.13) and
(2.30).(2.31)~:

\smallskip \noindent   $  \displaystyle  
{{1}\over{2}} \varphi_j \,\, \leq \,\, {{1}\over{2}} \alpha \,\, \leq \,\, 1 \,\, \leq \,\,
1 \,+\, {{1}\over{2}}   {{\varphi_{j+1}}\over{r_{j+1}}} \,$

\smallskip \noindent 
which proves (3.27) in this case. On the other hand, if$\, r_{j+1} < 0 , \,$  we have, 
due to (3.12)~:

\smallskip \noindent   $  \displaystyle  
{{1}\over{2}} \varphi_j \,- \, {{1}\over{2}}  {{\varphi_{j+1}}\over{r_{j+1}}} 
\,\, \leq \,\, {{1}\over{2}} \alpha \,-\, {{1}\over{2}}  (\alpha - 2) \,\,=\,\, 1 \,$

\smallskip \noindent 
and (3.27) is established. Inequality (3.28) is equivalent , according to (3.18), to~:

\smallskip \noindent   $  \displaystyle  
1 - {{1}\over{2}} \varphi(\lambda_{j+1}) \,+\, {{1}\over{2}} {{1}\over{\lambda_{j}}}
\,  \varphi \Big({{1}\over{\lambda_{j}}} \Big) \,, \qquad 
\lambda_j \,=\, {{1}\over{r_j}} \,$

\smallskip \noindent 
which corresponds exactly to (3.27).

\smallskip \noindent $   \bullet $ \qquad 
	The third and fifth terms of the right hand side of (3.26) are bounded by

\smallskip \noindent   (3.29)  $\qquad   \displaystyle  
\gamma \,\, = \,\,  2 \, M \, (2 + K ) \, \theta \,$

\smallskip \noindent 
according to (3.11), (3.7), (2.30), (2.31), (2.39) and the constant $\, \gamma \,$  is
independent of $k$ and $\, \ell .\,$  When $k$ tends to  $\, - \infty \,$   and  $\, \ell \,$ 
to $\, + \infty ,\,$ the previous inequalities give the following estimate~: 

\smallskip \noindent   (3.30)  $\qquad   \displaystyle  
{\rm TV}(t) \,-\, {\rm TV}(t+\theta)  \,\, \geq \,\,-  4 \, M \, (2 + K) \, \theta \,. \,$

\smallskip \noindent 
By taking $tÊ=Ê0$ and $ \theta \leq T ,\,$   this inequality proves that $\,  {\rm
TV}( \theta)  \,$ is {\bf well defined} for all $\, \theta  \leq T . \,$  Moreover, the
third and fifth terms of the right hand side  of (3.26) tend to zero from the previous
considerations and the Lebesgue dominated convergence Theorem. Taking $\, kÊ\longrightarrow
- \infty \,$ and $\, \ell \longrightarrow +\infty \,$   in (3.26), we finally deduce~:

\smallskip \noindent   $ \displaystyle  
{\rm TV}(t) \,-\, {\rm TV}(t+\theta)  \,\, \geq \,\, 0 \,,\qquad \theta \geq 0 \,$ 

\smallskip \noindent 
which is exactly (3.17). $\hfill \square $

\bigskip \noindent  {\bf  Remark 3.3. }   

\noindent 
	Since the work of Harten [1983], the TVD property is derived from an incremental 
form of the scheme (3.1)  Moreover, the inequalities (3.27).(3.28) have been previously
considered among others by Sweby [1984], Chakravarthy-Osher [1983] and Spekreijse [1987]. Here
these inequalities are obtained essentially as a {\bf consequence} of other considerations
concerning the limiter, and not as necessary conditions to assume the TVD property. In
particular, the limiter function $\, \varphi(r) \, $ is {\bf not} necessarily null for the
negative values of the argument $r$ as it is usually supposed (see Figure 3.1). To be complete,
the sufficient conditions proposed by Spekreijse [1987] (see also Cahouet-Coquel [1989]) allow
possible non-null values for  $\, \varphi(r) \, $  when $\, rÊ<Ê0 \,$  but their conditions are
{\bf not}  {\it a priori} compatible with the hypothesis of monotonicity. 

\bigskip   \noindent 
We turn now to an important result of this paper.

\bigskip \noindent  {\bf   Theorem  3.2. \quad $\ell^{\infty}$ stability. }   

\noindent 
Let $\,  u_j(t) \,$  be the solution of the differential system (3.1).(3.3); we assume that all
the hypotheses of the Theorem 3.1 are satisfied. Then we have 

\smallskip \noindent   (3.31)  $\qquad   \displaystyle  
\inf_{\ell \in \Z} \, u_\ell^0 \,\, \leq \,\,  u_j(t) \,\, \leq \,\, \sup_{\ell \in \Z} 
 \, u_\ell^0 \,,\qquad j \in \Z \,,\quad 0 \,\leq \, t \, \leq \, T \, .\,$

\bigskip
\vfill \eject    %%%%%%%%%%%%%   fd le 15 juin 2010 
 \noindent  {\bf  Remark 3.4. }   

\noindent 
An immediate corollary of the Theorem 3.2 is the fact that the solution $\, u_j(t)  \,$ is (for
each $h$ fixed) defined globally in time, {\it i.e.} the time $\, TÊ\geq Ê0\,$  introduced at
the Proposition 3.1 can be taken equal to $\, +\infty .\,$ 

\bigskip  \noindent  $   \bullet $ \qquad 
The proof of Theorem 3.2 contains essentially the ideas given previously by Sanders [1983] and
Osher [1985]. The new difficulty here is to work with an infinite set of indices $ j $  (Osher
supposed the periodicity in space). So we first prove the following property~:

\bigskip \noindent  {\bf Proposition 3.2. }   

\noindent 
For each time $t ,\, $ $ 0 \leq t \leq T ,\,$   the sequence $\, (u_j(t))_{j \in \Z} \,$  
defined by (3.1).(3.3) has a limit $ M(t) $ (respectively  $ m(t)$) for $\, j
\longrightarrow + \infty $   (respectively  $\, j \longrightarrow - \infty $) and we have 

\smallskip \noindent   (3.32)  $\qquad   \displaystyle  
M(t) \, \equiv \, M(0) \,,\qquad  m(t) \, \equiv \, m(0) \,,\quad  0 \, \leq \, t \, \leq \,
T \, .\,$

\bigskip \noindent  {\bf Proof of Proposition 3.2. }   

\noindent   $   \bullet $ \qquad 
The existence of $M(t)$ and $m(t)$ is a direct consequence of the convergence of the series
that defines the total variation for each time t ({\it c.f.} Theorem 3.1).

\smallskip \noindent $   \bullet $ \qquad 
 We now consider the mean value of $\, u_j(t) \, : \,$ 

\smallskip \noindent  $ \displaystyle  
v_N(t) \,\,=\,\,  {{1}\over{N}} \, \sum_{j=1}^N  u_j(t) \,,\qquad N \geq 1 \,,\quad  
0 \, \leq \, t \, \leq \, T \, .\,$ 

\smallskip \noindent 
The function $\, v_N \,$  is derivable and we have

\smallskip \noindent   (3.33)  $\qquad   \displaystyle  
{{{\rm d}v_N}\over{{\rm d}t}} \,+\, {{1}\over{N}} \, \big( f_{N+{{1}\over{2}}} - 
f_{{{1}\over{2}}} \big) \,.\,$ 

\smallskip \noindent 
From the Remark 3.1 and the continuity of the flux relatively to 
$  \, \{u_j(t) \} \,$  (see the proof of the proposition 3.1), there exists some constant
$\, \kappa \,$ such that

\smallskip \noindent   (3.34)  $\qquad   \displaystyle  
\mid {{{\rm d}v_N}\over{{\rm d}t}} \mid \,\, \leq \,\, {{\kappa}\over{N}} \,.\,$ 

\smallskip \noindent 
Then  $\, {{{\rm d}v_N}\over{{\rm d}t}} \,$  tends uniformly to zero as $N$ tends to $ +
\infty .\,$  Moreover from the Cesaro convergence Theorem,$\, v_N \,$  tends to $M(t)$ as
$\, N \longrightarrow \infty .\,$ We deduce from these two facts that $M(t)$ is derivable, and 

\smallskip \noindent   (3.35)  $\qquad   \displaystyle 
{{{\rm d}M}\over{{\rm d}t}} \,\,=\,\, 0 \,.\,$ 

\smallskip \noindent
The first equality of (3.32) is established. The proof is analogous for $\, m(t).	  \hfill
\square $

\bigskip 
\vfill \eject    %%%%%%%%%%%%%   fd le 15 juin 2010 
\noindent  {\bf Proof of Theorem 3.2. }   

\noindent   $   \bullet $ \qquad
 It is sufficient to prove (3.31) for $t$  small enough to establish the property. For
$\, t \in  [0,T], \,$ we define

\smallskip \noindent  $ \displaystyle 
\mu(t) \,\,=\,\, \inf_{j \in \Z} \, u_j(t) \,.\,$ 

\smallskip \noindent 
Then $\, \longmapsto \mu(t) \,$  is continuous and we distinguish two cases~:

\smallskip \noindent {\bf (i) } \qquad 
$\, \mu < m .\,$   Then for some $\, j_0(t) ,\,$  we have

\smallskip \noindent   (3.36)  $\qquad   \displaystyle 
\mu(t) \,\,=\,\, u_{j_0} (t) \,.\,$ 

\smallskip \noindent 
With the notations introduced in the proof of Theorem 3.1, we have~:

\setbox21=\hbox {$\displaystyle  
h \,  {{{\rm d}u_j}\over{{\rm d}t}} \,=\, - \beta_{j+{{1}\over{2}}} \, (u_{j+1} - u_j) \,
 \Big( 1 - {{1}\over{2}}  \widetilde{\varphi}_{j+1} +
{{1}\over{2}}   \widetilde{\varphi}_{j} \, r_{j} \Big) \,$ }
\setbox22=\hbox {$\displaystyle  \qquad \qquad \qquad 
+ \,\, \alpha_{j-{{1}\over{2}}} \, (u_{j-1} - u_j) \,  \Big( 1 + {{1}\over{2}} 
{{\varphi_j}\over{r_j}} - {{1}\over{2}} \varphi_{j-1}  \Big) \,  $}
\setbox30= \vbox {\halign{#&# \cr \box21 \cr \box22    \cr   }}
\setbox31= \hbox{ $\vcenter {\box30} $}
\setbox44=\hbox{\noindent  (3.37) $\displaystyle \qquad  \left\{ \box31 \right. $}   
\smallskip \noindent   $ \box44 $ 

\smallskip \noindent   
and (3.37) is non-negative for $\, jÊ=Ê j_0. \,$ Consequently for $\, t'Ê\geq Êt \,$  and $\, 
t'-t \,$ sufficiently small, we have

\smallskip \noindent   (3.38)  $\qquad   \displaystyle  
\mu(t') \,\, \geq \,\, \mu(t) \,,\qquad \forall \, t' \geq t \,,\quad t'-t \,$  small enough. 

\smallskip \noindent {\bf (ii) } \qquad 
$\, \mu \geq  m .\,$  If the alternative (i) is realized at $tÊ=Ê0,$ the property (3.38)
corresponds exactly to the desired result. Otherwise, we have

\smallskip \noindent   (3.39)  $\qquad   \displaystyle  
\mu(0) \,\, = \,\, m \,\, \equiv \,\, \inf_{j \in \Z} \, u_j^0 \,$ 

\smallskip \noindent  
and if (3.38) is in default, we consider $\, \theta \,$  such that $\, \mu(\theta) <
m \,$   and $ \, \tau \, $ defined by 

\smallskip \noindent   (3.40)  $\qquad   \displaystyle  
\tau \,\,=\,\, \sup \, \{ t,\,\, t \,\leq \, \theta \,, \,\, \mu(t) \geq m \, \} \,.$ 

\smallskip \noindent  
We have  $\, 0Ê\leq \tau < \theta \,$ Ê   and $\, \mu(\tau)  Ê=Êm \,$  by continuity. According
to (3.38), $\, \mu \,$  is non-decreasing inside the interval $\,  [\tau ,\, \theta  ]. \,$ So
that

\smallskip \noindent   (3.41)  $\qquad   \displaystyle  
\mu(\theta) \, \geq \, \mu(\tau) \,$ 

\smallskip \noindent 
which establishes the contradiction.
	The proof concerning the maximum bound is analogous.	$\, \hfill \square $

% \vfill \eject
\bigskip   \bigskip  
\noindent {\bf 4. \quad Convergence, Order of Accuracy }

\smallskip \noindent
In this section the parameter $\,  h > 0 \,$  tends to zero. Two questions are naturally
related to this problem~: 1)	What is the order of the truncation error of the scheme
(3.1)(3.2), (2.27).(2.29) as $h$ tends to zero ? 2)	Does the function $\, u_h(x,\,t) \,$
defined by

% \titredroite={\pecaps  Convergence, Order of Accuracy   }
\toppagetrue  
\botpagetrue    

\smallskip \noindent   (4.1)  $\qquad   
u_h(x,\,t) \,\,=\,\, u_j^h (t) \,,\qquad \big(j-{{1}\over{2}}\big) \,h \,\leq x \, \leq 
 \big(j+{{1}\over{2}}\big) \, h \,$

\smallskip \noindent 
(where the superscript $h$ recalls the dependence towards $h$ of the solution $u_j(t)$ of the
method of lines) converge to the unique entropy solution of (1.1).(1.2)  ? The following
property answers the first point.

\bigskip \noindent  {\bf Proposition 4.1. }   

\noindent 
We suppose that the entropy solution $u(x,t)$ of (1.1).(1.2) is regular around some point
$\, (x_0Ê,\,Êt_0). \,$  Let $\, v_j^h(t) \,$ be defined by

\smallskip \noindent   (4.2)  $\qquad   \displaystyle  
 v_j^h(t) \,\,=\,\, u(j \,h \,,\, t) \,,\qquad j \in \Z \,,\quad h > 0 \,,\quad t \geq 0
\,.\,$ 

\smallskip \noindent 
The residual of the semi-discrete scheme (3.1)(3.2), (2.27).(2.29), {\it i.e. } the quantity
obtained when we apply the scheme (3.1) to the exact ``nodal'' values   $\, v_j^h(t) \,:\,$ 

\smallskip \noindent   (4.3)  $\qquad   \displaystyle  
{\rho}_j^h (t) \,\,  \equiv \,\, {{{\rm d} v_j^h}\over{{\rm d}t}} \,-\, {{1}\over{h}} \, \Big(
\, \Phi \big ( v_{j+{{1}\over{2}}}^- \,,\,  v_{j+{{1}\over{2}}}^+  \big) \,-\, 
\Phi \big ( v_{j-{{1}\over{2}}}^- \,,\,  v_{j-{{1}\over{2}}}^+  \big) \, \Big) \,$

\smallskip \noindent 
is second order accurate in space :

\smallskip \noindent   (4.4)  $\qquad   \displaystyle  
\rho_j^h (t) \,\, = \,\, {\rm O} (h^2) \,,\qquad h \longrightarrow 0 \,,\quad \exists \, j_0
\,,\,\, x_0 \,=\, j_0 \, h \,$ 

\smallskip \noindent  
if the limiter function $\, \varphi \,$  of the relation (2.29) satisfies~:

\smallskip \noindent   (4.5)(a)  $\qquad   \displaystyle  
\varphi(1) \,=\, 1 \,,\qquad \varphi \,$ derivable on both sides of $\, +1 \,$ 
\smallskip \noindent   (4.5)(b)  $\qquad   \displaystyle  
\varphi(-1) + \varphi(3) \,\,=\,\, 2 \,$ 

\smallskip \noindent  
and if the monotone flux function $ \, \Phi \,$  is absolutely continuous as in Theorem 3.1.

\bigskip \noindent  {\bf Proof of Proposition 4.1. }   

\noindent $\bullet $ \qquad 
The interpolated values $\, v_{j+{{1}\over{2}}}^\pm \,$   are constructed from the $\, v_j
\,'$s   according to the relations (2.27).(2.29), {\it i.e.}

\smallskip \noindent   (4.6)  $\qquad   \displaystyle  
r_j \,\,=\,\, {{ v_j^h - v_{j-1}^h }\over{ v_{j+1}^h -  v_{j}^h }} \,\, \equiv \,\, {{ 
\Delta _{j-{{1}\over{2}}} }\over{ \Delta_{j+{{1}\over{2}}} }} \,$ 

\smallskip \noindent   (4.7)(a)  $\qquad   \displaystyle 
 v_{j+{{1}\over{2}}}^- (t) \,\, = \,\, v_{j}^h  \,+\, {{1}\over{2}} \, \varphi(r_j) \, 
\Delta_{j+{{1}\over{2}}} \,\, \equiv \,\, v_j  \,+\, {{1}\over{2}} \, \varphi_j \, \,
\Delta_{j+{{1}\over{2}}}  \,$ 
\smallskip \noindent   (4.7)(b)  $\qquad   \displaystyle 
 v_{j-{{1}\over{2}}}^+ (t) \,\, = \,\, v_{j}^h  \,-\, {{1}\over{2}} \, \varphi 
\Big( {{1}\over{ r_j}} \Big ) \, \,
\Delta_{j-{{1}\over{2}}} \,\, \equiv \,\, v_j  \,-\, {{1}\over{2}} \, 
\widetilde{\varphi}_j  \,  \Delta_{j-{{1}\over{2}}}  \, . \,$ 

\smallskip \noindent  
We develop $\, {\rho}_j^h (t) \,$  around $\,(x_0,\, t_0) \,$ assuming that  $\, x_0 \,$ is a
mesh point (see (4.4)) :

\smallskip \noindent   (4.8)  $\qquad   \displaystyle 
x_0 \,=\, j \, h \,\,$ \qquad for some $\, j \in \Z \,.\,$  

\smallskip \noindent 
Two cases can occur, according to the nullity of the gradient  $\,\, {{\partial
u}\over{\partial x}} (x_0,\, t_0) \,.$ 

\smallskip \noindent $\bullet $ \qquad 
$\, \displaystyle u'_j \, \equiv \,  {{\partial u}\over{\partial x}} (x_0,\, t_0)  \,\, \neq
\,\, 0 \,. \,$  We develop the numerical flux 

\smallskip \noindent   (4.9)  $\qquad   \displaystyle 
f_{j+{{1}\over{2}}} \,\,=\,\, \Phi \Big( \, v_j + {{1}\over{2}} \, \varphi_j \, 
\Delta_{j+{{1}\over{2}}} \,,\,  v_j + \big( 1 - {{1}\over{2}} \,\widetilde { \varphi}_{j+1}
\big) \, \Delta_{j+{{1}\over{2}}}  \, \Big) \,$

\smallskip \noindent  
by the Taylor formula, with the notations 

\smallskip \noindent  $  \displaystyle 
\alpha_j \,=\, {{\partial \Phi}\over{\partial u}} (u=v_j \,,\, v=v_j) \,,\quad 
\beta_j \,=\, {{\partial \Phi}\over{\partial v}} (u=v_j \,,\, v=v_j) \,,\quad 
f_j \,=\, \Phi(v_j \,,\, v_j) \,.\, $ 

\smallskip \noindent  
We get 

\smallskip \noindent   (4.10)  $\qquad   \displaystyle 
f_{j+{{1}\over{2}}} \,\,=\,\, f_j \,+\, \alpha_j \, \varphi_j \, \Delta_{j+{{1}\over{2}}}
\,\,+\,\, \beta_j \, \big( 1 - {{1}\over{2}} \widetilde { \varphi}_{j+1} \big) \,\,
\Delta_{j+{{1}\over{2}}} \,+\, {\rm O} \big( \Delta_{j+{{1}\over{2}}}^2 \big) \,$

\smallskip \noindent   
We must now specify the behavior of the $\varphi$'s as $h$ tends to zero. We have on one hand

\smallskip \noindent   (4.11)  $\qquad   \displaystyle 
r_j \,\,=\,\, {{ \displaystyle 1 - {{1}\over{2}} h \,  {{u_j''}\over{u_j'}}  \,+\, {\rm O}
(h^2)  }  \over   {
 \displaystyle 1 + {{1}\over{2}} h \,  {{u_j''}\over{u_j'}}  \,+\, {\rm O} (h^2) }  } 
 \,\,=\,\, 1 - h \,  {{u_j''}\over{u_j'}}  \,+\, {\rm O} (h^2)   \,$

\smallskip \noindent 
with $  \displaystyle  \,\, u_j'' = {{\partial^2 u}\over{\partial x^2}} (x_0 ,\, t_0)  \,\,$ 
and on the other hand

\smallskip \noindent   (4.12)  $\qquad   \displaystyle 
{{1}\over{r_j}} \,\,=\,\, 1 + h \,  {{u_j''}\over{u_j'}}  \,+\, {\rm O} (h^2)   \,.\,$

\smallskip \noindent 
When we substitute these developments into (4.10), with the hypotheses ({\it c.f.} (4.5)) 

\smallskip \noindent   (4.13)(a)  $\qquad   \displaystyle 
\varphi(1+\xi) \,\,=\,\, 1 \,+\, \xi \, \varphi_+' \,+\,  {\rm O} (\xi^2) \,,\qquad \xi
\longrightarrow 0 \,,\quad \xi > 0 \,$ 
\smallskip \noindent   (4.13)(b)  $\qquad   \displaystyle 
\varphi(1-\xi) \,\,=\,\, 1 \,-\, \xi \, \varphi_-' \,+\,  {\rm O} (\xi^2) \,,\qquad \xi
\longrightarrow 0 \,,\quad \xi > 0 \,.\,$ 

\smallskip \noindent 
We obtain (assuming  $\, {{u_j''}\over{u_j'}} > 0 \, ; \,$ the other case is analogous) :

\setbox21=\hbox {$\displaystyle  
f_{j+{{1}\over{2}}} \,\,=\,\, f_j \,+\, {{1}\over{2}} h \, u_j' \, (\alpha_j + \beta_j) \,+\, 
{{1}\over{2}} h^2 \, u_j'' \, ( - \alpha_j \, \varphi_-' \,-\,   \beta_j \,\varphi_+'  )
 \, \,+ \,$ }
\setbox22=\hbox {$\displaystyle  \qquad \qquad \qquad 
+ \, {{h^2}\over{4}} \,  u_j'' \, (\alpha_j + \beta_j) \, 
+ \, {{h^2}\over{8}} \,  (u_j')^2 \, f_j'' \,+\,   {\rm O} (h^3) 
\,  $}
\setbox30= \vbox {\halign{#&# \cr \box21 \cr \box22    \cr   }}
\setbox31= \hbox{ $\vcenter {\box30} $}
\setbox44=\hbox{\noindent  (4.14) $\displaystyle \qquad  \left\{ \box31 \right. $}   
\smallskip \noindent   $ \box44 $ 

\smallskip \noindent 
where  $\,  f_j''  \,$   is defined by the formula :

\smallskip \noindent   (4.15)  $\qquad   \displaystyle 
 f_j'' \,\,=\,\, ( \Phi_{uu} + \Phi_{uv}  + \Phi_{vu} + \Phi_{uu} ) ( v_j \,,\, v_j) \,\,=\,\,
{{{\rm d}^2f}\over{{\rm d}u^2}} (v_j) \, \,.\,$ 

\smallskip \noindent   
Similarly, we easily obtain 

\smallskip \noindent   (4.16)  $\qquad   \displaystyle 
r_{j-1} \,\,\equiv\,\, {{ v_{j-1}^h  -  v_{j-2}^h }\over{ v_{j}^h  -  v_{j-1}^h }}
\,\,= \,\, 1 - h   \,  {{u_j''}\over{u_j'}}  \,+\, {\rm O} (h^2)   \, \,$

\smallskip \noindent   (4.17)(a)  $\qquad   \displaystyle 
\varphi_{j-1} \,\,=\,\, 1  - h   \,  {{u_j''}\over{u_j'}}  \, \varphi_-' \,+\,  {\rm O}
(h^2)   \, \,$
\smallskip \noindent   (4.17)(b)  $\qquad   \displaystyle 
\widetilde{ \varphi}_{j} \,\,=\,\, 1  + h   \,  {{u_j''}\over{u_j'}}  \, \varphi_+' \,+\, 
{\rm O} (h^2)   \, .\,$

\smallskip \noindent 
We insert the latter expressions in the following expression of the numerical flux at 
$\, (j-{{1}\over{2}}) \, h \,:\,$ 

\smallskip \noindent   (4.18)  $\qquad   \displaystyle 
f_{j-{{1}\over{2}}} \,\,=\,\, \Phi \Big( \, v_j + \Big(  {{1}\over{2}}   \varphi_{j-1} - 1
\Big) \,  \Delta_{j-{{1}\over{2}}} \,,\,  v_j   - {{1}\over{2}} \widetilde { \varphi}_{j}
\, \, \Delta_{j-{{1}\over{2}}}  \, \Big) \,$

\smallskip \noindent 
and we obtain finally:
 the other case is analogous) :

\setbox21=\hbox {$\displaystyle  
f_{j-{{1}\over{2}}} \,\,=\,\, f_j \,-\, {{1}\over{2}} h \, u_j' \, (\alpha_j + \beta_j) \,-\, 
{{1}\over{2}} h^2 \, u_j'' \, ( \alpha_j \, \varphi_-' \,+\,   \beta_j \,\varphi_+'  ) 
\,\,+\,$ }
\setbox22=\hbox {$\displaystyle  \qquad \qquad \qquad 
+ \, {{h^2}\over{4}} \,  u_j'' \, (\alpha_j + \beta_j) \, 
+ \, {{h^2}\over{8}} \,  (u_j')^2 \, f_j'' \,+\,   {\rm O} (h^3)  \,. \,  $}
\setbox30= \vbox {\halign{#&# \cr \box21 \cr \box22    \cr   }}
\setbox31= \hbox{ $\vcenter {\box30} $}
\setbox44=\hbox{\noindent  (4.19) $\displaystyle \qquad  \left\{ \box31 \right. $}   
\smallskip \noindent   $ \box44 $ 

\smallskip \noindent 
We subtract (4.19) from (4.14)  and we also notice that

\smallskip \noindent   (4.20)  $\qquad   \displaystyle 
 (\alpha_j + \beta_j)  \,  u_j' \,\,=\,\, f' (v_j(t)) \,  u_j' \,\,=\,\, - {{\rm d}\over{{\rm
d}t}} \, u_j(t)  \, .\, $ 

\smallskip \noindent 
Finally we get 

\smallskip \noindent   (4.21)  $\qquad   \displaystyle 
\rho_j^h(t) \,\,=\,\, {\rm O}(h^2) \,$ 

\smallskip \noindent 
as desired.

\smallskip  \noindent $   \bullet $ \qquad 
We consider now the case of an extremum :  $ \displaystyle \quad 
\smash { u_j' \,\equiv \, {{\partial u}\over{\partial x}} (x_0 \,,\, t_0) } 
\,= \,  0 .\, $ 
Then $ \displaystyle \, {{\partial u}\over{\partial t}} (x_0 \,,\, t_0) \,= \,  0 \,$ due to
the conservation law (1.1). 
 We focus on the development of the ratio $\, r_j \,$  of the slopes :

\smallskip \noindent   (4.22)  $\qquad   \displaystyle 
r_j \,\,=\,\, {{  - {{1}\over{2}}  u_j''   \,+\, {\rm O} (h)  }  \over   {
{{1}\over{2}}  u_j''   \,+\, {\rm O} (h)  }  } 
 \,\,=\,\, -1  \,+\, {\rm O} (h)   \,. \,$

\smallskip \noindent  
In a similar manner, we have:

\smallskip \noindent   (4.23)  $\qquad   \displaystyle 
{{1}\over{r_{j+1}}}  \,\,=\,\, {{    {{3}\over{2}}  u_j''   \,+\, {\rm O} (h)  }  \over   {
{{1}\over{2}}  u_j''   \,+\, {\rm O} (h)  }  } 
 \,\,=\,\, 3 \,+\, {\rm O} (h)   \,. \,$

\smallskip \noindent 
Then the flux at the interface $\, (j+{{1}\over{2}}) \,$ admits the following form :

\smallskip \noindent   (4.24)  $\qquad   \displaystyle 
f_{j+{{1}\over{2}}} \,\,=\,\, f_j \,+\, {{1}\over{2}} h^2 \, u_j'' \, 
\Big(  {{1}\over{2}} \alpha_j \,  \varphi(-1) \,+\,  \beta_j \, \big( 1 - {{1}\over{2}} 
\varphi(3) \big) \Big) \,+\,    {\rm O} (h^3)  \,. \,  $

\smallskip \noindent 
The expressions   $\, r_{j-1} \,$  and $\, {{1}\over{r_j}} \,$   are computed according to
(4.23) and (4.22) respectively. Then the numerical flux at  $\, (j-{{1}\over{2}}) \,h \,$ can
be developed as

\smallskip \noindent   (4.25)  $\qquad   \displaystyle 
f_{j-{{1}\over{2}}} \,\,=\,\, f_j \,-\, {{1}\over{2}} h^2 \, u_j'' \, 
\Big(  \big( 1 - {{1}\over{2}}  \varphi(3) \big)  \alpha_j \, \,-\,   {{1}\over{2}}  
   \varphi(-1) \, \beta_j \, \Big) \,+\,    {\rm O} (h^3)  \,  \,  $

\smallskip \noindent 
and we finally get in this case :

\smallskip \noindent  $   \displaystyle 
\rho_j^h(t)  \,\,=\, \, -{{1}\over{h}} \,( f_{j+{{1}\over{2}}} -  f_{j-{{1}\over{2}}} ) \,$ 

\smallskip \noindent 
that is 

\smallskip \noindent   (4.26)  $\qquad   \displaystyle 
\rho_j^h(t) \,\,=\, \,  - {{1}\over{2}} h \,  u_j'' \, \Big( {{\varphi(-1) + \varphi(3)
}\over{2}} - 1 \Big) \, (\alpha_j - \beta_j ) \,+\,  {\rm O}(h^2) \,$ 

\smallskip \noindent 
and the proposition is established. $\hfill \square $ 

\bigskip \noindent  {\bf  Remark 3.1. }   

\noindent In a previous paper, Osher-Chakravarthy [1984] claimed that if a semi-discrete
scheme which takes the form

\smallskip \noindent   (4.27)  $\qquad   \displaystyle 
{{{\rm d}u_j}\over{{\rm d}t}} \,\,=\,\, {{1}\over{h}} \, \big( C_{j+{1\over2}} \, 
\Delta_{j+{1\over2}} \,- \, D_{j-{1\over2}} \, \Delta_{j-{1\over2}} \big) \,,\qquad j \in \Z
\,$ 

\smallskip \noindent 
with 

\smallskip \noindent   (4.28)  $\qquad   \displaystyle 
 C_{j+{1\over2}} \, \geq \, 0 \,,\quad  D_{j+{1\over2}} \, \geq \, 0 \,,\qquad j \in \Z
\,$ 

\smallskip \noindent 
is necessarily first order near a non sonic critical point  (${{\partial u}\over{\partial x}}
= 0 ,\, $ $   f'(u) \neq 0 $).  The numerical semi-discrete scheme (3.1)(3.2),
(2.27).(2,29) studied herein admits the form (4.27)(4.28). More precisely, we have, according
to (3.20).(3.22) :

\smallskip \noindent   (4.29)(a)  $\qquad   \displaystyle 
 C_{j+{1\over2}} \, \,=\,\, - \beta_{j+{1\over2}}  \, \Big( \,1 - {{1}\over{2}}
\widetilde{\varphi}_{j+1} +  {{1}\over{2}} \widetilde{\varphi}_{j} r_j  \, \Big) \,$
\smallskip \noindent   (4.29)(b)  $\qquad   \displaystyle 
D_{j+{1\over2}} \, \, \,=\,\,   \alpha_{j+{1\over2}}  \, \Big( \,
1 + {{1}\over{2}}  {{\varphi_{j+1}}\over{r_{j+1}}} - {{1}\over{2}} \varphi_j  \, \Big) \,. \,$

\smallskip \noindent 
Near a critical point, the expressions (3.27) and (3.28) vanish with $h$ and according to
(4.5), the factors $\, \Delta_{j\pm{1\over2}} \,$   are $\, {\rm O}(h^2) .\,$   Moreover the
coefficients $\,  \beta_{j+{1\over2}}\, $ and $\,  \alpha_{j-{1\over2}}\, $  are
bounded. The second order is clear. In fact the consistency property assumed in the previous
reference (formula (2.11)) is {\bf not} relevant for a critical point while the term 
$\, {{\partial u}\over{\partial x}} \,$  in (2.10) is null. The result can also be stated
intuitively : near a regular extremum, the curve looks like a parabola and its interpolation
by a constant is second order accurate in that case.

\smallskip \noindent 
Our condition (4.5)(b) which is the key for this second-order accuracy has been
independently proposed by Wu [1989].

\bigskip \noindent $   \bullet $ \qquad 
We study now the convergence of the sequence $\, u_h(x,\,t) \,$   defined in (4.1) to the
entropy solution of (1.1)(1.2).We have not proved that the scheme studied previously converges
to the entropy solution ; therefore we give known results concerning the convergence of
 $\, u_h(x,\,t) .\,$   First we have the following property :

\bigskip \noindent  {\bf Proposition 4.2. }   

\noindent 
The sequence of functions $\, u_h(\smb,\,\smb) \,$ defined for $\, (x,t) \in \, \R \times ]0,\,
+\infty [ \,$   from the line values $ u_j(t) $ according to the relation (4.1) are
bounded in each compact set in $\, {\rm L}^{\infty} \,\cap $ $ {\rm BV}(\R \times 0,\, +\infty
[) \,$  where the space BV is defined according to Volpert [1967]. This bound is independent of
$ hÊ>Ê0. \,$

\bigskip \noindent  {\bf Proof of Proposition 4.2. }   

\noindent 
In Part 3, we have proved the BV property relatively to the space variable $ x $  (Theorem 3.1)
and the  $\, {\rm L}^{\infty} $  stability (Theorem 3.2). It is now sufficient to prove the BV
stability for the time variable. Let $T$  and $ \Delta t > 0 \,$   be two fixed times. We
have :

\smallskip \noindent  $  \displaystyle 
\int_0^T \,{\rm d}t \,   \Big( \int_{-\infty}^{+\infty} {\rm d}x \mid u_h(x,\, t+\Delta t) - 
u_h(x,\, t ) \mid \Big) \,\, \leq \,\,$  
\smallskip \noindent  $  \displaystyle  \qquad \qquad \qquad \qquad  
\leq \,\,  \int_0^T \,{\rm d}t \,  \sum_{j \in \Z} \, \int_0^{\Delta t} \, {\rm d} \theta \mid
{{\partial u_j}\over{\partial t}} (t+\theta)  \mid \,$
\smallskip \noindent  $  \displaystyle  \qquad \qquad \qquad \qquad  
\leq \,\,  h \,   \sum_{j \in \Z} \, T \, {{4 C}\over{h}} \, 
\int_0^{\Delta t} \, {\rm d} \theta  \, \big( \mid u_{j+1} - u_j \mid \,+\, 
\mid u_{j} - u_{j-1} \mid  \big) (t+\theta)    \,$

\smallskip \noindent 
where the constant $C$ has been introduced in (3.7). 
 Therefore the BV property in time follows from the BV stability in space (Theorem 3.1). $
\hfill \square $

\bigskip \noindent $   \bullet $ \qquad 
From Proposition 4.2 and the compactness of the injection $\, {\rm L}^{\infty} \,\cap $ $ {\rm
BV} $ $\, \longrightarrow {\rm L}^{\infty} \, $  (see Volpert [1967])
we deduce that a subsequence of  uh  converges to a weak solution of (1.1)(1.2) but not
necessarily to the entropy solution. To prove that a subsequence of  uh converges to the
entropy solution it is sufficient to establish a discrete entropy inequality (Harten-Hyman-Lax
[1976]). Recently, Osher [1984] has proposed a sufficient discrete entropy condition to
establish the discrete entropy inequality :

\smallskip \noindent   (4.30)  $\qquad   \displaystyle 
\int_{\displaystyle  u_j}^{\displaystyle  u_{j+1}} \, \eta'(w) \, \big( \, f_{j+{1\over2}} - 
f(w) \, \big) \, {\rm d}w \,\, \leq \,\, 0 \,,\qquad j \in \Z \,.\,$ 

\smallskip \noindent 
When  $\,  f_{j+{1\over2}}  \, $    is computed according to the MUSCL method ({\it e.g.} the
relations (3.2), (2.27)-(2.29)), Osher [1985] proved that  when the flux function $\,\Phi
\,$  in (3.2) is monotone (which has been always assumed in this paper), the sequence $ u_h $
converges almost everywhere to the unique entropy solution of (1.1) (1.2) if the limiter
function $\, \varphi(\smb) \,$ satisfies the monotonicity property 2.4 and the following severe
restrictions :  the function $\, \varphi(\smb) \,$   is null for negative arguments and 
the function $\, \varphi(\smb) \,$   satisfies the following property :

\setbox21=\hbox {$\displaystyle  
0 \, \, \leq \,\, \varphi_j \,\, \leq \,\, 2 \, \max \, \Big[ 0 \,,\, \min \Big( 
{{u_j - \widetilde{u}_{j+{1\over2}}}\over{u_j - u_{j+1}}} \, , \, 
{{ \widetilde{u}_{j-{1\over2}} -u_j}\over{u_j - u_{j+1}}} \Big) \, \Big] \,$ }
\setbox22=\hbox { \qquad \qquad \qquad \qquad  \qquad \qquad if \quad 
$ u_j > u_{j+1} \,$, for each $\, j \in \Z \,   $}
\setbox30= \vbox {\halign{#&# \cr \box21 \cr \box22    \cr   }}
\setbox31= \hbox{ $\vcenter {\box30} $}
\setbox44=\hbox{\noindent  (4.31) $\displaystyle \qquad  \left\{ \box31 \right. $}   
\smallskip \noindent   $ \box44 $ 

\smallskip \noindent  
with the following definition for the mean value at $ (j+{1\over2})\, h \, : \,$ 

\smallskip \noindent   (4.32)  $\qquad   \displaystyle 
f \big(  \widetilde{u}_{j+{1\over2}} \big) \,\,=\,\, {{1}\over{u_j - u_{j+1}}} \,
\int_{\displaystyle u_{j+1}}^{\displaystyle u_j} f(w) {\rm d} w \,.\,$ 

\smallskip \noindent  
The condition (4.31) couples the interpolation at the interface values and the 
(physical) flux
function $f.$ This result  is at our knowledge the best one concerning the 
convergence of the
sequence  $u_h $ towards the entropy solution of the conservation law if one uses the MUSCL
approach. For other studies concerning this problem we refer also to 
Osher-Tadmor [1988] and
Vila [1988].

%%%%%%%%%%%%%%%%%%%%%%%%%%%%%%%%%%%%%%%%%%%%%%%%%%%%%%%%%%%%%%%%%%%%%%%%%%%%%%%%%%%%%%%%%%%%%
% \titredroite={\pecaps  Time Discretization  }
\toppagetrue  
\botpagetrue    

\bigskip   \bigskip  
\noindent {\bf 5. \quad  Time Discretization} 

\smallskip \noindent
We focus now on the discretization in time of the semi-discrete scheme (3.1). (3.2)
(associated with  (2.27).(2.29)). We assume that this semi-discrete 
scheme is Total Variation
Diminishing, {\it i.e.} that  the sufficient conditions (2.30)-(2.32) and (3.12).(3.13) are
satisfied. We consider the following three schemes in time $\, (\Delta t = t^{n+1} - t^{n}
\,$  is the time step and we denote by $ \, \lambda \,$  the ratio $\,
 {{\Delta t}\over{\Delta x}}) \,$~:

\smallskip \noindent
{\bf First Order Scheme} 

\smallskip \noindent

\smallskip \noindent   (5.1)  $\qquad   \displaystyle 
{{1}\over{\Delta t}} (u_j^{n+1} - u_j^n ) + {{1}\over{\Delta x}} \Big( \, 
\Phi \big( u_{j+{1\over2}}^{n,\,-} \,,\, u_{j+{1\over2}}^{n,\,+} \big) \,-\, 
\Phi \big( u_{j-{1\over2}}^{n,\,-} \,,\, u_{j-{1\over2}}^{n,\,+} \big)  \, \Big) \,\,=\,\, 0
\,. \,$ 

\smallskip \noindent
{\bf   Modified Euler Scheme (Heun scheme) }

\smallskip \noindent   (5.2)(a)  $\qquad   \displaystyle 
{{1}\over{\Delta t}} (\widetilde{u}_j^{n+1} - u_j^n ) + {{1}\over{\Delta x}} \Big( \, 
\Phi \big( u_{j+{1\over2}}^{n,\,-} \,,\, u_{j+{1\over2}}^{n,\,+} \big) \,-\, 
\Phi \big( u_{j-{1\over2}}^{n,\,-} \,,\, u_{j-{1\over2}}^{n,\,+} \big)  \, \Big) \,\,=\,\, 0
\,$ 
\smallskip \noindent   (5.2)(b)  $\qquad   \displaystyle 
{{1}\over{\Delta t}} (u_j^{n+1} - u_j^n ) + {{1}\over{\Delta x}} \, \Big( \,  \big( \, 
f_{j+{1\over2}}^n + \widetilde{f}_{j+{1\over2}}^{n+1} \big) \,-\, ( 
f_{j-{1\over2}}^n + \widetilde{f}_{j-{1\over2}}^{n+1} \big) \, \Big) \,\,=\,\, 0 \,  \,$ 

\smallskip \noindent
with the fluxes given by the relations 

\smallskip \noindent   (5.3)  $\qquad   \displaystyle 
f_{j+{1\over2}}^n \,\,=\,\, \Phi \big( u_{j+{1\over2}}^{n,\,-} \,,\, u_{j+{1\over2}}^{n,\,+}
\big) \,,\, \qquad 
\widetilde{f}_{j+{1\over2}}^{n+1} \,\,=\,\, \Phi \big( 
\widetilde{u}_{j+{1\over2}}^{n+1,\,-} \,,\, \widetilde{u}_{j+{1\over2}}^{n+1,\,+} \big) \,.\,$

\smallskip \noindent
{\bf    Predictor-Corrector Scheme} 

\smallskip \noindent   (5.4)(a))  $\qquad   \displaystyle 
{{2}\over{\Delta t}} \big(u_j^{n+{1\over2}} - u_j^n \big ) + {{1}\over{\Delta x}} \big( \, 
f_{j+{1\over2}}^n - f_{j-{1\over2}}^n \, \big) \,\, = \,\, 0 \,$ 
\smallskip \noindent   (5.4)(b))  $\qquad   \displaystyle 
{{1}\over{\Delta t}} (u_j^{n+1} - u_j^n ) + {{1}\over{\Delta x}} \big( \, 
f_{j+{1\over2}}^{n+{1\over2}}  - f_{j-{1\over2}}^{n+{1\over2}} \, \big) \,\,=\,\, 0 \,.\,$ 

\smallskip \noindent
As in the previous cases we have the same type of notations :

\smallskip \noindent   (5.5)  $\qquad   \displaystyle 
f_{j+{1\over2}}^n \,\,=\,\, \Phi \big( u_{j+{1\over2}}^{n,\,-} \,,\, u_{j+{1\over2}}^{n,\,+}
\big) \,,\, \qquad 
f_{j+{1\over2}}^{n+{1\over2}} \,\,=\,\, \Phi \big( 
u_{j+{1\over2}}^{n+{1\over2},\,-} \,,\,u_{j+{1\over2}}^{n+{1\over2},\,+} \big) \,.\,$

\smallskip \noindent
In the following, we study the TVD property of the scheme (5.1).(5.4) under some {\it ad hoc
}restrictions.

\bigskip \noindent  {\bf Property 5.1 }   

\noindent 
We suppose that the limiter function $\, \varphi \,$  defined in (2.9) satisfies the following
inequality :

\smallskip \noindent   (5.6)  $\qquad   \displaystyle 
\exists \, M > 0 \,,\, \,\, \forall \, r > 0 \,,\,\, 
{{\varphi(r)}\over{r}} \,\, \leq \,\, M \, .\,  $

\smallskip \noindent
From the Figure 3.1, we have clearly  $\, M \geq  1. \,$ 

\bigskip \noindent  {\bf Proposition  5.1 } 

\noindent  
The first order scheme (5.1) joined with (2.27).(2.29) is TVD under the hypotheses
(2.30).(2.32), (3.12).(3.13) and (5.6) for the limiter function if the following
Courant-Friedrichs-Lewy condition holds : 

\smallskip \noindent   (5.7)  $\qquad   \displaystyle 
{{\Delta t}\over{\Delta x}} \, \big( \mid \alpha_{j+{1\over2}}^{n} \mid \,+\, 
\mid \beta_{j+{1\over2}}^{n} \mid \big) \,\, \leq \,\, {{1}\over{1 + {M\over2}}} \,,\qquad j
\in \Z \, \, $ 

\smallskip \noindent  
with the partial gradients of the fluxes $\,\alpha_{j+{1\over2}} \, $  and 
$\,\beta_{j+{1\over2}} \, $ defined by the formulae (3.21).(3.22).

\bigskip \noindent  {\bf Proof of  Proposition  5.1 }   

\noindent   $   \bullet $ \qquad
The incremental form of the scheme (5.1) has the expression 

\setbox21=\hbox {$\displaystyle  
u_j^{n+1} \,\,=\,\, u_j \,+\, \lambda \, (-\beta_{j+{1\over2}}) \, \Big( \,1 - {{1}\over{2}}
\widetilde{\varphi}_{j+1} +  {{1}\over{2}} \widetilde{\varphi}_{j} \, r_j  \, \Big) \,
\Delta_{j+{1\over2}}\, $  }
\setbox22=\hbox {$\displaystyle  \qquad  \qquad  \qquad  \qquad 
- \,  \lambda \, \alpha_{j-{1\over2}} \,   \Big( \,
1 + {{1}\over{2}}  {{\varphi_{j}}\over{r_{j}}} - {{1}\over{2}} \varphi_{j-1}  \, \Big)\,
\Delta_{j-{1\over2}}\, $  }
\setbox30= \vbox {\halign{#&# \cr \box21 \cr \box22    \cr   }}
\setbox31= \hbox{ $\vcenter {\box30} $}
\setbox44=\hbox{\noindent  (5.8) $\displaystyle \qquad  \left\{ \box31 \right. $}   
\smallskip \noindent   $ \box44 $ 

\smallskip \noindent  
(we drop the index $n$ when there is no ambiguity). Therefore it is under the form considered by
Harten [1983] :

\smallskip \noindent   (5.9)  $\qquad   \displaystyle 
u_j^{n+1} \,\,=\,\, u_j \,+\, C_{j+{1\over2}} \, \Delta_{j+{1\over2}} \,-\, 
 D_{j-{1\over2}} \, \Delta_{j-{1\over2}} \,$ 

\smallskip \noindent  
with

\smallskip \noindent   (5.10)(a)  $\qquad   \displaystyle 
C_{j+{1\over2}} \,\,=\,\,  \lambda \, (-\beta_{j+{1\over2}}) \, \Big( \,1 - {{1}\over{2}}
\widetilde{\varphi}_{j+1} +  {{1}\over{2}} \widetilde{\varphi}_{j} r_j  \, \Big) \, $ 
\smallskip \noindent   (5.10)(b)  $\qquad   \displaystyle 
 D_{j-{1\over2}}  \,\,=\,\, \lambda \, \alpha_{j-{1\over2}} \,   \Big( \,
1 + {{1}\over{2}}  {{\varphi_{j}}\over{r_{j}}} - {{1}\over{2}} \varphi_{j-1}  \, \Big) \,.\,$ 

\smallskip \noindent  
Harten's conditions must be satisfied :
\smallskip \noindent   (5.11)(a)  $\qquad   \displaystyle 
C_{j+{1\over2}} \,\geq\, 0 \,,\quad $ (b) $  \qquad D_{j+{1\over2}} \,\geq\, 0 \,,\quad $ (c)
$  \qquad C_{j+{1\over2}} + D_{j+{1\over2}} \, \leq \, 1 \,.\, $

\smallskip  \noindent   $   \bullet $ \qquad
The inequalities (5.11)(a)(b) are a direct consequence of the Theorem 3.1. On the other hand,
we have :

\smallskip \noindent   (5.12)(a)  $\qquad   \displaystyle 
1 + {{1}\over{2}}  {{\varphi_j}\over{r_{j}}} - {{1}\over{2}} \varphi_{j-1} \,\, \leq
\,\,1 + {M\over2} \,$ 
\smallskip \noindent   (5.12)(b)  $\qquad   \displaystyle 
1 - {{1}\over{2}} \widetilde{\varphi}_{j+1} +  {{1}\over{2}} \widetilde{\varphi}_{j} \, r_j
\,\, \leq \,\,1 + {M\over2} \,$ 

\smallskip \noindent  
and the inequality (5.11)(c) is a direct consequence of the CFL condition (5.7). $\hfill
\square$ 

\bigskip \noindent  {\bf Proposition  5.2 }  (Shu-Osher [1988])

\noindent  
The modified Euler scheme (5.2).(5.3) is total variation diminishing under the same hypotheses
as in Proposition 5.1 and the same CFL restriction (5.7).

\bigskip \noindent   $   \bullet $ \qquad
The predictor-corrector scheme is more complicated to study. We begin with the following
lemma~:

\bigskip \noindent  {\bf Proposition  5.3 }   

\noindent  
	We consider a numerical scheme given by the following incremental form :

\smallskip \noindent   (5.13)  $\qquad   \displaystyle 
u_j^{n+1} \,\,=\,\, u_j \,+\,  a_{j+{3\over2}} \, \Delta_{j+{3\over2}} \,+\, 
c_{j+{1\over2}} \, \Delta_{j+{1\over2}} \,-\,  d_{j-{1\over2}} \, \Delta_{j-{1\over2}} \,-\, 
 b_{j-{3\over2}} \, \Delta_{j-{3\over2}} \,.\,  $ 

\smallskip \noindent  
This scheme is TVD if the following inequalities hold for each integer $j$ :

\smallskip \noindent   (5.14)(a)  $\qquad   \displaystyle 
 a_{j+{1\over2}} \, \geq \, 0 \,$ 
\smallskip \noindent   (5.14)(b)  $\qquad   \displaystyle 
 c_{j+{1\over2}}  - a_{j+{1\over2}} \, \geq \, 0 \,$ 
\smallskip \noindent   (5.14)(c)  $\qquad   \displaystyle 
1 -  c_{j+{1\over2}}  - d_{j+{1\over2}} \, \geq \, 0 \,$ 
\smallskip \noindent   (5.14)(d)  $\qquad   \displaystyle 
 d_{j+{1\over2}}  - e_{j+{1\over2}} \, \geq \, 0 \,$ 
\smallskip \noindent   (5.14)(e)  $\qquad   \displaystyle 
e_{j+{1\over2}} \, \geq \, 0 \, . \,$ 

\smallskip \noindent  
This proposition generalizes the one proposed by Harten [1983] and used previously in this
paper. The proof is straightforward. We omit it.

\bigskip \noindent   $   \bullet $ \qquad
We now analyze the TVD property for the predictor-corrector scheme (5.4). (5.5). This scheme
does not belong to the general class considered by Shu-Osher [1988]. In fact, we have not
proved that under the form (5.4).(5.5) the predictor-corrector scheme is total variation
diminishing but no conter-examples have been found. Nevertheless we modify the scheme
(5.4).(5.5) in the following way :

\bigskip \noindent  {\bf Definition 5.1 }   

\noindent  
For $k$ equal to  $n$  and  $n+{1\over2} \, $   let  $\, C_{j+{1\over2}}  \,$ 
 and  $\, D_{j-{1\over2}}  \,$   be defined by the relations

\smallskip \noindent   (5.15)(a)  $\qquad   \displaystyle 
C_{j+{1\over2}}^k \,\,=\,\, -\lambda \, \,{{ 
\Phi \Big( \, u_{j+{1\over2}}^{k,\,-} \,,\, u_{j+{1\over2}}^{k,\,-} \, \Big) \,-\, 
\Phi \Big( \, u_{j+{1\over2}}^{k,\,-} \,,\, u_{j-{1\over2}}^{k,\,+} \, \Big)  }\over{
u_{j+1} - u_j}} \,$
\smallskip \noindent   (5.15)(b)  $\qquad   \displaystyle 
D_{j-{1\over2}}^k \,\,=\,\,  \lambda \,\, {{ 
\Phi \Big( \, u_{j+{1\over2}}^{k,\,-} \,,\, u_{j-{1\over2}}^{k,\,+} \, \Big) \,-\, 
\Phi \Big( \, u_{j-{1\over2}}^{k,\,-} \,,\, u_{j-{1\over2}}^{k,\,+} \, \Big)  }\over{
u_{j+1} - u_j}} \,. \,$

\smallskip \noindent
We define {\bf two} versions of the {\bf modified} predictor-corrector scheme by the
incremental relations

\smallskip \noindent   (5.16)(a)  $\qquad   \displaystyle 
u_j^{n+{1\over2}} \,\,=\,\, u_j \,+\,  {1\over2} \, \Big( \, 
C_{j+{1\over2}}^k \, \Delta_{j+{1\over2}} \,-\,   D_{j-{1\over2}}^k \, \Delta_{j-{1\over2}} 
\, \Big) \,$ 
\smallskip \noindent   (5.16)(b)  $\qquad   \displaystyle 
u_j^{n+1} \,\,=\,\, u_j \,+\,  {1\over2} \, \Big( \, 
C_{j+{1\over2}}^k \, \Delta_{j+{1\over2}}^{n+{1\over2}} \,-\,   D_{j-{1\over2}}^k \,
\Delta_{j-{1\over2}}^{n+{1\over2}}  \, \Big) \, . \,$ 

\smallskip \noindent  
If $k$ is equal to $n$ (respectively $n+{1\over2}$) in relation (5.16) the predictor-corrector
scheme is explicit in time for the two steps (respectively implicit for the predictor step and
explicit for the corrector).
Notice that in equation (5.16) the index $k$ is the same in both steps (5.16)(a) and (5.15)(b)
whereas in the original predictor-corrector scheme (5.4) (5.5)  $k$  is equal to  $n$  for the
predictor and $k$ is equal to  $n+{1\over2}$ in the corrector. We have the following property :

\bigskip \noindent  {\bf  Proposition 5.4   }   

\noindent  
The predictor-corrector schemes defined in Definition 5.1 are TVD under the hypotheses of
Proposition 5.1.

\bigskip \noindent  {\bf  Proof of Proposition 5.4   }   

\noindent      $   \bullet $ \qquad
 We first develop $\, \Delta_{j+{1\over2}}^{n+{1\over2}}  \,$  and we drop the exponent $n$ for
all the variables defined at time $\, n \Delta t$   and the exponent $k$ of (5.16) associated
with the coefficients $C$ and $D$ :

\smallskip \noindent  $  \displaystyle 
\Delta_{j+{1\over2}}^{n+{1\over2}}  \,\,\equiv \,\, u_{j+1}^{n+{1\over2}} - 
 u_{j}^{n+{1\over2}}  \,$

\smallskip \noindent   (5.17)  $\quad   \displaystyle 
\Delta_{j+{1\over2}}^{n+{1\over2}}  \,\,=\,\, \Big( 1 - {1\over2} \big( C_{j+{1\over2}}
+ D_{j+{1\over2}}  \big) \Big) \, \Delta_{j+{1\over2}} \,+\,  {1\over2} C_{j+{3\over2}} \,
\Delta_{j+{3\over2}} \,+\,  {1\over2} D_{j-{1\over2}} \, \Delta_{j-{1\over2}} \,.\,$

\bigskip \noindent      $   \bullet $ \qquad
 We substitute (5.17) in equation (5.16)(b). We obtain in this manner :

\smallskip \noindent   $  \displaystyle 
u_j^{n+1} \,\,=\,\, u_j + {1\over2}\, C_{j+{1\over2}} \, C_{j+{3\over2}} \,
\Delta_{j+{3\over2}} \,+\, \,$
\smallskip \noindent   $  \displaystyle \qquad \qquad \qquad \qquad 
\,+\, \Big[ C_{j+{1\over2}} \, \big( 1 -  {1\over2} (C_{j+{1\over2}} 
  + D_{j+{1\over2}} ) \big) -  {1\over2}   D_{j-{1\over2}}  \,   C_{j+{1\over2}} 
\, \Big] \, \Delta_{j+{1\over2}} \,$
\smallskip \noindent   $  \displaystyle  
- \, \Big[ D_{j-{1\over2}} \, \big( 1 -  {1\over2} (C_{j-{1\over2}} 
  + D_{j-{1\over2}} ) \big) -  {1\over2}   D_{j-{1\over2}}  \,   C_{j+{1\over2}} 
\, \Big] \, \Delta_{j-{1\over2}} \,-\, {1\over2} D_{j-{1\over2}}  \,  D_{j-{3\over2}}  \, 
 \Delta_{j-{3\over2}}  \,.\,$

\smallskip \noindent  
This representation is of the type proposed in (5.13). The inequalities (5.14) take the form :

\smallskip \noindent   (5.18)(a)  $\qquad   \displaystyle 
C_{j+{1\over2}} \, C_{j-{1\over2}} \,\, \geq \,\, 0 \,$ 
\smallskip \noindent   (5.18)(b)  $\qquad   \displaystyle 
C_{j+{1\over2}} \, \big( 1 -  {1\over2} (C_{j+{1\over2}} 
  + D_{j+{1\over2}} ) \big) -  {1\over2} \big(   D_{j-{1\over2}} +  C_{j-{1\over2}} \big) 
C_{j+{1\over2}} \,\, \geq \,\, 0 \,$ 

\setbox21=\hbox {$\displaystyle  
( C_{j+{1\over2}}  + D_{j+{1\over2}} ) \, \Big( 1 - {1\over2} (C_{j+{1\over2}} 
  + D_{j+{1\over2}} ) \Big) \,$ }
\setbox22=\hbox {$\displaystyle  \qquad \qquad \qquad \qquad 
-  {1\over2} \Big (  C_{j+{1\over2}}  \,  D_{j-{1\over2}}  +  C_{j+{3\over2}}  \, 
D_{j+{1\over2}} \Big) \,\, \leq \,\, 1  \,   $}
\setbox30= \vbox {\halign{#&# \cr \box21 \cr \box22    \cr   }}
\setbox31= \hbox{ $\vcenter {\box30} $}
\setbox44=\hbox{\noindent  (5.18)(c)  $\displaystyle \qquad  \left\{ \box31 \right. $}   
\smallskip \noindent   $ \box44 $ 

\smallskip \noindent   (5.18)(d)  $\qquad   \displaystyle 
( C_{j+{1\over2}}  + D_{j+{1\over2}} ) \, \Big( 1 - {1\over2} (C_{j+{1\over2}} 
  + D_{j+{1\over2}} ) \Big) \,-\,  {1\over2}  D_{j-{1\over2}} \,(C_{j+{1\over2}}   +
D_{j+{1\over2}} ) \,\geq \, 0 \,$ 

\smallskip \noindent   (5.18)(e)  $\qquad   \displaystyle 
D_{j+{1\over2}} \, D_{j-{1\over2}} \,\, \geq \,\, 0 \, . \,$ 

\smallskip \noindent 
From the inequalities (5.11) the inequalities (5.18) (a) and (e) hold clearly. On the other
hand, we have :

\smallskip \noindent   left hand side of (5.18)(b)  $  \displaystyle \quad 
\geq \,\, C_{j+{1\over2}} \, \big( 1 - {1\over2} \big) \,-\,  {1\over2} C_{j+{1\over2}} \,=\,
0 \,$

\smallskip \noindent  
and (5.18)(b) is established. The inequality (5.18)(d) holds by the same argument. Finally,
the left hand side of (5.18)(c) is dominated by $\, ( C_{j+{1\over2}} + D_{j+{1\over2}}) \,$ 
and this last inequality is established according to (5.11)(c). $\hfill \square $ 

\bigskip \noindent  {\bf  Remark 5.1  }   

\noindent  
As mentioned previously, we have not proved that under the natural form (5.4).(5.5) the
predictor-corrector scheme is TVD. In fact, the inequalities of the type (5.18) contain in
that case a mixing of the $C$s and $D$s with time steps $n$ and $n+{1\over2}$. The simple
algebra derived in the last proposition is not possible anymore. We also remark that
Definition 5.1 proposes a scheme under its incremental form. Therefore this scheme is not a
priori conservative. These two remarks illustrate the difficulty in maintaining both the
conservative form of the numerical scheme and the TVD property for a predictor-corrector
scheme constructed from the first order scheme (5.1).

\bigskip \noindent   $   \bullet $ \qquad
In the last part of this section we study the order of accuracy of the predictor-corrector
scheme (5.4).(5.5) and the modified Euler scheme (5.2).(5.3) in the sense of the truncation
error near a regular local extremum of the solution $\, (t,x)\longmapsto  u(t,\,x) \,$ of the
conservation law (1.1). We recall that in Proposition 4.1 we have proved that the truncation
error (4.4) for the method of lines is second order accurate in space if the limiter function
$\, \varphi \,$  satisfies the condition 

\smallskip \noindent   (5.19)  $\qquad   \displaystyle 
\varphi(-1) + \varphi(3) \,=\, 2 \,.\,  \,$ 

\smallskip \noindent 
We first make the following two hypotheses concerning the numerical flux function $\, \Phi
\,$  of relation (3.2) and the limiter function $\, \varphi(-1) \, $  for the interface
interpolations.

\bigskip \noindent  {\bf Property 5.2 }   

\noindent  
We suppose that the monotone flux function $\, \Phi \equiv  \Phi (u,\,v) \, $ satisfies the
condition :

\smallskip \noindent   (5.20)(a)  $\qquad   \displaystyle 
{\rm d}f(w) \, > \, 0 \,\,\, \Longrightarrow \,\,\, {{\partial   \Phi }\over{\partial v}}
(w,\,w) \,\,=\,\, 0 \,$ 
\smallskip \noindent   (5.20)(b)  $\qquad   \displaystyle 
{\rm d}f(w) \, < \, 0 \,\,\, \Longrightarrow \,\,\, {{\partial   \Phi }\over{\partial u}}
(w,\,w) \,\,=\,\, 0 \, .\,$

\bigskip \noindent   $   \bullet $ \qquad
We note that the monotone flux functions proposed by Godunov [1959] and Engquist-Osher [1980]
satisfy the hypothesis 5.1 whereas the numerical flux associated with the Lax-Wendroff [1960]
scheme  :

\smallskip \noindent   (5.21)  $\qquad   \displaystyle 
\Phi_{\mu} (u,\,v) \,\,=\,\, {1\over2} \big( f(u) + \mu \, u  \big) \,+\, 
{1\over2} \big( f(v) - \mu \, v  \big) \,$ 

\smallskip \noindent  
(where $\, \mu \,$  is chosen such that $\, {\rm d}f(w) Ê\leq Ê-\mu \,$  in the domain of
variation of $w$) does {\bf not}.

\bigskip \noindent  {\bf Property 5.3 }   

\noindent  
The limiter function   $\, r \longmapsto  \varphi(r) \, $  satisfies the identity

\smallskip \noindent   (5.22)  $\quad   \displaystyle 
(1 - \sigma) \,\, \varphi \Big( {{1 + \sigma}\over{-1+ \sigma}} \Big) \,\,+\,\, 
(1 + \sigma) \,\, \varphi \Big( {{3 + \sigma}\over{1 + \sigma}} \Big)\,\,\equiv \,\, 2
\,,\qquad \forall \, \sigma \in \, [0,\, \delta ]  \,$

\smallskip \noindent 
for some parameter $\, \delta  > 0 \,$ and the discrete relations

\smallskip \noindent   (5.23)  $\quad   \displaystyle 
\big( k + {1\over2} \big) \,\, \varphi \Big( {{k - {1\over2}}\over{k + {1\over2}}} \Big) \,\,
- \,\, \big( k - {1\over2} \big) \,\, \varphi \Big( {{k - {3\over2}}\over{k - {1\over2}}}
\Big) \,\,= \,\, 1 \,,\,\,  \, k = -2,\,-1 ,\, 0 ,\, +1 ,\, +2   \,. \,$

\bigskip \noindent  {\bf Proposition 5.5.  Accuracy of the predictor-corrector scheme}   

\noindent 
	If on one hand the flux function $\, \Phi \,$  is regular ({\it i.e.} absolutely continuous
relatively to each variable, monotone (relations (3.17)) and satisfies the property 5.2, and
if on the other hand the limiter function $\, \varphi \,$  satisfies the monotonicity property
2.4, the property 5.3, and is derivable at the particular value $r=1,\,$  then the
predictor-corrector scheme (5.4).(5.5) admits a truncation error 

\smallskip \noindent   (5.24)  $\qquad   \displaystyle 
\rho(h,\,u) \, \equiv \, {{1}\over{\Delta t}} \Big( u(x_j,\, t^{n+1}) - u(x_j,\, t^{n}) \Big) 
\,+\, {{1}\over{\Delta x}} \, \Big(  f_{j+{1\over2}}^{n+{1\over2}} (u) \,-\, 
f_{j-{1\over2}}^{n+{1\over2}}(u)  \Big) \,$ 

\smallskip \noindent 
of second order accuracy :  $ \quad  
\rho(h,\,u) \,=\, {\rm O}(h^2) \, \,$ 
if the following CFL restriction is satisfied :

\smallskip \noindent   (5.25)  $\qquad   \displaystyle 
{{\Delta t}\over{\Delta x}} \,\, \sup_w \, \mid {\rm d}f(w) \mid \,\, \leq \,\, \delta
\,.\,$ 

\bigskip \noindent  {\bf Proof of Proposition 5.5. }   

\noindent 
Three cases are distinguished by the fact that the first derivative of $u$ is null or not. We
use the same notations as in Proposition 4.1.

\smallskip \noindent   $   \bullet $ \qquad $\,\, u_x \,\equiv \, {{\partial u}\over{\partial
x}} (x_j ,\, t^n) \neq 0 \,.\,$ 
We have

\smallskip \noindent   (5.26)  $\qquad   \displaystyle 
r_{j+k} \,\,\equiv \,\, {{v_{j+k}^n - v_{j+k-1}^n}\over{v_{j+k+1}^n - v_{j+k}^n}} \,\,=\,\, 1
- {{u_{xx}}\over{u_x}} \,h \,\, +\,\, {\rm O}(h^2) \,,\qquad j,\, k \, $ integers. 

\smallskip \noindent 
When we introduce this relation into (4.7) we have :

\setbox21=\hbox {$\displaystyle  
v_{j+k+{1\over2}}^{n,\,\pm} \,\,=\,\, v_j \,+\, \big( k+{1\over2} \big) \, h \, u_x 
\,+\, {1\over2} \big( k^2 + k+{1\over2} - \varphi'(1) \big) \, h^2 \, u_{xx}   \,$  }
\setbox22=\hbox {$\displaystyle  \qquad \qquad  \qquad \qquad 
- \, {1\over2} \Delta t \, f_x \,-\, {1\over2} \Delta t \,h \, \big( k+{1\over2} \big) \, 
f_{xx} \,+\, {\rm O}(h^3) \, \, $}
\setbox30= \vbox {\halign{#&# \cr \box21 \cr \box22    \cr   }}
\setbox31= \hbox{ $\vcenter {\box30} $}
\setbox44=\hbox{\noindent  (5.27)  $\displaystyle \qquad  \left\{ \box31 \right. $}   
\smallskip \noindent   $ \box44 $ 

\smallskip \noindent 
and

\smallskip \noindent   (5.28)  $\qquad   \displaystyle 
v_{j+k}^{n+{1\over2}} \,\,=\,\, v_j + k h \, u_x \, {1\over2} k^2 h^2 \, u_{xx} \,-\,
 {1\over2} \Delta t \, f_x - {1\over2} \Delta t \, k h \, f_{xx} \,+\, {\rm O}(h^3) \, .\, $

\smallskip \noindent 
So that the ratio  $\, r_{j+k}^{n+{1\over2}} \,$    at  $\, t = (n+{1\over2} ) \Delta t \,$ 
  also satisfies the relation (5.26) (!), {\it i.e.}~:

\smallskip \noindent   (5.29)  $\qquad   \displaystyle 
r_{j+k}^{n+{1\over2}} \,\,\equiv \,\, {{v_{j+k}^{n+{1\over2}} - v_{j+k-1}^{n+{1\over2}}}
\over{v_{j+k+1}^{n+{1\over2}} - v_{j+k}^{n+{1\over2}}}} \,\,=\,\, 1 - {{u_{xx}}\over{u_x}} \,h
\,\, +\,\, {\rm O}(h^2) \,.\,$

\smallskip \noindent  
We insert again this relation into (4.7). We deduce

\setbox21=\hbox {$\displaystyle  
v_{j+k+{1\over2}}^{n,\,\pm} \,\,=\,\, v_j \,+\, \big( k+{1\over2} \big) \, h \, u_x 
\,+\, {1\over2} \big( k^2 + k+{1\over2} - \varphi'(1) \big) \, h^2 \, u_{xx}   \,$  }
\setbox22=\hbox {$\displaystyle  \qquad \qquad  \qquad \qquad 
- \, {1\over2} \Delta t \, f_x \,-\, {1\over2} \Delta t \,h \, \big( k+{1\over2} \big) \, 
f_{xx} \,+\, {\rm O}(h^3) \, \, $}
\setbox30= \vbox {\halign{#&# \cr \box21 \cr \box22    \cr   }}
\setbox31= \hbox{ $\vcenter {\box30} $}
\setbox44=\hbox{\noindent  (5.30)  $\displaystyle \qquad  \left\{ \box31 \right. $}   
\smallskip \noindent   $ \box44 $ 

\smallskip \noindent 
and the numerical flux can finally be developed at the point $\, (j+k+{1\over2}) Êh \,$ as 

\smallskip \noindent   $  \displaystyle 
f_{j+k+{1\over2}}^{n+{1\over2}} \,\, = \,\, f \big( 
 v_{j+k+{1\over2}}^{n,\,\pm} \big) \,+\,
{\rm O}(h^3) \,$ 

\setbox21=\hbox {$\displaystyle  
f_{j+k+{1\over2}}^{n+{1\over2}} \,\, = \,\,  f_j \,+\, \big( k+{1\over2} \big) \, h \, f_x 
\,-\, {1\over2} \Delta t \, f' \,f_x \,+\,  \big(  k+{1\over2}  
\big)^2 \, h^2 \, f_{xx} \,$}
\setbox22=\hbox {$\displaystyle  \qquad \qquad  \qquad \qquad 
\,-\, {1\over2} \Delta t \,  \big(  k+{1\over2}  \big) \, h^2 \, \big( f' \,  
f_{xx} + f' \,
f_x \, u_x  \big) \, h^2 \, f_{xx}  \,$  }
\setbox23=\hbox {$\displaystyle  \qquad \qquad    
+ {1\over8} \Delta t^2 \, \big[ f' \big( f'' \, u_x \, f_x + f' \,
 f_{xx} \big) \,+\, f'' \,
(f_x)^2 \,\big] \,+\, {\rm O}(h^3) \,.\, $}
\setbox30= \vbox {\halign{#&# \cr \box21 \cr \box22  \cr \box23    \cr   }}
\setbox31= \hbox{ $\vcenter {\box30} $}
\setbox44=\hbox{\noindent  (5.31)  $\displaystyle \qquad  \left\{ \box31 \right. $}   
\smallskip \noindent   $ \box44 $ 

\smallskip \noindent 
Then the finite difference of fluxes at the point $\, x_j\,$  is given by : 

\smallskip \noindent   (5.32)  $\qquad   \displaystyle 
f_{j +{1\over2}}^{n+{1\over2}} \,-\,  f_{j -{1\over2}}^{n+{1\over2}} \,\,=\,\, f_x \, h
\,-\,{1\over2} \Delta t \, \, h \, \big( f' \, f_{xx} + f'' \, u_x \, f_x \big) \,+\, {\rm
O}(h^3) \,.\, $

\smallskip \noindent 
Finally the truncation error is evaluated by 

\smallskip \noindent   (5.33)  $\qquad   \displaystyle 
\rho(h,\,u) \,\,=\,\, u_t \,+\, {1\over2} \, \Delta t \, u_{tt} \,+\, 
f_x \,-\, {1\over2} \, 
 \, \Delta t \, \big( f' \, f_{xx} + f' \, u_x \, f_x \big) \,+\, {\rm O}(h^2) \,\,$ 

\smallskip \noindent  
and the result is in this particular case a direct consequence of elementary 
algebra from the
partial differential equation (1.1).

\smallskip \noindent   $   \bullet $ \qquad $\,\, u_x \,= \, 0 \,,\,\, u_{xx} \,\ne
\, 0 \,.\,$
We observe first that $\, \smash { v_j^{n+1} }\,$   is defined by the intermediate values 
$\, v_j^{n+1/2} \,$  for $\, k = -2,\, -1 ,\, 0, +1,\, +2 \,$ because the interpolation (4.7)
is a three-point scheme. Moreover these nodal values at $\, t = (n+{1\over2}) \Delta t \,$ are
defined from the given values $\, v_{j+k}^n \equiv u(x_{j+k},\,t^n) \,$  for $\, k = -3,\,
-2,\, -1 ,\,0, +1,\, +2,\, +3.\, $ We have :

\smallskip \noindent   (5.34)  $\qquad   \displaystyle 
r_{j+k}^n \,\,=\,\, {{k-{1\over2}}\over{k+{1\over2}}} \,+\, {\rm O}(h) \,$ 

\smallskip \noindent   (5.35)(a)  $\qquad   \displaystyle 
v_{j+k+{1\over2}}^{n,\,-} \,\,=\,\, v_j \,+ \, {1\over2}\, h^2 \, u_{xx} \, \Big[ \,  k^2 \,+\,
\big (k+{1\over2} \big) \, \varphi \Big( {{k-{1\over2}}\over{k+{1\over2}}} \Big) 
\, \Big] \,+\, {\rm O}(h^3 )  \,$
\smallskip \noindent   (5.35)(b)  $\qquad   \displaystyle 
v_{j+k+{1\over2}}^{n,\,+} \,\,=\,\, v_j \,+ \, {1\over2}\, h^2 \, u_{xx} \, \Big[ \,  (k+1)^2
\,-\, \big (k+{1\over2} \big) \, \varphi \Big( {{k+{3\over2}}\over{k+{1\over2}}} \Big) 
\, \Big] \,+\, {\rm O}(h^3 ) \,.\,$

\smallskip \noindent  
Therefore the flux difference $\, \smash {\delta f_{j+k}^n \,\equiv \, f_{j+k+{1\over2}}^n  - 
f_{j+k-{1\over2}}^n  } \,$  are evaluated by the developments 

\smallskip \noindent   (5.36)  $\qquad   \displaystyle 
 f_{j+k}^n \,\,=\,\, k \, h^2 \, f' \, u_{xx} \, + \, {\rm O}(h^3) \,,\qquad 
k = -2,\, -1 ,\,0, +1,\, +2  \, $

\smallskip \noindent
because the coefficients of the partial derivatives of the flux function 
$\, \Phi\,$  in Taylor
formula are null according to the property 5.3 (relations (5.23)). The values at the
intermediate time admit the development 

\smallskip \noindent   (5.37)  $\qquad   \displaystyle 
v_{j+k}^{n+{1\over2}} \,\,=\,\, v_j \,+ \, {1\over2}\, k^2 \, h^2 \, u_{xx} \,-\,
{1\over2} \, \Delta t \, k \, h^2 \, f' \, u_{xx} \,+\, {\rm O}(h^3) \,$

\smallskip \noindent  
and we have 

\smallskip \noindent   (5.38)  $\qquad   \displaystyle 
r_{j+k}^{n+{1\over2}} \,\,=\,\, {{-1 \,+\, 2 k \,-\, \lambda \, f'}\over{
 1 \,+\, 2 k \,-\, \lambda \, f'}} \, + \, 
{\rm O}(h^3) \,,\qquad  k =   -1 ,\,0, +1   \,. \, $
 
\smallskip \noindent 
Using the Taylor formula again, we finally get : 

\setbox21=\hbox {$\displaystyle  
f_{j+ {1\over2}}^{n+ {1\over2}}  \,-\,  f_{j -{1\over2}}^{n+ {1\over2}} \,\,=\,\, {{\partial
\Phi}\over{\partial u}} (v_j ,\, v_j) \,\, \Big[ \,  -1 + {1\over2} \, (1 - \lambda f')
\,\,   \varphi \Big( {{1+ \lambda f'}\over{-1+ \lambda f'}} \Big) \,+\, \, $  }
\setbox22=\hbox {$\displaystyle   \qquad  \qquad  \qquad  \qquad \qquad \qquad 
\,+\, {1\over2}\, (1 + \lambda f')
\,\, \varphi \Big( {{3 + \lambda f'}\over{1+ \lambda f'}} \Big) \, \Big] \,+\,  $  }
\setbox23=\hbox {$\displaystyle  \qquad    \qquad   
\,+\,  {{\partial \Phi}\over{\partial v}} (v_j ,\, v_j) \,\, \Big[   \,   1 - {1\over2} \, (1 -
\lambda f') \,\,   \varphi \Big( {{ 3 -  \lambda f'}\over{1 -  \lambda f'}} \Big) 
\,+\, \, $  }
\setbox24=\hbox {$\displaystyle  \qquad  \qquad  \qquad  \qquad \qquad
\,+\, {1\over2}\, (-1 -  \lambda f')
\,\, \varphi \Big( {{ -1 +\lambda f'}\over{1+ \lambda f'}} \Big) \, \Big]  \, + \, 
{\rm O}(h^3) \,.\, $  }
\setbox30= \vbox {\halign{#&# \cr \box21 \cr \box22  \cr \box23   \cr  \box24   \cr  }}
\setbox31= \hbox{ $\vcenter {\box30} $}
\setbox44=\hbox{\noindent  (5.39) $\displaystyle \qquad  \left\{ \box31 \right. $}   
\smallskip \noindent   $ \box44 $ 

\smallskip \noindent 
To fix the ideas, we suppose $ f' > 0.\,$  Then $\, {{\partial \Phi}\over{\partial v}} \,$ 
 is null according to Property 5.2 of the numerical flux function. Moreover, the coefficient
in front of  $\, {{\partial \Phi}\over{\partial u}} \,$ 
is also null thanks to Property 5.3 (relation (5.22)) if we have $\, \lambda \, f' \leq
\delta \,$  which is exactly the CFL condition (5.25). The proof is analogous when 
$ f' < 0.\,$    The sonic case $\, f'Ê=Ê0 \,$ is very simple because the brackets of (5.39) are
both null according to (5.19) (which can be seen as a particular case of both (5.22) and (5.23)
with $\, \sigma = 0 \,$   $\, k = 0\,$  respectively). The truncation error (5.24) is again
evaluated according to relation (5.33) and the property is established in this second case.

\bigskip \noindent   $   \bullet $ \qquad $\,\, u_x \,= \, 0 \,,\,\, u_{xx} \,=
\, 0 . \, \,$  We have clearly 

\smallskip \noindent   $  \displaystyle 
  v_{j+k+{1\over2}}^{n,\,\pm} \,\,=\,  v_j \,+\, {\rm O}(h^3)  \,\, \,$  

\smallskip \noindent   
because the function $\, \varphi(\smb) \,$   is bounded and  $\, \Delta_{j+{1\over2}}
\,$   is O($h^3$). It is straightforward to develop the truncation error (5.24) and to
verify the property. This ends the proof.  $\hfill \square$

\bigskip \noindent  {\bf  Remark 5.2  }   

\noindent  
If the limiter function satisfies Property 5.3 and is affine, it takes necessarily the form 

\smallskip \noindent   (5.40)  $\qquad   \displaystyle 
\varphi(r) \,\,=\,\, 1 \,+\, a \, (r-1) \,.\,$ 

\smallskip \noindent 
This kind of function is not compatible with the monotony-convexity conditions (2.30)-(2.32)
and the sufficient TVD conditions (3.12).(3.13) if (5.40) holds for all the real values of the
argument r. Therefore the actually known limiters ({\it e.g.} minmod, superbee (see Roe
[1985]), MUSCL (see Van Leer [1977]), {\it etc.} ) do not satisfy both (5.22).(5.23).
Nevertheless , the function defined as follows

\smallskip \noindent   (5.41)(a)  $\qquad   \displaystyle 
\varphi(r) \,\,=\,\, 0 \,,\qquad \qquad \,\,\,  r \leq -3 \,$ 

\smallskip \noindent   (5.41)(b)  $\qquad   \displaystyle 
\varphi(r) \,\,=\,\, {1\over4} \, (r+3)  \,,\quad -3 \leq  r \leq -1 \,$ 

\smallskip \noindent   (5.41)(c)  $\qquad   \displaystyle 
\varphi(r) \,\,=\,\, -{{r}\over2}    \,,\qquad  \quad \,\,  -1 \leq  r \leq 0 \,$ 

\smallskip \noindent   (5.41)(d)  $\qquad   \displaystyle 
\varphi(r) \,\,=\,\,  {{5 \, r}\over2}    \,,\qquad \quad \,\,\,\, 0 \leq  r \leq {1\over3}
\,$ 

\smallskip \noindent   (5.41)(e)  $\qquad   \displaystyle 
\varphi(r) \,\,=\,\, {1\over4} \, (r+3)  \,,\quad   \,\, {1\over3} \leq  r \leq 3 \,$ 

\smallskip \noindent   (5.41)(f)  $\qquad   \displaystyle 
\varphi(r) \,\,=\,\,  {3\over2}   \,,\qquad  \qquad  \, \,\,  r \geq 3 \,$ 

\smallskip \noindent  
is affine and satisfies (5.40) for $\, r \in [ -3 ,\, -1] \, \cup \,  [{1\over3}, \, 3]. \,$
Then (5.23) is satisfied. Moreover (5.22) holds for $\, \delta = {1\over2} \,$   and
Property (5.1) is also satisfied with $\, M = {5\over2}.\,$  This limiter satisfies also the
TVD conditions (3.12).(3.13) of Property 3.1 with $\, \alpha = {3\over2} \,$ and the
properties 2.4 and 2.5 of monotonocity and convexity. Finally we observe that the interpolate
value 

\smallskip \noindent  $  \displaystyle 
u_{j+{1\over2}}^{-} \,\,=\,\, u_j \,+\, {1\over2} \, \varphi(r_j) \, (u_{j+1} - u_j) \,\,
\equiv \,\, L( u_{j-1} \,,\, u_j \,,\, u_{j+1} ) \,$ 

\smallskip \noindent   
(according to the notation proposed in (2.1)) corresponds to the Lagrange polynomial
interpolation of degree 2 if and only if $\, \varphi \,$  satisfies (5.40) with  $\, a =
{1\over4}. \,$ In the following we will denote this new limiter function (5.41)
 by ``Lagrange
limiter''. We note also that the idea of using Lagrange interpolation for the 
construction of
limiter functions was at our knowledge first proposed by Leonard [1979] with the so-called
Quick limiter. For a modification of this limiter that does not satisfy the monotonicity
condition for negative values of the argument $r$, we refer to Cahouet-Coquel [1989].

\bigskip 
\centerline { \epsfysize=6cm    \epsfbox  {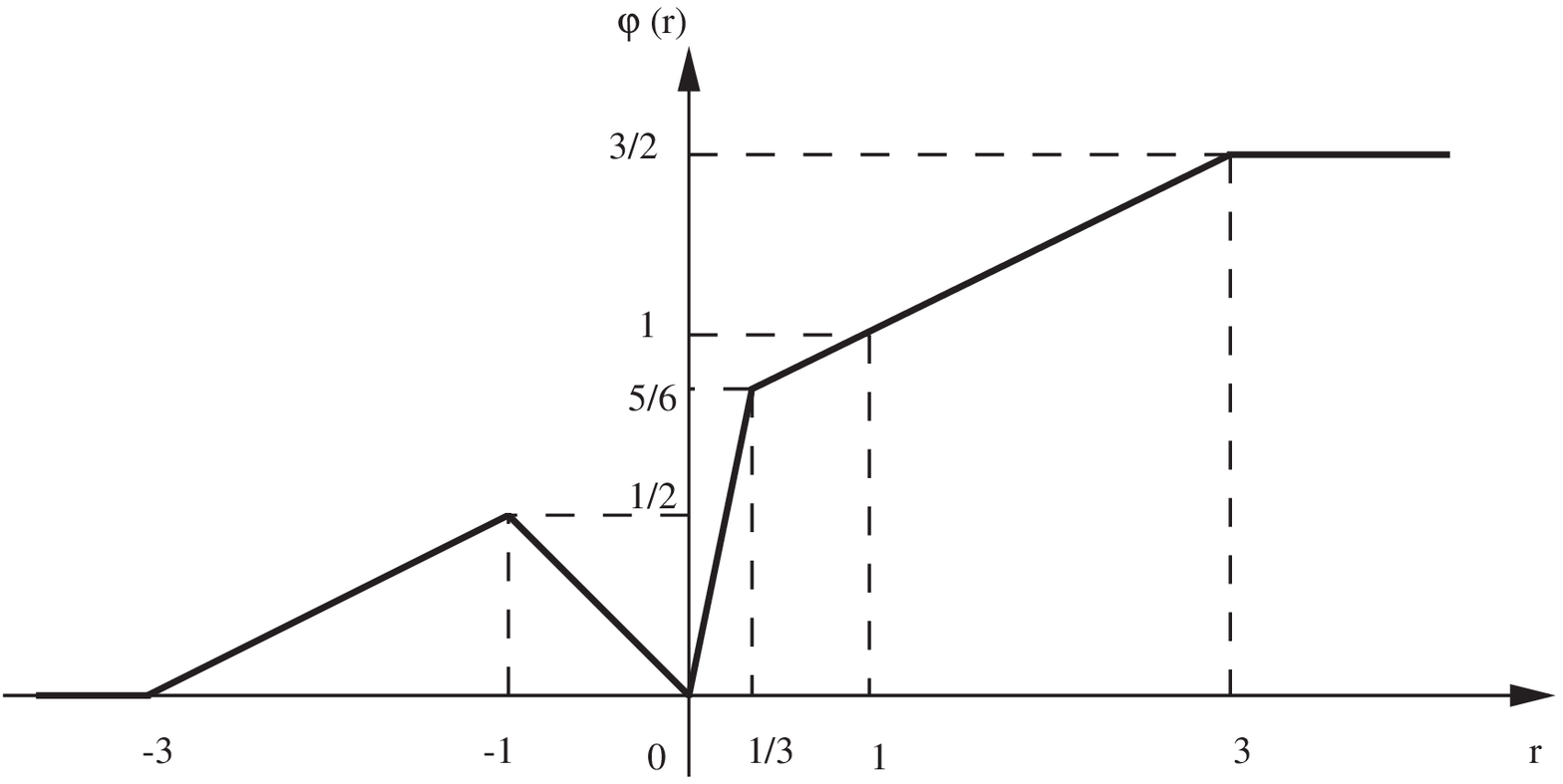} }
\smallskip  \smallskip
\noindent  {\bf Figure 5.1.} \quad  The Lagrange limiter (defined by the relations (5.41)).

\bigskip 
% \vfill \eject   %%%%%%%%%   fd 15 juin 2010
\noindent  {\bf  Proposition 5.6.	Accuracy of the Heun scheme}

\noindent
We assume that the flux function $\, \Phi \,$  and the limiter function 
$\, \varphi \,$ 
satisfy the same hypotheses as in Proposition 5.5. Then the truncation error 

\setbox21=\hbox {$\displaystyle  
\rho(h,\,u) \,\,\equiv \,\, {{1}\over{\Delta t}} \, \big( u(x_j ,\, t^{n+1}) - 
u(x_j ,\, t^{n}) \big) \,+\, $  }
\setbox22=\hbox {$\displaystyle   \qquad   \qquad  
\,+\,   {{1}\over{2 \, \Delta x}} \,\,  \Big( \, 
f_{j+{1\over2}}^n (u)  + \widetilde{f}_{j+{1\over2}}^{n+1} (u) - 
f_{j-{1\over2}}^n (u)   - \widetilde{f}_{j-{1\over2}}^{n+1} (u)  \, \Big)  \,  $  }
\setbox30= \vbox {\halign{#&# \cr \box21 \cr \box22  \cr   }}
\setbox31= \hbox{ $\vcenter {\box30} $}
\setbox44=\hbox{\noindent  (5.42) $\displaystyle \qquad  \left\{ \box31 \right. $}   
\smallskip \noindent   $ \box44 $ 

\smallskip \noindent 
of the modified Euler scheme (5.2).(5.3) is of {\bf second order accuracy} if the following
more restrictive CFL condition is satisfied :

\smallskip \noindent   (5.43)  $\qquad   \displaystyle 
{{\Delta t}\over{\Delta x}} \, \sup_{w} \, \mid f'(w) \mid \,\, \leq \,\, {1\over2} \,
\delta \,.  \,$ 

\bigskip \noindent  {\bf Proof of Proposition 5.6}  

\noindent
As in Proposition 5.5 we divide the proof in three steps, following the value of $\, u_x \,$
and $\, u_{xx} .\,$ 

\smallskip  \noindent   $   \bullet $ \qquad $\, u_x (x_j \,,\, t^n) \, \neq 0 .\,$  
We develop the numerical flux at time  $\, t^n \,$   and  $\, x = (k+j+{1\over2}) \,h \,$  as
in Proposition 5.5 :

\setbox21=\hbox {$\displaystyle  
f_{k+j+{1\over2}}^n \,\,=\,\, f_j \,+\, \big( k+{1\over2} \big) \, h \, u_x \, f' \,+\, 
{1\over2} \,  \big( k+{1\over2} \big)^2  \, h^2 \, (u_x)^2 \, f''  \,+\,   $  }
\setbox22=\hbox {$\displaystyle   \qquad   \qquad  \qquad  
+\, {1\over2} \,  h^2 \, u_{xx} \, \big( k^2 + k + {1\over2} - \varphi'(1) \big) \, f' \,+\, 
{\rm O}(h^3) \,.\,$   }
\setbox30= \vbox {\halign{#&# \cr \box21 \cr \box22  \cr   }}
\setbox31= \hbox{ $\vcenter {\box30} $}
\setbox44=\hbox{\noindent  (5.44) $\displaystyle \qquad  \left\{ \box31 \right. $}   
\smallskip \noindent   $ \box44 $ 

\smallskip \noindent 
On one hand, we have 

\smallskip \noindent   (5.45)  $\qquad   \displaystyle 
\delta f_j^n \,\,=\,\, h \, f_x \,+\, {\rm O}(h^3) \,$ 

\smallskip \noindent 
and on the other hand :

\smallskip \noindent   (5.46)  $\qquad   \displaystyle 
\widetilde{u}_{j+k}^{n+1} \,\,=\,\, u_j \,+\, (k \, h \, u_x \,-\, \Delta t \, f_x ) \,+\,
{1\over2} \, h \, \big ( k^2 \,h \, u_{xx} \,-\, 2 k \, \Delta t \, f_{xx} \, \big)  \,+\, {\rm
O}(h^3) \,. \,$ 

\smallskip \noindent 
Therefore the ratio of gradients after the ``tilda'' step (5.2)(a) is given by (5.26) or
(5.29). The extrapolated values at the interface  $\, (j+k+{1\over2}) \,$ are again evaluated
according to (5.30) except that $\, \Delta t \,$  must be replaced by $\, 2 \, \Delta t .\,$   In an
analogous way, the tilda flux defined in (5.3) admits a development (5.31) with $\,  \Delta t
.\,$  replaced by  $\, 2 \, \Delta t .\,$ We deduce :

\smallskip \noindent   (5.47)  $\qquad   \displaystyle 
\widetilde{\delta f}_j \,\, =\,\, h \, f_x \,-\, \Delta t \, h \, \big( f' \, f_{xx} + f'' \,
u_x \, f_x  \big) \, \,+\, {\rm O}(h^3) \,. \,$ 

\smallskip \noindent  
The truncation error admits exactly the form (5.33) and the property is proved in this case.

\smallskip \noindent   $   \bullet $ \qquad $\,\, u_x (x_j ,\, t^n) = 0 \,,\,  u_{xx} 
(x_j ,\, t^n) \neq 0 .\,$ 
The difference of fluxes around the point $\,  (j+k) h \,$  is given according to relation
(5.36) due to Property 5.3. We have :

\smallskip \noindent   (5.48)  $\qquad   \displaystyle 
\widetilde{r}_{j+k} \,\, \equiv \,\, {{ \widetilde{v}_{j+k}^{n+1}  \,-\, 
 \widetilde{v}_{j+k-1}^{n+1}  } \over{ \widetilde{v}_{j+k+1}^{n+1}  \,-\, 
\widetilde{v}_{j+k}^{n+1} }} \,\,=\,\, {{ -1 \,+\, 2 \,k \,-\, 2 \, \lambda \, f'}\over{
1 \,+\, 2 \,k \,-\, 2 \, \lambda \, f'}} \, + \, {\rm O}(h) \,\,$ 

\smallskip \noindent 
according to (5.43) and (5.22).  Therefore the development of the flux difference at point 
$\, j \, h \,$    and at time  $\, (n+1) \Delta t \,$   is simply developed :

\smallskip \noindent   (5.49)  $\qquad   \displaystyle 
\widetilde{\delta f}_j^{,+1} \,\, \equiv \,\, \widetilde{f}_{j+{1\over2}}^{n+1} \,-\, 
\widetilde{f}_{j-{1\over2}}^{n+1} \,\,=\,\, f' \, f_{xx} \,+\,  \, {\rm O}(h^3) \,\,$

\smallskip \noindent  
and the truncation error admits finally the form (5.33). $ \hfill \square $

%%%%%%%%%%%%%%%%%%%%%%%%%%%%%%%%%%%%%%%%%%%%%%%%%%%%%%%%%%%%%%%%%%%%%%%%%%%%%%%%%%%%%%%%%%%%%%%%
% \titredroite={\pecaps  Numerical Tests with the Advection Equation  }
\toppagetrue  
\botpagetrue

\bigskip   \bigskip  
\noindent {\bf 6. \quad  Numerical Tests with the Advection Equation} 

\smallskip \noindent
In this section we present numerical comparisons of the schemes studied
 in Part 5 for a simple model problem. We show that the truncation errors studied at
propositions 5.5 and 5.6 are numerically second order accurate for the Lagrange limiter (5.41)
especially at a nonsonic extremum. We study also the effect of the restrictive CFL condition
(5.43) on the accuracy of the two step Heun scheme. Moreover we compare the accuracy measured
in the discrete $\, \ell^1  ,\,$ $\, \ell^2  ,\,$ $\, \ell^\infty  ,\,$ norms and in terms of
the truncation error (see also Harten et al [1987]).

\smallskip \noindent
	We focus on the one-dimensional linear advection equation

\smallskip \noindent   (6.1)  $\qquad   \displaystyle 
{{\partial u}\over{\partial  t}} \,+\, a \, {{\partial u}\over{\partial  x}} \,\,=\,\, 0
\,\,,\qquad 0 \, \leq \, x  \, \leq \, 1 \,,\quad t \, \geq \, 0 \,,\quad a = 1 \,$ 

\smallskip \noindent
with periodic boundary conditions :

\smallskip \noindent   (6.2)  $\qquad   \displaystyle 
u(0) \,=\, u(1) \,.\,$ 
 
\smallskip \noindent
We look at the advection of a regular profile with a (local) extremum at $\, x \,=\,
{1\over2} \,$   without annulation of the third derivative at this point (which is the
case with a sinus-type profile) giving unexpected superconvergence for the truncation error
at  $\, x \,=\, {1\over2} \,$ :

\smallskip \noindent   (6.3)  $\qquad   \displaystyle 
{{{\rm d} u^0}\over{{\rm d}  x}} \big({1\over2} \big) \,\,=\,\, 0 \,,\qquad 
{{{\rm d}^3 u^0}\over{{\rm d}  x^3}} \big({1\over2} \big) \,\,\ne\,\, 0 \,.\, $ 

\smallskip \noindent  
We define our profile $\, u^0 \, $  as a $\, {\cal C}^2 \,$   class periodic function on the
interval  $\, [0,\, 1] \,$  according to the conditions

\smallskip \noindent   (6.4)(a)   $\qquad   \displaystyle 
u^0(x) \,\,=\,\, 1 \,-\, \big( x - {1\over2} \big)^2 \,+\, \big( x - {1\over2} \big)^3
\,,\qquad {1\over4} \,\leq \, x \, \leq \, {3\over4} \,$

\smallskip \noindent   (6.4)(b)   $\qquad   \displaystyle 
u^0(0) \,=\, 
{{{\rm d} u^0}\over{{\rm d}  x}} (0) \,=\,
{{{\rm d}^2 u^0}\over{{\rm d}  x^2}} (0) \,=\,
u^0(1) \,=\, 
{{{\rm d} u^0}\over{{\rm d}  x}} (1) \,=\,
{{{\rm d}^2 u^0}\over{{\rm d}  x^2}} (1) \,=\, 0 \,$

\smallskip \noindent   (6.4)(c)   $\qquad   \displaystyle 
u^0 \,$ polynomial of degree $\, 5 \,$ in $\displaystyle  \,\, \big[ 0 ,\, {1\over4} \big] \,$
and in   $\displaystyle  \,\, \big[ {3\over4} ,\, 1 \big] \, . \,$ 

\smallskip \noindent   
The uniqueness of$\, u^0 \,$ defined by (6.4) is clear due to the use of the classical finite
element interpolation in the intervals $ \, [ 0 ,\, {1\over4} ] \,$ and $ \, [ {3\over4} ,\,
1 ] . \,$ 

\smallskip  \noindent   $   \bullet $ \qquad
	We study the Cauchy problem (6.1)(6.2) associated with the initial condition 

\smallskip \noindent   (6.5)  $\qquad   \displaystyle 
u(0,\, x) \,\,=\,\, u^0(x) \,,\qquad   0 \, \leq \, x  \, \leq \, 1 \,$ 

\smallskip \noindent 
during one period, {\it i.e.} for time $t$ in $[0,\,1]$. At time $t = 1,$ the exact solution of
(6.1).(6.2).(6.5) is equal to the initial condition $\, u^0 \,$   and the error is therefore
easy to evaluate. We discretize the equation (6.1) on meshes with cell diameters   $\,
\Delta x \, \equiv \, h \,$    that  are constant on each mesh and chosen as powers of
$ {1\over2} \, : $ 

\smallskip \noindent   (6.6)  $\qquad   \displaystyle 
h \,\,=\,\, 2 ^{-k} \,\qquad k \,=\, 3, \, 4, \, 5, \, 6, \, 7, \, 8, \, 9,\, 10 \,,\, $   

\smallskip \noindent 
and the computational nodes $\, x_j \,$    (the cell centers in the finite volume approach)
are located at 

\smallskip \noindent   (6.7)  $\qquad   \displaystyle 
x_j \,\,=\,\, (j-1) \,h \,,\qquad j \,=\, 1 ,\, \dots ,\, N(h) \,=\, {{1}\over{h}} \,.\,$ 

\smallskip \noindent  	
We have tested three types of temporal schemes : the so-called ``five points TVD scheme'', the
predictor-corrector scheme, and the Heun scheme (second order Runge-Kutta). All these schemes
are parameterized by the interpolation $\, u_{j-1/2}^- \,$   (at the left side of the
interface  $\, x_{j-1/2}  $) from the cell center values (6.7) in the spirit described in Part
2 of this paper and detailed in the following. The temporal schemes are associated with the
Courant number 

\smallskip \noindent   (6.8)  $\qquad   \displaystyle 
\sigma \,\,=\,\, a \, {{\Delta t}\over{h}} \,.\,$ 

\smallskip \noindent
To be precise we give the formulae for this temporal integration.

\noindent  $ \smb $ 
Five points TVD scheme (see {\it e. g.} Sweby [1984]) :

\smallskip \noindent   (6.9)  $\qquad   \displaystyle 
u_j^{n+1} \,\,=\,\, u_j \,-\, \sigma \, \Big[ \, \sigma \, (u_j - u_{j-1}) \,+\, (1-\sigma) \, 
\big(  u_{j+1/2}^- \,-\,  u_{j-1/2}^- \big) \, \Big] \,.\, $ 

\smallskip \noindent  $ \smb $ 
 Predictor-Corrector Scheme (see also (5.4).(5.5) ) :

\smallskip \noindent   (6.10)(a) $\qquad   \displaystyle 
u_j^{n+1/2} \,\,=\,\, u_j \,-\, {{\sigma}\over{2}}  \, \big(  u_{j+1/2}^- \,-\,  u_{j-1/2}^-
\big)  \, $ 

\smallskip \noindent   (6.10)(b) $\qquad   \displaystyle 
u_j^{n+1} \,\,=\,\, u_j \,-\,  \sigma  \, \Big(  u_{j+1/2}^{n+1/2,\,-} \,-\, 
u_{j-1/2}^{n+1/2,\,-} \Big)  \,. \, $ 

\smallskip \noindent  $ \smb $ 
 Heun scheme (see also (5.2).(5.3))

\smallskip \noindent   (6.11)(a) $\qquad   \displaystyle 
{\widetilde u}_j^{n+1} \,\,=\,\, u_j \,-\,  \sigma   \, \big(  u_{j+1/2}^- \,-\, 
u_{j-1/2}^- \big)  \, $ 

\smallskip \noindent   (6.11)(b) $\qquad   \displaystyle 
u_j^{n+1} \,\,=\,\, u_j \,-\,  {{\sigma}\over{2}}  \bigg[ \,  \Big(  u_{j+1/2}^- \,+\, 
{\widetilde u}_{j+1/2}^{n+1,\,-} \Big) \,-\, 
\Big(  u_{j-1/2}^- \,+\,  {\widetilde u}_{j-1/2}^{n+1,\,-} \Big) \, \bigg] \,.\,$ 

\smallskip \noindent  
The interpolation at point  $\,  x_{j+1/2}\,$  is chosen according to a three point
interpolation as in the previous sections :

\smallskip \noindent   (6.12)  $\qquad   \displaystyle 
 u_{j+1/2}^- \,\,=\,\, u_j \,+ \, {1\over2} \, \varphi \bigg( {{u_j - u_{j-1} }\over {
u_{j+1}  - u_{j}  }}  \bigg) \, \big( u_{j+1}  - u_{j}  \big)  \,$

\smallskip \noindent  
or with the UNO2 interpolation proposed by Harten-Osher [1987]. Concerning the limiter
functions we have restricted ourselves to the following ones : 

\smallskip \noindent    {\it (i)} \quad First-order upwind scheme 

\smallskip \noindent   (6.13)  $\qquad   \displaystyle 
\varphi^{\rm upwind}(r) \,\,\equiv  \,\, 0 \,$

\smallskip \noindent    {\it (ii)} \quad  Lax Wendroff type scheme

\smallskip \noindent   (6.14)  $\qquad   \displaystyle 
\varphi^{\rm LW}(r) \,\,\equiv  \,\, 1 \, . \,$ 

\smallskip \noindent  
We remark that (6.14) joined with (6.12) defines the Lax-Wendroff [1960] scheme 
only when it is associated with the temporal scheme (6.9). The two other cases justify the
``type'' restriction in our denomination.

\smallskip \noindent     {\it (iii)} \quad 
Van Leer's MUSCL limiter [1977]

\smallskip \noindent   (6.15)(a)  $\quad   \displaystyle 
\varphi^{\rm MUSCL}(r) \,\,\equiv  \,\, 0 \,\,\,\, {\rm if } \,\,\,\, r \,< \, 0 \,,\qquad 
\varphi^{\rm MUSCL}(r) \,\,=  \,\, 2 \, r  \,\,\,\, {\rm if } \,\,\,\, 0 \, \leq \, r \, \leq
\, {1\over3} \,$

\smallskip \noindent   (6.15)(b)  $\quad   \displaystyle 
\varphi^{\rm MUSCL}(r) \,=\, {1\over2} \, (1 + r)  \,\,\,\, {\rm if } \,\,\,\, 
{1\over3} \, \leq \, r \, \leq \,  3  \,,\quad 
\varphi^{\rm MUSCL}(r) \,=\, 2   \,\,\,\, {\rm if } \,\,   r \, \geq \, 3
\,.\,$

\smallskip \noindent     {\it (iv)} \quad 
Min-mod limiter (Harten [1983])

\setbox21=\hbox {$\displaystyle  
\varphi^{\rm minmod}(r) \,\, \equiv  \,\,  0  \,\,\,\, {\rm if } \,\,\,\, r \,< \, 0 \,,\, $  }
\setbox22=\hbox {$\displaystyle  
\varphi^{\rm minmod}(r) \,\,=\,\, r  \,\,\,\, {\rm if } \,\,\, \, 0 \, \leq \, r \, \leq
\, 1 \,, \,   $  }
\setbox23=\hbox {$\displaystyle  
\varphi^{\rm minmod}(r) \,\equiv \, 1  \,\,\,\, {\rm if } \,\,\,   r \, \geq \, 1  \,$ }
\setbox30= \vbox {\halign{#&# \cr \box21 \cr \box22  \cr \box23    \cr   }}
\setbox31= \hbox{ $\vcenter {\box30} $}
\setbox44=\hbox{\noindent  (6.16) $\displaystyle \qquad  \left\{ \box31 \right. $}   
\smallskip \noindent   $ \box44 $

\smallskip \noindent  {\it (v)} \quad 
Superbee limiter (Roe [1985])

\smallskip \noindent   (6.17)(a)  $\qquad   \displaystyle 
\varphi^{\rm superB}(r) \,\,\equiv\,\, 0   \,\,\,\, {\rm if } \,\,\,\, r \,< \, 0 \,,\,\,
\quad 
\varphi^{\rm superB}(r) \,\,=\,\, 2 \, r     \,\,\,\, {\rm if } \,\,\,\, 0 \,\leq \, r \, \leq
\,  {1\over2} \,, \,$

\smallskip \noindent   (6.17)(b)  $\qquad   \displaystyle 
\varphi^{\rm superB}(r) \,\, \equiv \,\, 1  \,\,\,\, {\rm if } \,\,\,\, {1\over2} \,\leq \, r
\,\leq \, 1   \,,\,\, $

\smallskip \noindent   (6.17)(c)  $\qquad   \displaystyle 
\varphi^{\rm superB}(r) \,\, = \,\, r  \,\,\,\, {\rm if } \,\,\,\, 1 \,\leq \, r \,\leq \, 2 
\,,\,\, \quad \, 
\varphi^{\rm superB}(r) \,\, \equiv \,\, 2  \,\,\,\, {\rm if } \,\,\,\,   r \,\geq \, 2 \,.
\,$

\smallskip \noindent  {\it (vi)} \quad 
The Lagrange limiter proposed in Part 5 (relations (5.41)).

\smallskip \noindent   {\it (vii)} \quad 
	Min Mod absolute value limiter, that we define according to 

\smallskip \noindent   (6.18)  $\qquad   \displaystyle 
\varphi^{\rm MMA}(r) \,\, = \,\, \varphi^{\rm minmod} \big( \mid r\mid \big)  \,. \,$ 

\smallskip \noindent   {\it (viii)} \quad 
The UNO2 interpolation, which is a {\bf five} point interpolation that is defined according to 

\smallskip \noindent   (6.19)  $\quad   \displaystyle 
{\rm minmod} (a ,\, b ) \,\, \equiv \,\, a \,  \varphi^{\rm minmod} \big( {{b}\over{a}}  \big)
\,$ 

\smallskip \noindent   (6.20)  $\quad   \displaystyle 
\partial^2 u_{j+1/2} \,\,=\,\, {\rm minmod} \, \big(  \, u_{j-1} -2 \, u_j + u_{j+1} \,,\, 
  u_{j} -2 \, u_{j+1}  + u_{j+2} \, \big) \,$

\smallskip \noindent   (6.21)  $\quad   \displaystyle 
u_{j+1/2}^-  \,=\, u_j \,+\, {1\over2} \, {\rm minmod} \, \big( u_j - u_{j-1} + {1\over2}
\, \partial^2 u_{j-1/2}  \,,\,  u_{j+1} - u_j - {1\over2} \, \partial^2 u_{j+1/2} \, \big) 
\,. \, $

\smallskip \noindent 
We note that this UNO2 scheme has a stencil of a priori 7 points for the one step 
temporal scheme (6.9) (only 5 points in the interpolations (i)-(vii))  and 13 points for the
temporal schemes (6.10).(6.11)  (compared to 9 with the other schemes described herein).

\smallskip \noindent 
	Known results established in Sweby [1984] are precised by the following proposition~:

\bigskip \noindent  {\bf Proposition  6.1 } 

\noindent  
The discretization of the scalar advection equation (6.1) joined with the five
point scheme (6.9).(6.12) and some limiter function $  \, \varphi \,$ admits a truncation
error

\setbox21=\hbox {$\displaystyle  \quad 
{{1}\over{\Delta t}} \, \big( u(x_j ,\, t^{n+1}) \,-\,  u(x_j ,\,
t^{n}) \big) \,+\, \,$   }
\setbox22=\hbox {$\displaystyle    
\,+\, {{1}\over{h}} \, \big( {{a \, \Delta t}\over{\Delta x}} \big) ^2 \, 
 \big( u(x_j ,\, t^{n}) \,-\,  \big( u(x_{j-1} ,\, t^{n+1}) \big) \,$   }
\setbox23=\hbox {$\displaystyle     
-\, {{1}\over{h}} \, \Big( {{a \, \Delta t}\over{\Delta x}} \Big) \,
 \Big( 1 -  {{a \, \Delta t}\over{\Delta x}} \Big) \, \big( v_{j+1/2}^- (u,\, t^n) \, -\, 
 v_{j-1/2}^- (u,\, t^n)  \,$ }
\setbox30= \vbox {\halign{#&# \cr \box21 \cr \box22  \cr \box23    \cr   }}
\setbox31= \hbox{ $\vcenter {\box30} $}
\setbox44=\hbox{\noindent  (6.22) $\displaystyle  \quad 
\rho(h,\,u) \,\,\equiv \,\, \left\{ \box31 \right. $}   
\smallskip \noindent   $ \box44 $ 

\smallskip \noindent 
(where  $\,  v_{j+1/2}^- (u,\, t^n) \,$   is a notation for interpolation (6.12) applied
with the values of the solution interpolated at the grid points $\, x_{j+1/2}$)  which is
second order accurate 

\smallskip \noindent   (6.23)  $\quad   \displaystyle 
\rho(h,\,u) \,\,=\,\, {\rm O}(h^2) \,$ 

\smallskip \noindent 
at a local extremum of the function $\, u \,$   :

\smallskip \noindent   (6.24)  $\quad   \displaystyle 
{{\partial u}\over{\partial x}} \big (x=x_j \,,\, t=t^n \big) \,\, = \,\, 0 \,,\qquad \forall
\, h \,$

\smallskip \noindent  
if the following condition (proposed in Part  4) holds :

\smallskip \noindent   (6.25)  $\quad   \displaystyle 
\varphi(-1) \,+\, \varphi(3) \,\,=\,\, 2 \,.\,$

\bigskip \noindent  {\bf Proof of  Proposition  6.1 }   

\noindent   $   \bullet $ \qquad
From the relations (5.35)(a) we have : 

\smallskip \noindent   (6.26)  $\quad   \displaystyle 
v_{j+1/2}^- \,-\, v_{j-1/2}^- \,=\, {1\over2} \, h^2 \, {{\partial^2 u}\over{\partial x^2}} 
\, \Big( 1 - {1\over2} \,\varphi(-1) \,-\, {1\over2} \,\varphi(3) \Big) \,+\, {\rm O}(h^3) \,$

\smallskip \noindent 
and  $\, \rho(h,\,u) \,$  admits the development

\setbox21=\hbox {$\displaystyle   
{{\partial u}\over{\partial t}} \,+\,{1\over2}\, \Delta t  \,  
{{\partial^2 u}\over{\partial t^2}} 
\,-\,{1\over2}\, \sigma^2 \, h \,  {{\partial^2 u}\over{\partial x^2}} \,+\,\,$   }
\setbox22=\hbox {$\displaystyle    
\,+\,  {1\over2} \, h \, \sigma \, ( 1 - \sigma) \, {{\partial^2 u}\over{\partial x^2}} \, 
\Big(1 -   {1\over2} \, \varphi(-1) \, -  {1\over2} \,\, \varphi(3) \Big)  
 \,+\, {\rm O}(h^3) \,. \,$   }
\setbox30= \vbox {\halign{#&# \cr \box21 \cr \box22     \cr   }}
\setbox31= \hbox{ $\vcenter {\box30} $}
\setbox44=\hbox{\noindent  (6.27) $\displaystyle  \quad 
\rho(h,\,u) \,\,=\,\, \left\{ \box31 \right. $}   
\smallskip \noindent   $ \box44 $ 

\smallskip \noindent 
Thus estimate (6.23) is a consequence of (6.25),(6.8) and (6.1). $\hfill \square $

%  \vfill \eject    
%  \smallskip
%%%%%%%%%%%%%%%%%%%%%%%%%%%%%%%%%%%%%%%%%%%%%%%%%%%%%%%%%%%%
%%%%%%%%%   attention    reglages delicats   %%%%%%%%%%%%%%%
%%%%%%%%%%%%%%%%%%%%%%%%%%%%%%%%%%%%%%%%%%%%%%%%%%%%%%%%%%%% 

% \vfill \eject 
\bigskip \noindent  
We can compare on figure 6.1 several profiles obtained for  different meshes at time  $t =
1$  with the Heun scheme. In order to quantify the errors, we set  

\smallskip \noindent   (6.28)  $\quad   \displaystyle 
\mid \mid u_h(t=1) - u(t=1) \mid \mid_{ \displaystyle  \ell^1} \,\,=\,\,  h \, 
\sum_{j=1}^{N(h)} \, \, \mid u_{j,\,h}(t=1) \,-\, u^0(x_j) \, \mid \, $ 

\smallskip \noindent   (6.29)  $\quad   \displaystyle 
\mid \mid u_h(t=1) - u(t=1) \mid \mid_{ \displaystyle  \ell^2} \,\,=\,\,  \bigg[ \,  h \, 
\sum_{j=1}^{N(h)} \, \, \mid u_{j,\,h}(t=1) \,-\, u^0(x_j) \, \mid ^2 \, \bigg]^{1/2} 
\, $ 

\smallskip \noindent   (6.30)  $\quad   \displaystyle 
\mid \mid u_h(t=1) - u(t=1) \mid \mid_{ \displaystyle  \ell^\infty} \,\,=\,\, \max
_{j=1,\, \dots ,\, N(h)} \, \mid u_{j,\,h}(t=1) \,-\, u^0(x_j) \, \mid \, . \, $ 

\vfill \eject    %%%% fd 15 juin 2010 

\centerline { \epsfxsize=11.0 cm    \epsfbox  {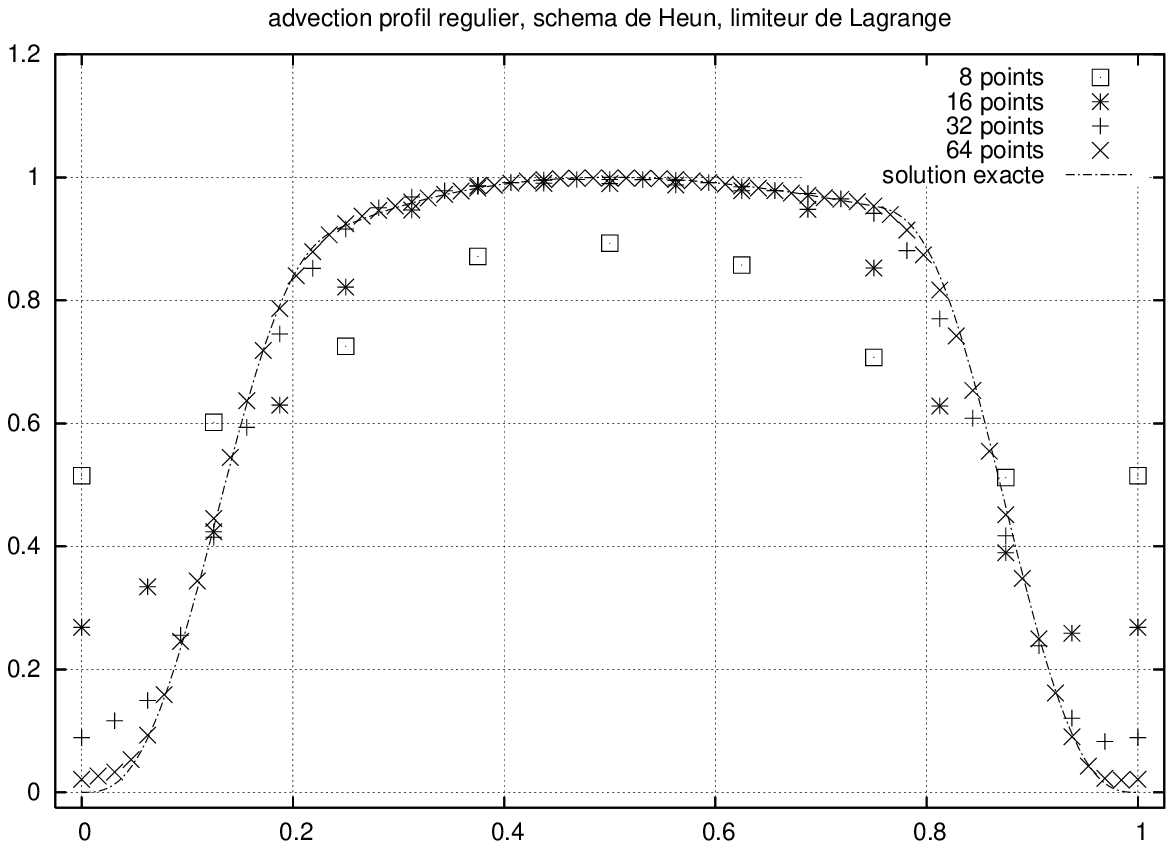} }
\smallskip 
\centerline { \epsfxsize=11.0 cm    \epsfbox  {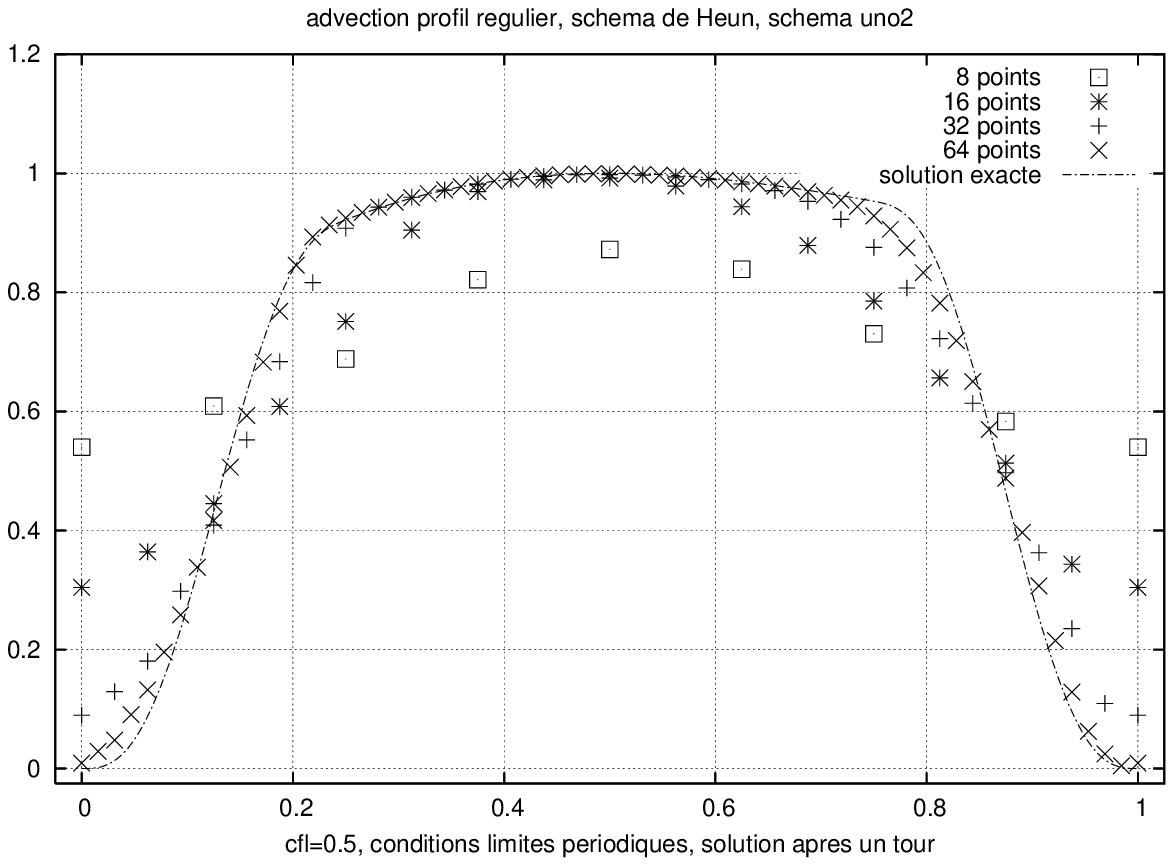} }
\smallskip  \smallskip
\noindent  {\bf Figure 6.1.} \quad  Advection of a regular profile. 
Solutions obtained with the Lagrange limiter and the 5 points uno2 scheme.

\bigskip  
\vfill \eject    %%%% fd 15 juin 2010 
\noindent   $   \bullet $ \qquad
We present four numerical experiments associated respectively to the following choice
concerning the temporal scheme : 

\smallskip  \noindent   exp 1	\quad Five points TVD scheme (6.9)  	\qquad  CFL : $ \,   \sigma  =
{1\over2} \,$

\smallskip \noindent 
 exp 2	\quad Predictor-corrector scheme (6.10)   	\qquad  CFL : $ \,   \sigma  = {1\over2} \,$

\smallskip \noindent 
 exp 3	\quad  Heun scheme (6.11)   \qquad  CFL : $ \,   \sigma  = {1\over2} \,$

\smallskip \noindent 
 exp 4	\quad  Heun scheme (6.11)   \qquad  CFL : $ \,   \sigma  = {1\over4} . \,$

\smallskip  \noindent 
For each of these numerical experiments we have advected the initial profile defined in (6.4)
during one period with seven different meshes defined in (6.6) and the eight different fluxes
(i) to (viii). We present on figures 6.2 to 6.5 the $\, \ell^1 ,\,$   $\, \ell^2 ,\,$ 
 $\, \ell^\infty \,$   errors defined in
(6.28)-(6.30) as functions of $\, {\rm  log}_2 (1/h) \,$   and the truncation error (defined
respectively in (6.22), (5.24) and (5.42)) at the particular point  $\, x = {1\over2} 
\,$ where the maximum occurs (it is always a mesh point). All the residuals are proposed in
terms of their base~2 logarithm. Therefore if some residual $\,  r(h) \,$  is ``at the order
$\alpha$''~:

\smallskip \noindent   (6.31)  $\quad   \displaystyle 
r(h) \,\, =\,\, C \, h^{\alpha} \, \big(1 + {\rm o}(h) \, \big) \,$ 

\smallskip \noindent 
the slope of the curve    $\, {\rm  log}_2 [r(h)] \,$
as a function of $\, {\rm  log}_2 (1/h) \,$  is simply equal to $-\alpha$ and is directly
readable on the figures. We notice first that the truncation errors at the extremum   ($ x =
{1\over2}$)  have exactly the behavior predicted in (6.23), (5.25) and (5.43); with $\,
\sigma = {1\over2} \, $  the order of the truncation error at the extremum is only 1 with the
Heun scheme (Fig. 6.4.d) therefore it is of order 2 with the same numerical scheme and   $\,
\sigma = {1\over4} \, $  (Fig. 6.5.d). The five point scheme (Fig. 6.2) show that the MUSCL
limiter and the UNO2 scheme are really second order accurate in the $\, L^\infty \,$  norm
(Fig. 6.2.c). When we use a two step scheme and a Courant number  $\, \sigma = {1\over2} \, $ 
(Fig. 6.3 and 6.4) the Lagrange limiter is ``better than second order accurate'' (order 2.2 in
the  $\, L^\infty \,$ norm) and the other limiters are second order accurate but ``only'' in
the norms L1 and L2. With the Heun scheme and a severe CFL restriction (Fig. 6.5) the results
of the Lagrange limiter are degraded despite the fact that the truncation error remains
everywhere second order accurate !

\vfill \eject

\noindent
$\!\!\!\!\!\!$
{ \epsfysize=7.5 cm  \epsfbox  {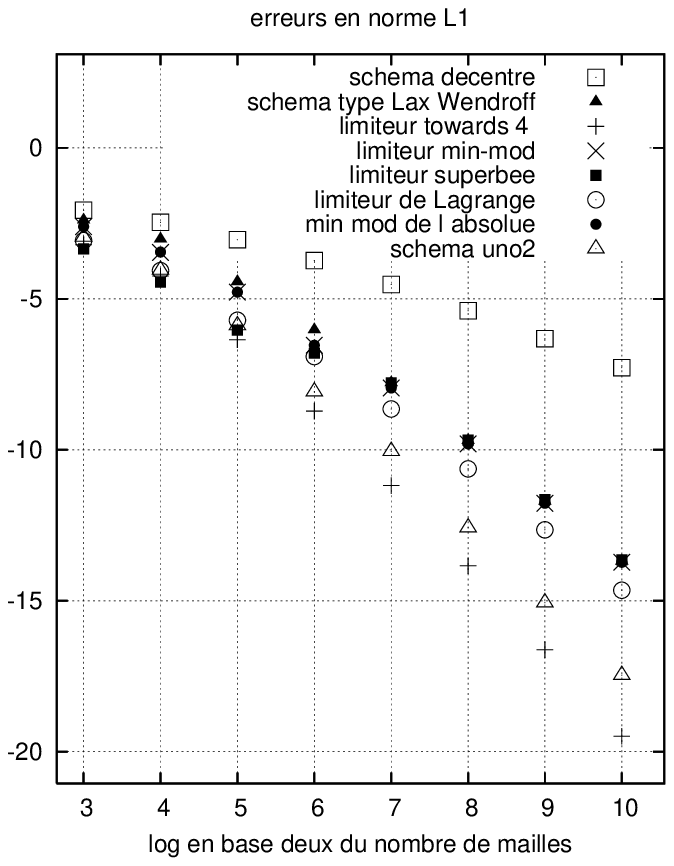}} 
$\!\!\!\!\!\!\!\!\!\!\!\!\!\!\!\!\!\!\!\!\!\!\!\!$
$\!\!\!\!\!\!\!\!\!\!\!\!\!\!\!\!\!\!\!\!\!\!\!\!$  
$\!\!\!\!\!\!\!\!\!\!\!\!\!\!\!\!\!\!\!\!\!\!\!\!$  
$\!\!\!\!\!\!\!\!\!\!$   
{ \epsfysize=7.5 cm  \epsfbox  {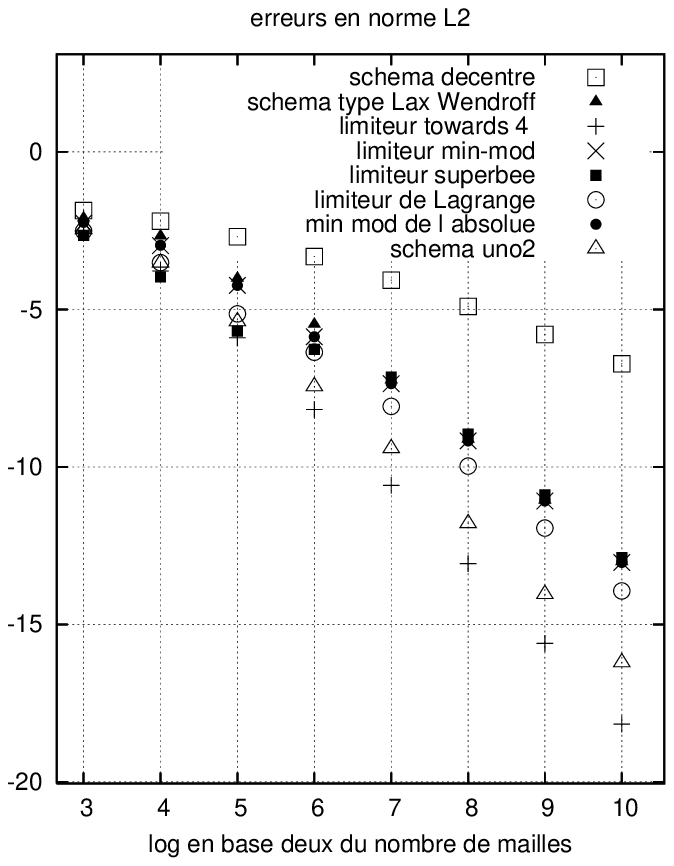}} 
\smallskip  \noindent
$\!\!\!\!\!\!$
{ \epsfysize=7.5 cm  \epsfbox  {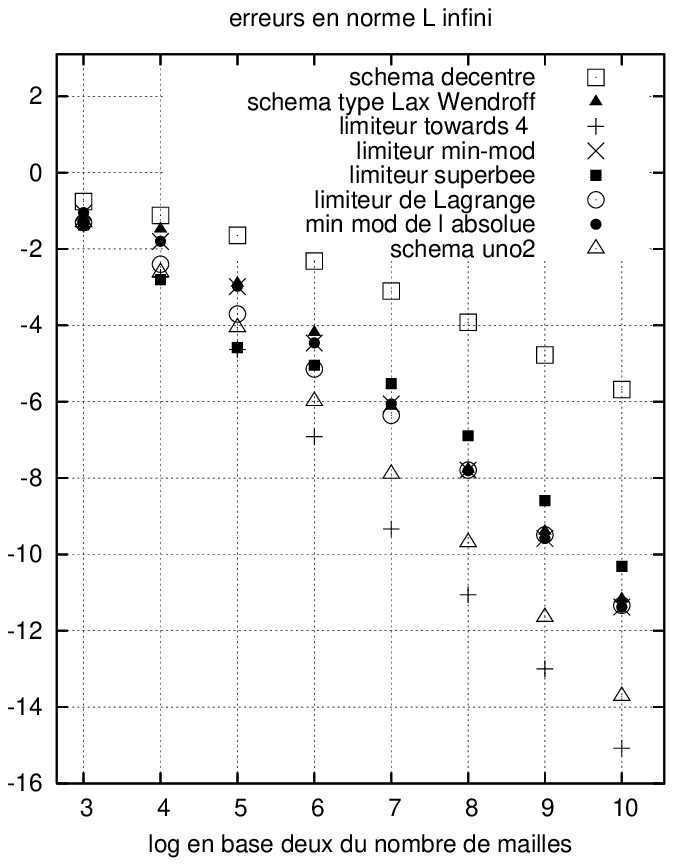}}
$\!\!\!\!\!\!\!\!\!\!\!\!\!\!\!\!\!\!\!\!\!\!\!\!$ 
$\!\!\!\!\!\!\!\!\!\!\!\!\!\!\!\!\!\!\!\!\!\!\!\!$
$\!\!\!\!\!\!\!\!\!\!\!\!\!\!\!\!\!\!\!\!\!\!\!\!$ 
$\!\!\!\!\!\!\!\!\!\!$   
{ \epsfysize=7.5 cm  \epsfbox  {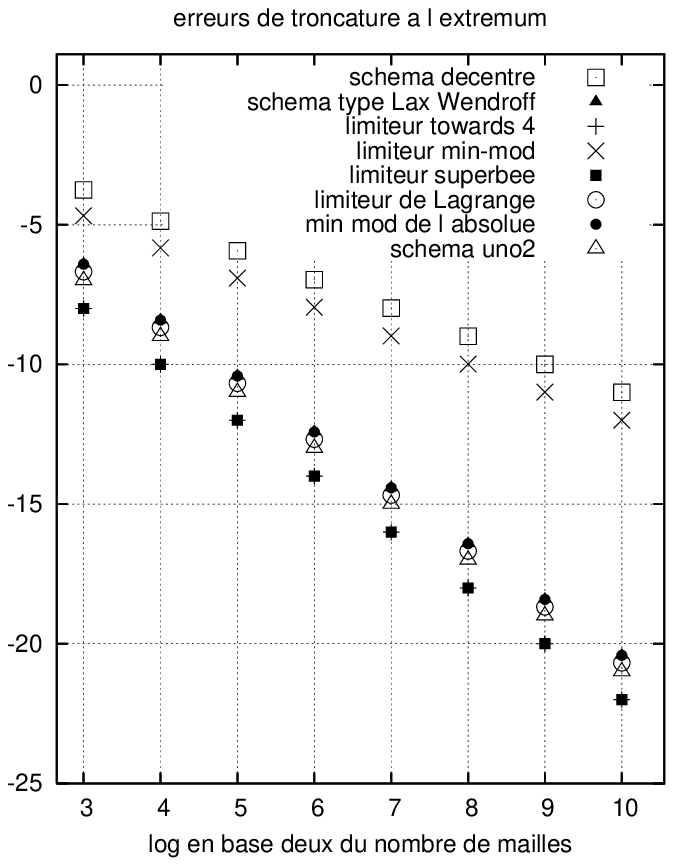}} 
\smallskip  \smallskip
\noindent  {\bf Figure 6.2.} \quad  Advection of a regular profile. 
Errors for the flux corrected transport type temporal scheme with Courant number = 0.5.

\noindent
$\!\!\!\!\!\!$
{ \epsfysize=7.5 cm  \epsfbox  {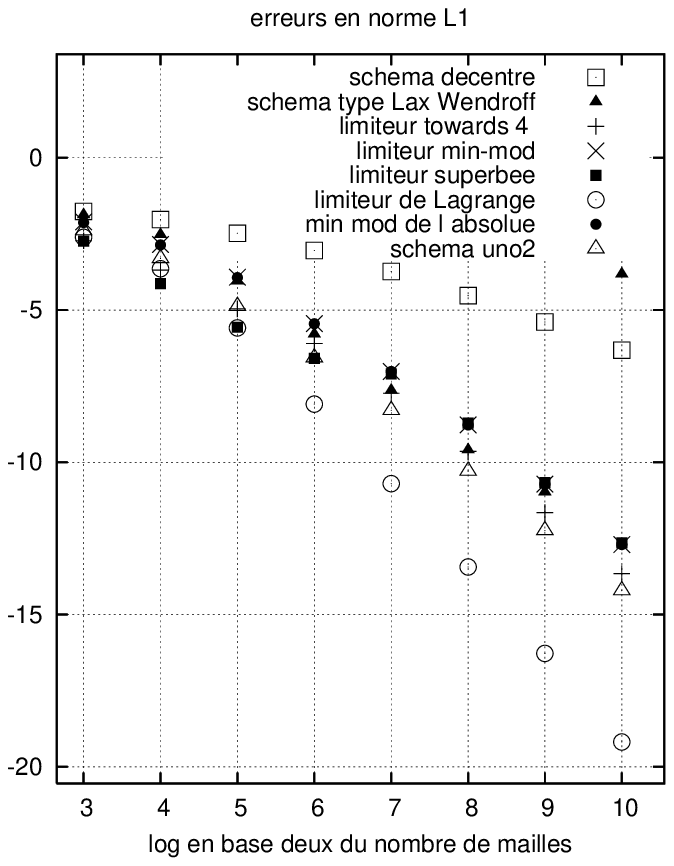}} 
$\!\!\!\!\!\!\!\!\!\!\!\!\!\!\!\!\!\!\!\!\!\!\!\!$
$\!\!\!\!\!\!\!\!\!\!\!\!\!\!\!\!\!\!\!\!\!\!\!\!$  
$\!\!\!\!\!\!\!\!\!\!\!\!\!\!\!\!\!\!\!\!\!\!\!\!$  
$\!\!\!\!\!\!\!\!\!\!$   
{ \epsfysize=7.5 cm  \epsfbox  {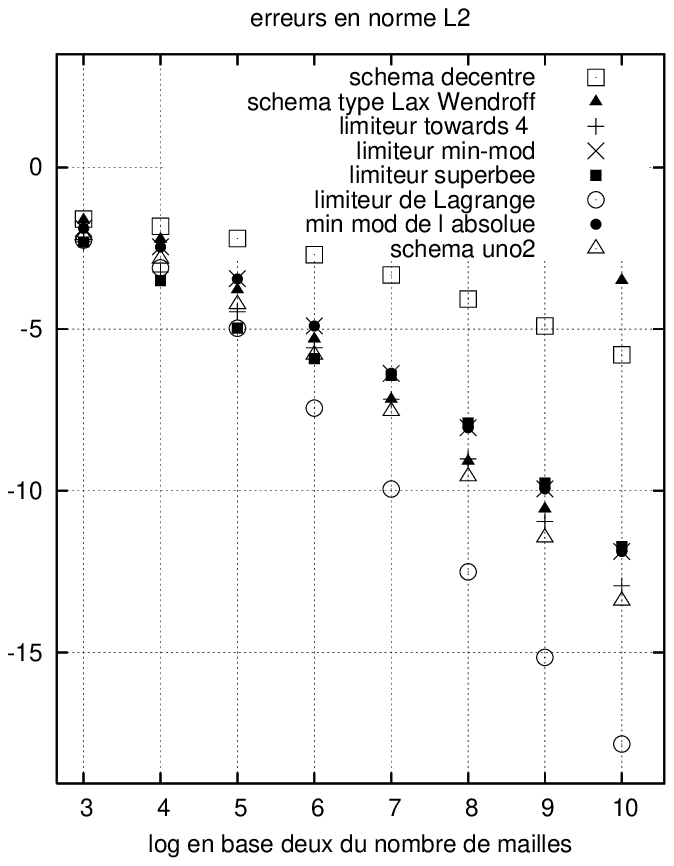}} 
\smallskip  \noindent
$\!\!\!\!\!\!$
{ \epsfysize=7.5 cm  \epsfbox  {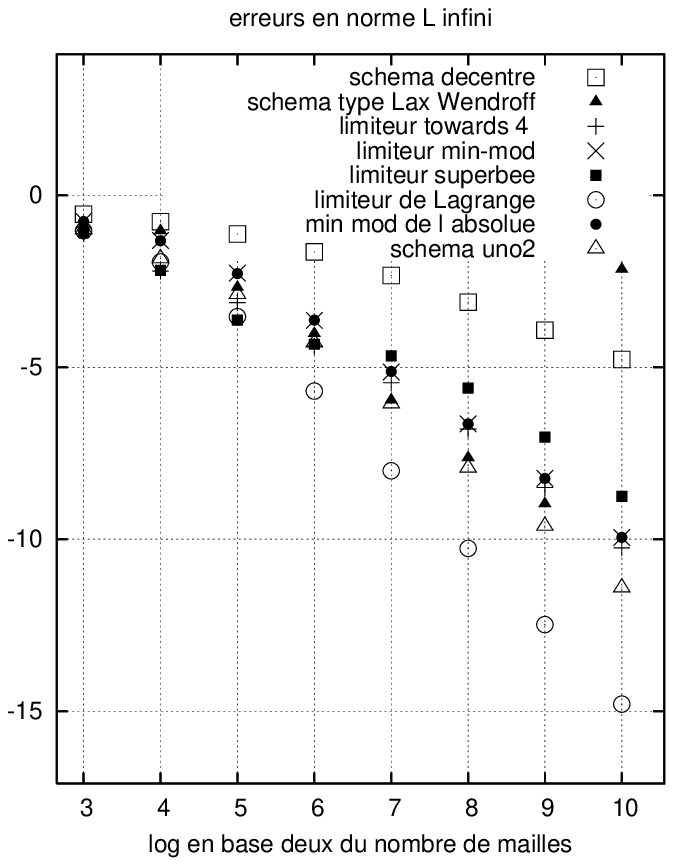}}
$\!\!\!\!\!\!\!\!\!\!\!\!\!\!\!\!\!\!\!\!\!\!\!\!$ 
$\!\!\!\!\!\!\!\!\!\!\!\!\!\!\!\!\!\!\!\!\!\!\!\!$
$\!\!\!\!\!\!\!\!\!\!\!\!\!\!\!\!\!\!\!\!\!\!\!\!$ 
$\!\!\!\!\!\!\!\!\!\!$   
{ \epsfysize=7.5 cm  \epsfbox  {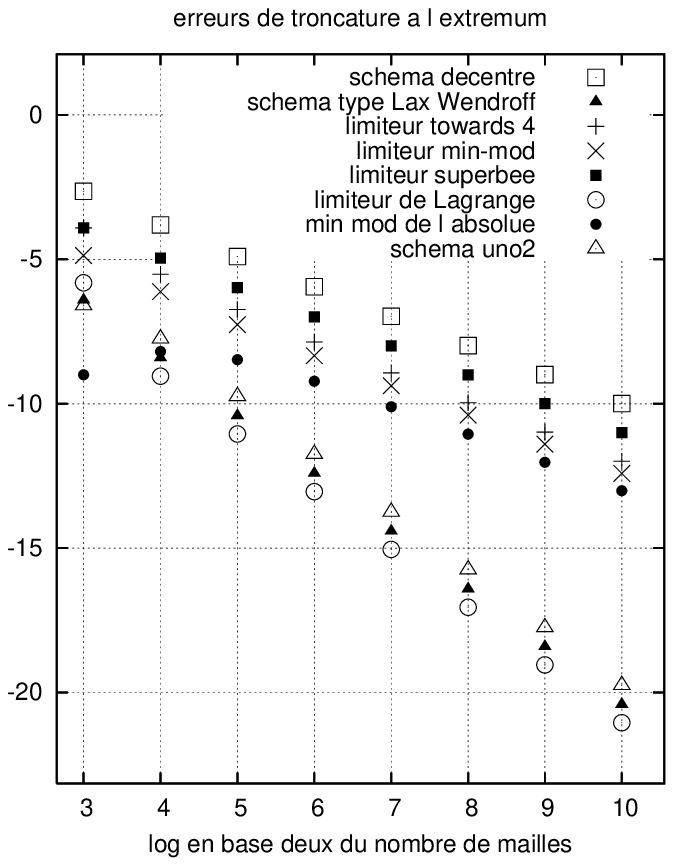}} 
\smallskip  \smallskip
\noindent  {\bf Figure 6.3.} \quad  Advection of a regular profile. 
Errors for the predictor-corrector temporal scheme with Courant number = 0.5.

\noindent
$\!\!\!\!\!\!$
{ \epsfysize=7.5 cm  \epsfbox  {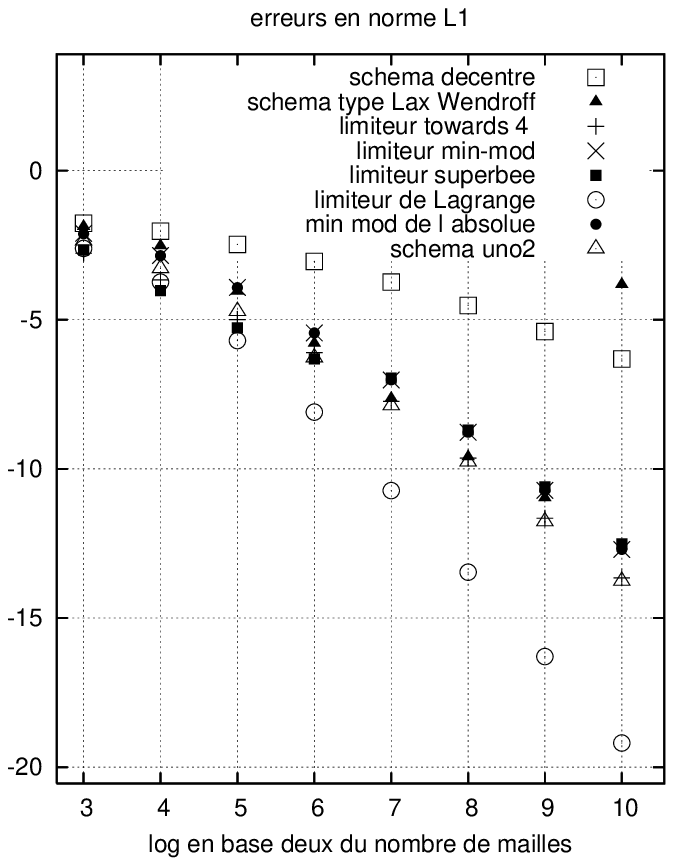}} 
$\!\!\!\!\!\!\!\!\!\!\!\!\!\!\!\!\!\!\!\!\!\!\!\!$
$\!\!\!\!\!\!\!\!\!\!\!\!\!\!\!\!\!\!\!\!\!\!\!\!$  
$\!\!\!\!\!\!\!\!\!\!\!\!\!\!\!\!\!\!\!\!\!\!\!\!$  
$\!\!\!\!\!\!\!\!\!\!$   
{ \epsfysize=7.5 cm  \epsfbox  {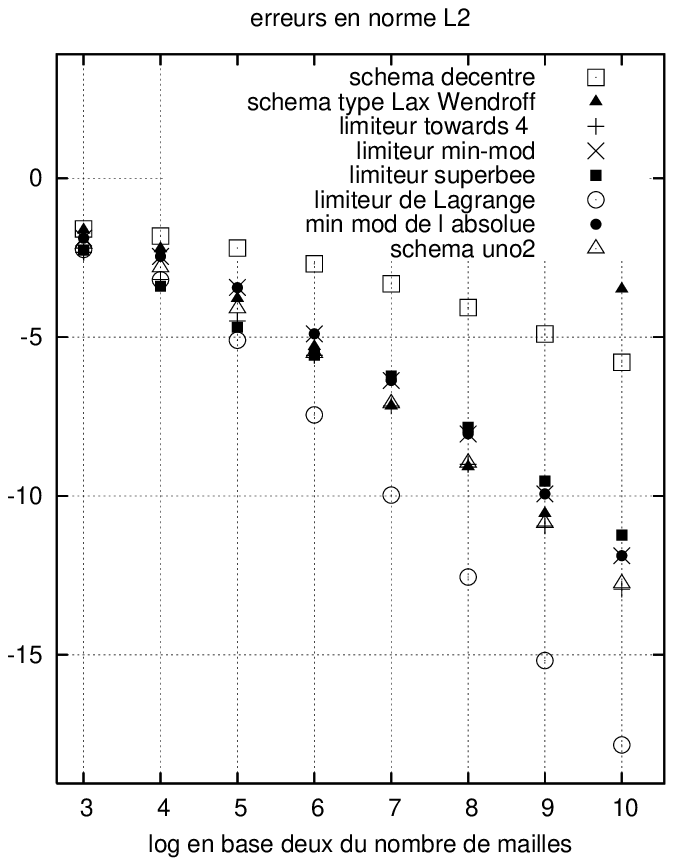}} 
\smallskip  \noindent
$\!\!\!\!\!\!$
{ \epsfysize=7.5 cm  \epsfbox  {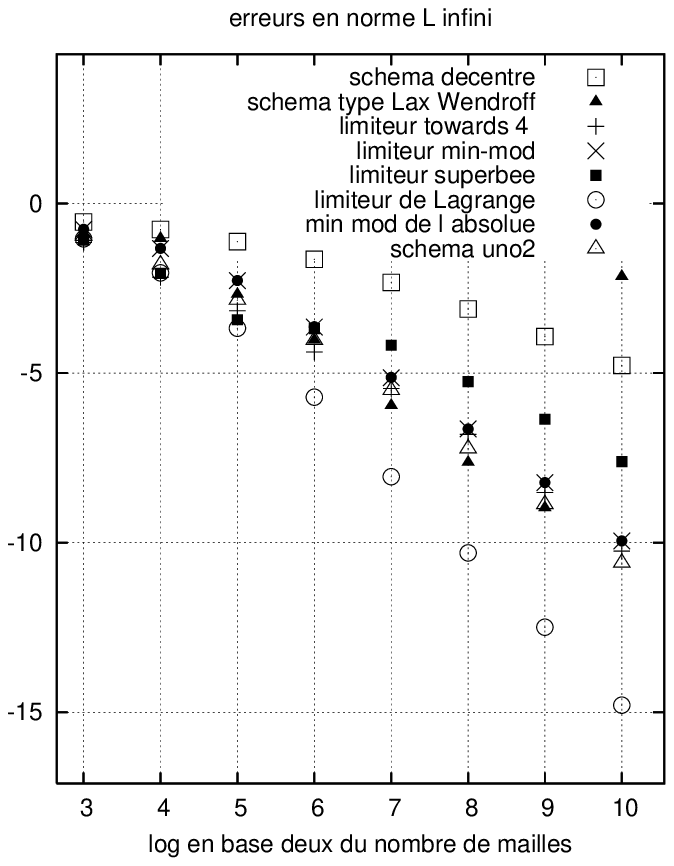}}
$\!\!\!\!\!\!\!\!\!\!\!\!\!\!\!\!\!\!\!\!\!\!\!\!$ 
$\!\!\!\!\!\!\!\!\!\!\!\!\!\!\!\!\!\!\!\!\!\!\!\!$
$\!\!\!\!\!\!\!\!\!\!\!\!\!\!\!\!\!\!\!\!\!\!\!\!$ 
$\!\!\!\!\!\!\!\!\!\!$   
{ \epsfysize=7.5 cm  \epsfbox  {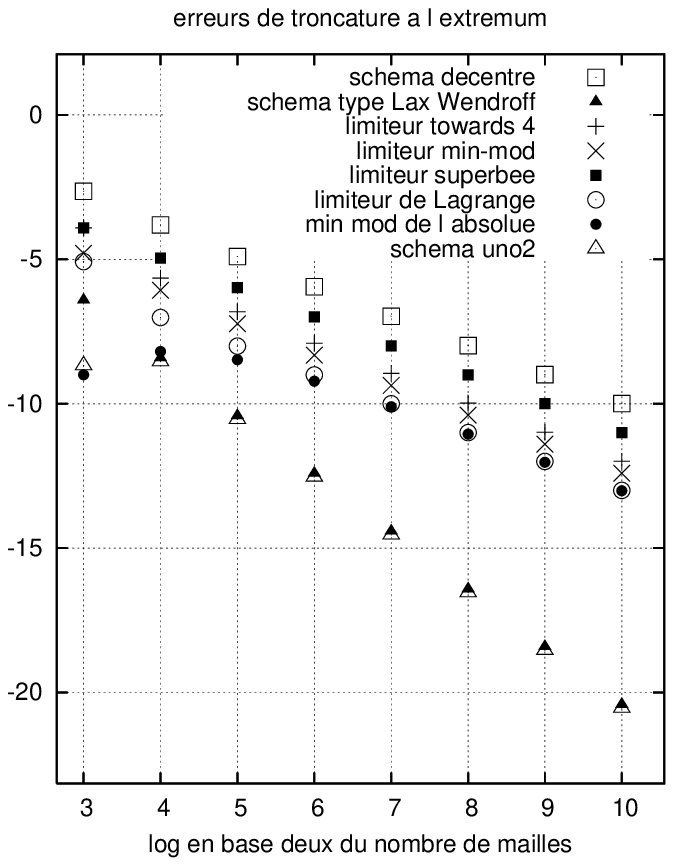}} 
\smallskip  \smallskip
\noindent  {\bf Figure 6.4.} \quad  Advection of a regular profile. 
Errors for the Heun temporal scheme with Courant number = 0.5.

\noindent
$\!\!\!\!\!\!$
{ \epsfysize=7.5 cm  \epsfbox  {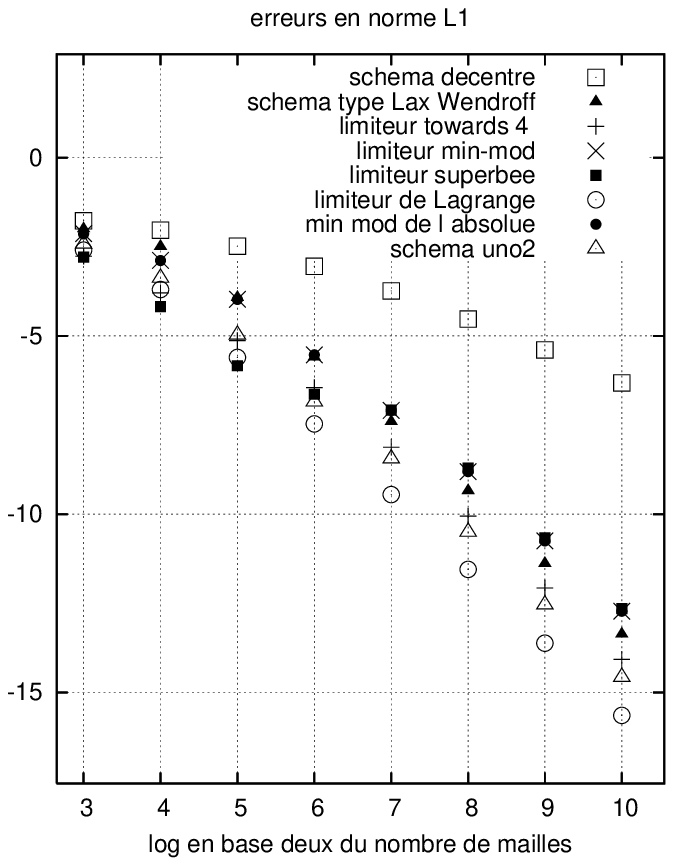}} 
$\!\!\!\!\!\!\!\!\!\!\!\!\!\!\!\!\!\!\!\!\!\!\!\!$
$\!\!\!\!\!\!\!\!\!\!\!\!\!\!\!\!\!\!\!\!\!\!\!\!$  
$\!\!\!\!\!\!\!\!\!\!\!\!\!\!\!\!\!\!\!\!\!\!\!\!$  
$\!\!\!\!\!\!\!\!\!\!$   
{ \epsfysize=7.5 cm  \epsfbox  {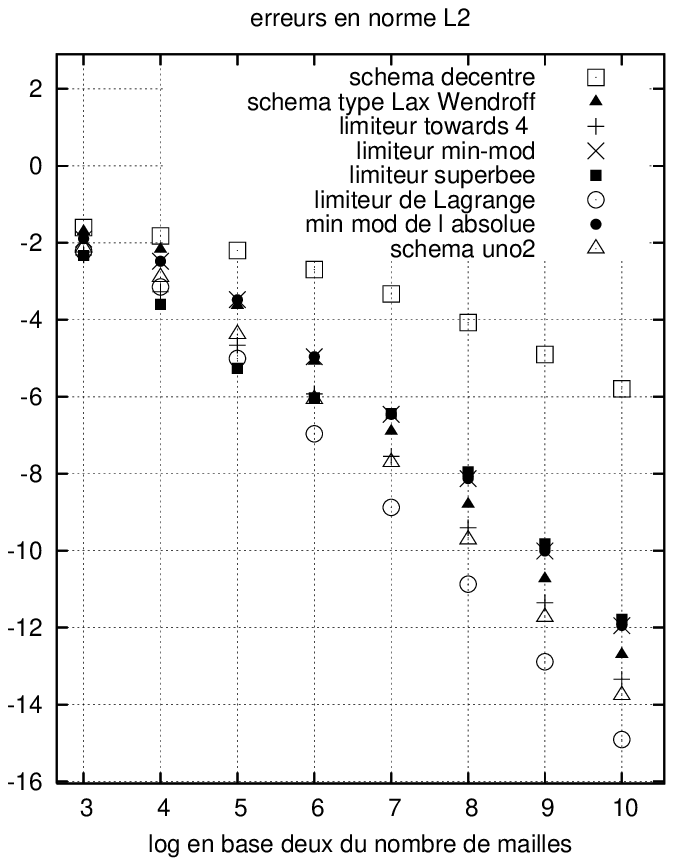}} 
\smallskip  \noindent
$\!\!\!\!\!\!$
{ \epsfysize=7.5 cm  \epsfbox  {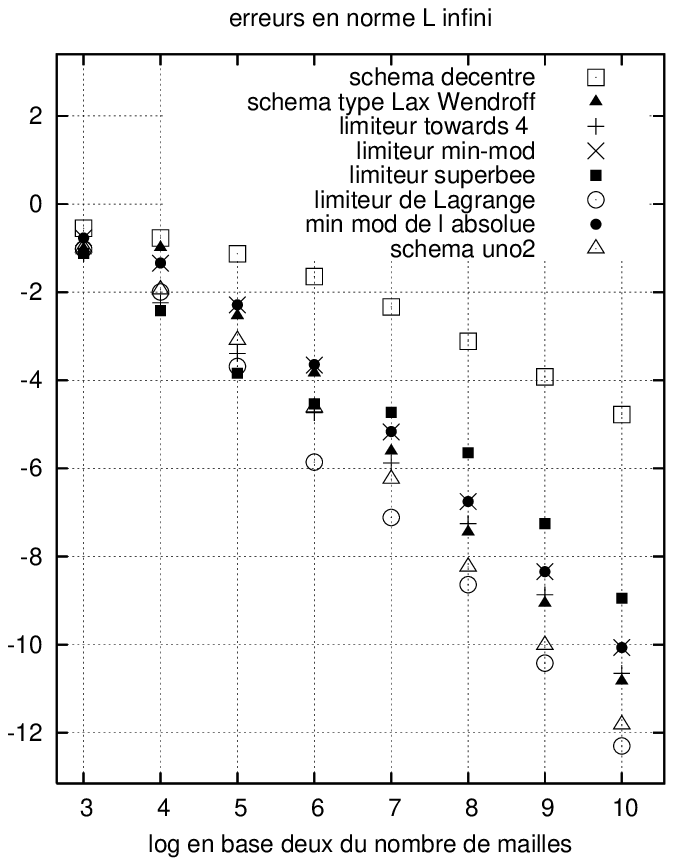}}
$\!\!\!\!\!\!\!\!\!\!\!\!\!\!\!\!\!\!\!\!\!\!\!\!$ 
$\!\!\!\!\!\!\!\!\!\!\!\!\!\!\!\!\!\!\!\!\!\!\!\!$
$\!\!\!\!\!\!\!\!\!\!\!\!\!\!\!\!\!\!\!\!\!\!\!\!$ 
$\!\!\!\!\!\!\!\!\!\!$   
{ \epsfysize=7.5 cm  \epsfbox  {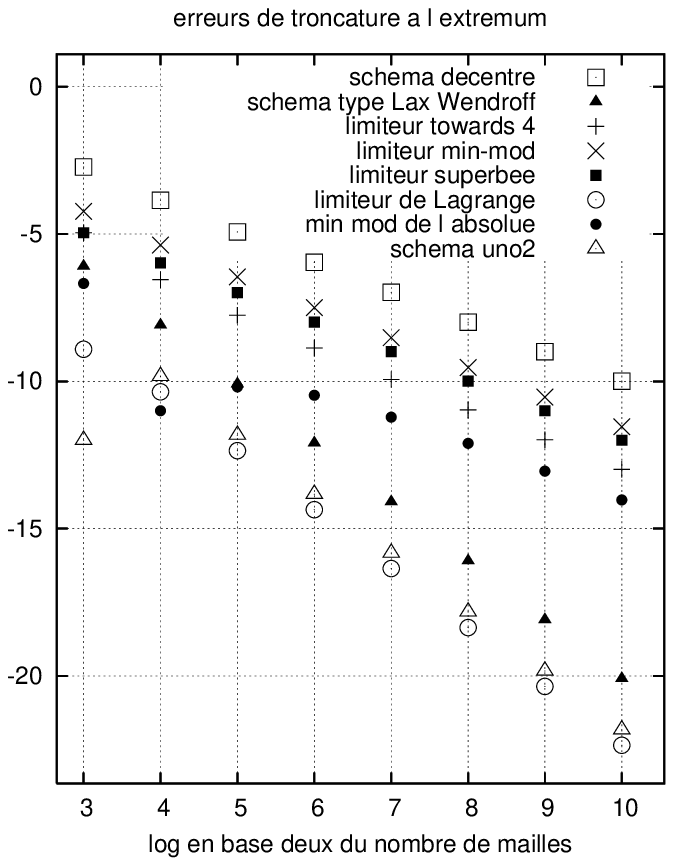}} 
\smallskip  \smallskip
\noindent  {\bf Figure 6.5.} \quad  Advection of a regular profile. 
Errors for the Heun temporal scheme with Courant number = 0.25.

%%%%%%%%%%%%%%%%%%%%%%%%%%%%%%%%%%%%%%%%%%%%%%%%%%%%%%%%%%%%%%%%%%%%%%%%%%%%%%%%%%%%%%%%%%%%%%%%
%                       \titredroite={\pecaps  Conclusion }
\toppagetrue  
\botpagetrue    

\bigskip  %  \bigskip  
\vfill \eject   %%% fd 15 juin 2010 
\noindent {\bf 7. \quad Conclusion}

\noindent
In this contribution, we have proposed a concept of convexity preserving  
property for the approximation of scalar conservation laws with the Van Leer's numerical
method. New restrictions have been derived for the construction of the associated limiter
functions. We have proved that the method of lines is well posed for $\, u_0 \,$ in $\,
\ell^1 \cap {\rm  BV} \,$   and new conditions on the limiter have been given to
prove the TVD property of the associated solution. We have proved that second order accuracy
(in the sense of the truncation error) can be maintained, even at a nonsonic extremum ; this
last property and Total Variation Diminishing are compatible with a discretization in time
with a two-level Runge-Kutta scheme under new restrictions on the limiter function that have
been detailed. An example of such a limiter (the so-called Lagrange limiter) has been
proposed and implemented in the case of the advection equation. Theoretical results on orders
of convergence have been confirmed numerically.
This concept of convexity preserving interpolation could be extended 
in two directions~: the treatment of multidimensional scalar problems and the approximation
of systems of conservation laws.

\toppagetrue  
\botpagetrue    

\bigskip %   \bigskip  
\noindent {\bf 8. \quad References}

% \vskip 0.07 cm  
 \hangindent=7mm \hangafter=1 \noindent 
Anderson W.K., Thomas J.L., Van Leer B.	\quad 
Comparison of Finite Volume Flux Vector Splittings for the Euler Equations, AIAA J., {\bf 24},
p.~1453-1460, 1986.

\vskip 0.065 cm   \hangindent=7mm \hangafter=1 \noindent 
Boris J.P., Book D.L. \quad 
Flux Corrected Transport I. SHASTA, A Fluid Transport Algorithm  That Works, J. Comp. Physics,
{\bf  11}, p.~38-69, 1973.

\vskip 0.065 cm   \hangindent=7mm \hangafter=1 \noindent 
Cahouet J., Coquel F. \quad 
Uniformly Second Order Convergent Schemes for 
Hyperbolic Conservation Laws including Leonard's Approach, in {\bf  Notes on Numer. Fluid
Mech., 24}, (Ballmann, Jeltsch eds), Vieweg Verlag, p.~51-62, 1989.

\vskip 0.065 cm   \hangindent=7mm \hangafter=1 \noindent 
Chakravarthy S.R., Szema K.Y., Goldberg U.C., Gorski J.J., Osher S. \quad 
Application of a New Class of High Accuracy TVD Schemes to the Navier Stokes Equations, AIAA
paper 85-0165, 1985.

\vskip 0.065 cm   \hangindent=7mm \hangafter=1 \noindent 
Colella  P. \quad 
A Direct Eulerian MUSCL Scheme for Gas Dynamics, SIAM J. Sci. Stat. Comput., {\bf 6}, pp.
104-117, 1985.

\vskip 0.065 cm   \hangindent=7mm \hangafter=1 \noindent 
Colella P, Woodward P.R. \quad 
The Piecewise Parabolic Method (PPM) for the Gas Dynamical Simulations, J. Comp. Physics, 54,
p.~174-201, 1984.

\vskip 0.065 cm   \hangindent=7mm \hangafter=1 \noindent 
Crandall M.G., Majda A.	\quad 
Monotone Difference Approximations for Scalar Conservation Laws, Math. of Comp., {\bf 34}, pp.
1-21, 1980.

\vskip 0.065 cm   \hangindent=7mm \hangafter=1 \noindent 
Di Perna R. \quad 
Convergence of Approximate Solutions to Conservation Laws, Arch. Rational Mech. Anal., {\bf
82}, p.~27-70, 1983.

\vskip 0.065 cm   \hangindent=7mm \hangafter=1 \noindent 
Dubois F. \quad 
Interpolation non lin\'eaire et sch\'emas \`a variation totale d\'ecroissante, internal
report, A\'erospatiale Les Mureaux ST/MI n$^{\rm o}$40 348, november 1988.

\vskip 0.065 cm   \hangindent=7mm \hangafter=1 \noindent 
Engquist B., Osher S. \quad 
Stable and Entropy Condition Satisfying Approximations for Transonic Flow Calculations, Math.
of Comp., {\bf 34}, p.~45-75, 1980.

\vskip 0.065 cm   \hangindent=7mm \hangafter=1 \noindent 
Godunov S.K. \quad 
A finite Difference Method for the Numerical Computation of Discontinuous Solutions of the
Equations of Fluid Dynamics, Math. Sb., {\bf 47}, p.~271-290, 1959.

\vskip 0.065 cm   \hangindent=7mm \hangafter=1 \noindent 
Harten A. \quad 
High Resolution Schemes for Hyperbolic Conservation Laws, J. Comp. Physics, {\bf 49}, pp.
357-393, 1983.

\vskip 0.065 cm   \hangindent=7mm \hangafter=1 \noindent 
Harten A., Engquist B., Osher S., Chakravarthy S. \quad 
High Order Accurate Essentielly Non-Oscillatory Schemes III, J. Comp. Physics, {\bf  71}, pp.
231-303, 1987.

\vskip 0.065 cm  \hangindent=7mm \hangafter=1 \noindent 
Harten A., Osher S. \quad 
Uniformly High Order Accurate Nonoscillatory Schemes I, SIAM J. Numer. Anal., {\bf 24}, pp.
279-309, 1987.

\vskip 0.065 cm  \hangindent=7mm \hangafter=1 \noindent 
Harten A., Hyman J.M., Lax P.D. quad 
On Finite Difference Approximations and Entropy Conditions For Shocks, Comm. Pure and Applied
Maths., {\bf 29}, p.~297-322, 1976.

\vskip 0.065 cm  \hangindent=7mm \hangafter=1 \noindent 
Harten A., Lax P.D., Van Leer B. \quad 
On Upstream Differencing and Godunov-type Schemes for Hyperbolic Conservation Laws, SIAM
Review, {\bf 25}, p.~35-61, 1983.

\vskip 0.065 cm  \hangindent=7mm \hangafter=1 \noindent 
Kruskov S. \quad 
First Order Quasilinear Equations with Several Space Variables, Math. Sb., {\bf 123}, pp.
228-255, 1970.

\vskip 0.065 cm  \hangindent=7mm \hangafter=1 \noindent 
Kusnezov N.N., Volosin  S.A. \quad 
On Monotone Difference Approximations for a First-order Quasi-linear Equation, Soviet Math.
Dokl., {\bf 17}, pp 1203-1206, 1976.

\vskip 0.065 cm  \hangindent=7mm \hangafter=1 \noindent 
Leonard B.P. \quad 
A Stable and Accurate Convective Modelling Procedure Based on Quadratic Upstream
Interpolation, Comp. Meth. in Appl. Mech. and Eng., {\bf 19}, p.~59-98, 1979.

\vskip 0.065 cm  \hangindent=7mm \hangafter=1 \noindent 
Leroux A.Y. \quad 
Convergence of the Godunov Scheme for First Order Quasilinear Equations, Proc. Japan Acad.,
{\bf 52}, p.~488-491, 1976.

\vskip 0.065 cm  \hangindent=7mm \hangafter=1 \noindent 
Lax P.D. \quad 
Shock Waves and Entropy, in {\bf Contributions to Nonlinear Functional Analysis}, edited by
Zarantonello, Academic Press, New York, 1971, p.~603-634.

\vskip 0.065 cm  \hangindent=7mm \hangafter=1 \noindent 
Lax P.D., Wendroff B. \quad 
Systems of Conservation Laws, Comm. Pure and Applied Maths., {\bf  13}, p.~217-237, 1960.

\vskip 0.065 cm  \hangindent=7mm \hangafter=1 \noindent 
Osher S.  \quad 
Riemann Solvers, the Entropy Condition, and Difference Approximations, SIAM J. Numer. Anal.,
{\bf 21}, p.~217-235, 1984.

\vskip 0.065 cm  \hangindent=7mm \hangafter=1 \noindent 
Osher S. \quad 
Convergence of generalized MUSCL schemes , SIAM J. of Numer. Anal., {\bf 22}, p.~947-961,
1985.

\vskip 0.065 cm  \hangindent=7mm \hangafter=1 \noindent 
Osher S., Chakravarthy S. \quad 
High Resolution Schemes and the Entropy Condition, SIAM J. Numer Anal., {\bf 21}, p.~955-984,
1984.

\vskip 0.065 cm  \hangindent=7mm \hangafter=1 \noindent 
Osher S., Tadmor E. \quad 
On the Convergence of Difference Approximations to Scalar Conservation Laws, Math. of Comp.,
{\bf 50}, n$^{\rm o}$ 181, p.~19-51, 1988.

\vskip 0.065 cm  \hangindent=7mm \hangafter=1 \noindent 
Roe P. \quad 
Some Contributions to the Modelling of Discontinuous Flows, in {\bf  Lectures in Applied
Maths.}, vol 22 (Engquist, Osher, Sommerville Eds), AMS, p.~161-193, 1985.

\vskip 0.065 cm  \hangindent=7mm \hangafter=1 \noindent 
Sanders R.  \quad 
On Convergence of Monotone Finite Difference Schemes With Variable Spatial Differencing, Math.
of Comp., {\bf 40}, p.~91-106, 1983.  

\vskip 0.065 cm   \hangindent=7mm \hangafter=1 \noindent 
Shu C.W., Osher S. \quad 
Efficient Implementation of Essentially Non-Oscillatory Shock-Capturing Schemes, J. Comp.
Physics, {\bf 77}, n$^{\rm o}$ 2, p.~439-471, 1988.

\vskip 0.065 cm   \hangindent=7mm \hangafter=1 \noindent 
Smoller J.  \quad 
{\bf Shock Waves and Reaction Diffusion Equations}, Springer Verlag, New York, 1983.  
 
\vskip 0.065 cm   \hangindent=7mm \hangafter=1 \noindent 
 Spekreijse S.  \quad 
Multigrid Solution of Monotone Second-Order Discretizations of Hyperbolic Conservation Laws,
Math. of Comp. {\bf 49}, p.~135-155, 1987.

\vskip 0.065 cm   \hangindent=7mm \hangafter=1 \noindent 
Sweby  P.K. \quad 
High Resolution Schemes Using Flux Limiters for Hyperbolic Conservation Laws, SIAM J. of
Numer. Anal., {\bf 21}, p.~995-1011, 1984.

\vskip 0.065 cm   \hangindent=7mm \hangafter=1 \noindent 
Van Albada  G.D., Van Leer B., Roberts W.W. \quad 
A Comparative Study of Computational Methods in Cosmic Gas Dynamics, Astronomy and
Astrophysics, {\bf 1065}, p.~76-84, 1982.

\vskip 0.065 cm   \hangindent=7mm \hangafter=1 \noindent 
Van Leer B. \quad 
Towards the Ultimate Conservative Scheme I. The Quest of Monotonicity, in Lectures Notes in
Physics, {\bf 18}, Springer Verlag, New York, p.~163-168, 1973.

\vskip 0.065 cm  \hangindent=7mm \hangafter=1 \noindent 
Van Leer B.  \quad 
Towards the Ultimate Conservative Scheme IV. A New Approach to Numerical Convection, J. Comp.
Phys., {\bf 23}, p.~276-299, 1977.

\vskip 0.065 cm   \hangindent=7mm \hangafter=1 \noindent 
Van Leer B. \quad 
Towards the Ultimate Conservative Scheme V. A Second Order Sequel to Godunov's Method, J.
Comp. Phys., {\bf 32}, p.~101-136, 1979.

\vskip 0.065 cm   \hangindent=7mm \hangafter=1 \noindent 
Vila J.P. \quad 
High Order Scheme and Entropy Condition for Nonlinear Hyperbolic Systems of Conservation Laws,
Math. of Comp., {\bf 50}, p.~53-73, 1988.

\vskip 0.065 cm   \hangindent=7mm \hangafter=1 \noindent 
Volpert A.I. \quad 
The spaces BV and quasilinear equations, Math. Sb. {\bf 73} (115), p.~255-302, 1967 and Math.
USSR Sb. {\bf 2}, p.~225-267, 1967. 
 
\vskip 0.065 cm   \hangindent=7mm \hangafter=1 \noindent 
Wu H. \quad 
MmB - A New Class of Accurate High Resolution Schemes for Conservation Laws in Two Dimensions,
8th GAMM Conference on Numerical Methods in Fluid Mechanics, Delft, 27-29 sept 1989.

\bye